\documentclass[onefignum,onetabnum]{siamart190516}

\usepackage{lipsum}
\usepackage{amsfonts}
\usepackage{graphicx}
\usepackage{epstopdf}
\usepackage{siunitx}

\usepackage{algorithmic}

\usepackage[letterpaper, margin=1.in]{geometry}
\usepackage{amsmath,amsbsy,amsfonts,amssymb}
\usepackage{graphicx}
\usepackage{xcolor}
\usepackage{mathtools}
\graphicspath{{freegs/figs/}{figs/}}
\usepackage{subfiles}
\usepackage{blindtext}

\usepackage{hyperref}
\hypersetup{bookmarks=true,colorlinks=true, linkcolor=blue,urlcolor=red,linktoc=all}
\usepackage{listings}
\usepackage{soul}

\usepackage{titlesec}
\usepackage{float}
\usepackage{esvect}
\usepackage{media9}
\usepackage{mwe}
\titlespacing*{\section}{0pt}{1.2\baselineskip}{\baselineskip}
\titlespacing*{\subsection}{0pt}{1.0\baselineskip}{0.9\baselineskip}
\titlespacing*{\subsubsection}{0pt}{1.0\baselineskip}{0.9\baselineskip}
\newcommand{\vect}[1]{\textbf{#1}}

\usepackage[utf8x]{inputenc}
\usepackage{amsmath}
\usepackage{comment}
\usepackage{color}
\graphicspath{{figs/}}
\usepackage{caption}
\usepackage{subcaption}
\usepackage{tikz}
\usetikzlibrary{arrows,calc,decorations.markings}
\usetikzlibrary{shapes,arrows}
\usepackage{pgfplots}
\usetikzlibrary{decorations.pathreplacing}
\usetikzlibrary{shapes,arrows}
\definecolor{MyLightMagenta}{cmyk}{0.1,0.8,0,0.1}
\definecolor{Grey}{cmyk}{0,0,0,0.1}

\ifpdf
  \DeclareGraphicsExtensions{.eps,.pdf,.png,.jpg}
\else
  \DeclareGraphicsExtensions{.eps}
\fi

\newsiamremark{remark}{Remark}
\newsiamremark{hypothesis}{Hypothesis}
\crefname{hypothesis}{Hypothesis}{Hypotheses}
\newsiamthm{claim}{Claim}

\headers{Parallel cut-cell algorithm for free-boundary Grad-Shafranov}{Shuang Liu, Qi Tang, and Xian-zhu Tang}

\title{A parallel cut-cell algorithm for the free-boundary Grad-Shafranov problem\thanks{Submitted to the editors \today.
\funding{{This work was supported by the U.S. Department of Energy through the Fusion
Theory Program of the Office of Fusion Energy Sciences, and the Tokamak Disruption Simulation
(TDS) SciDAC partnerships between the Office of Fusion Energy Science and the Office of Advanced
Scientific Computing.}}}}

\author{Shuang Liu\thanks{Los Alamos National Laboratory, Los Alamos, NM 87545, USA 
  (\email{shuangliu@lanl.gov}). Current address: Department of Mathematics, University of California, San Diego,  La Jolla, CA 92093, USA (\email{shl083@ucsd.edu}).}
\and Qi Tang\thanks{Los Alamos National Laboratory, Los Alamos, NM 87545, USA
  (\email{qtang@lanl.gov}).}
\and Xian-zhu Tang\thanks{Los Alamos National Laboratory, Los Alamos, NM 87545, USA (\email{xtang@lanl.gov}).}}

\usepackage{amsopn}

\DeclareMathOperator*{\argmin}{arg\,min}

\ifpdf
\hypersetup{
 pdftitle={Parallel cut-cell algorithm for free-boundary Grad-Shafranov},
 pdfauthor={Shuang Liu, Qi Tang, and Xian-zhu Tang}
}
\fi

\begin{document}

\maketitle

\begin{abstract}
A parallel cut-cell algorithm is described to solve the free-boundary problem of the Grad-Shafranov equation.
The algorithm reformulates the free-boundary problem in an irregular bounded domain and  its important aspects include
a searching algorithm for the magnetic axis and separatrix, 
a surface integral along the irregular boundary to determine the boundary values, 
an approach to optimize the coil current based on a targeting plasma shape,
Picard iterations with Aitken's acceleration for the resulting nonlinear problem,  and a Cartesian grid embedded boundary method to handle the complex geometry.
The algorithm is implemented in parallel using a standard domain-decomposition approach and a good parallel scaling is observed.
Numerical results verify the accuracy and efficiency of the free-boundary Grad-Shafranov solver. 
\end{abstract}



\section{Introduction}

Tokamak fusion relies on magnetic confinement of a plasma at   a
temperature of around 10~keV  and a particle density of
10$^{20}$/m$^3.$ The force balance is achieved by running a current
inside the plasma that produces a Lorentz force to counter the plasma
pressure gradient. In an axisymmetric configuration like a Tokamak,
such an equilibrium is described by an elliptic equation for the
poloidal magnetic flux, commonly known as the Grad-Shafranov
equation. In addition to the plasma current, electrical currents are
also carried in toroidal and poloidal magnetic field coils outside the
plasma chamber to produce the confining magnetic field.  The coil
current can be individually adjusted to accommodate a variety of
plasma pressure and current profiles in terms of plasma positioning
and shaping.

In a fixed-boundary Grad-Shafranov solver, both the location of the
computational boundary and the boundary condition for the poloidal
magnetic flux are known.  Many numerical approaches, including
spectral elements~\cite{li2020solving,howell2014solving}, hybridizable
discontinuous Galerkin methods~\cite{sanchez2019hybridizable,
  sanchez2020adaptive}, boundary integral
approaches~\cite{pataki2013fast, lee2015ecom}, Hermite finite
element~\cite{lutjens1996chease,huysmans1991isoparametric}{,}  and
discontinuous Petrov Galerkin methods~\cite{dpgGS}, have been extended
to solve the fixed-boundary problem.  A common use of the
fixed-boundary Grad-Shafranov solver is to set the computational
boundary to be the targeted last closed flux surface, so the plasma
shaping is enforced by the computational boundary on which the
poloidal flux is a constant.  The free-boundary Grad-Shafranov solver
aims to account for the coil currents as well as the plasma current on
plasma positioning and shaping.  It is an essential tool for tokamak
machine design and specific experimental discharge planning.  Most
tokamak facilities around the world have in-house development and/or
maintenance of free-boundary Grad-Shafranov solvers. There have also
been more recent interest that have a focus on the numerical aspects
of the free-boundary problem~\cite{honda2010simulation,
  faugeras2017fem, heumann2017finite}.

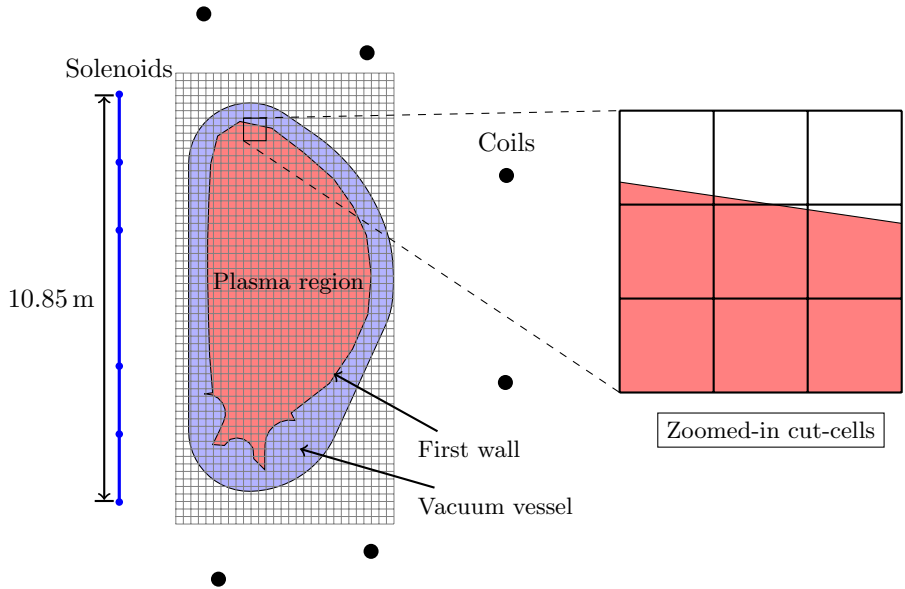
\begin{figure}[hbt]
\centering
\begin{tikzpicture}[scale=.5] 
 \useasboundingbox (1,-8) rectangle (22,8);  
 
\coordinate (rz-1) at (6.2670,   -3.0460); 
\coordinate (rz-2) at (7.2830,   -2.2570); 
\coordinate (rz-3) at (7.8990,   -1.3420); 
\coordinate (rz-4) at (8.3060,   -0.4210); 
\coordinate (rz-5) at (8.3950,    0.6330);  
\coordinate (rz-6) at (8.2700,    1.6810);  
\coordinate (rz-7) at (7.9040,    2.4640);
\coordinate (rz-8) at (7.4000,    3.1790);  
\coordinate (rz-9) at (6.5870,    3.8940);  
\coordinate (rz-10) at (5.7530,    4.5320);  
\coordinate (rz-11) at (4.9040,    4.7120);  
\coordinate (rz-12) at (4.3110,    4.3240);  
\coordinate (rz-13) at (4.1260,    3.5820);  
\coordinate (rz-14) at (4.0760,    2.5660);
\coordinate (rz-15) at (4.0460,    1.5490);  
\coordinate (rz-16) at (4.0460,    0.5330);  
\coordinate (rz-17) at (4.0670,   -0.4840); 
\coordinate (rz-18) at (4.0970,   -1.5000); 
\coordinate (rz-19) at (4.1780,   -2.5060); 
\coordinate (rz-20) at (3.9579,   -2.5384); 
\coordinate (rz-21) at (4.0034,   -2.5384);
\coordinate (rz-22) at (4.1742,   -2.5674); 
\coordinate (rz-23) at (4.3257,   -2.6514); 
\coordinate (rz-24) at (4.4408,   -2.7808); 
\coordinate (rz-25) at (4.5066,   -2.9410); 
\coordinate (rz-26) at (4.5157,   -3.1139); 
\coordinate (rz-27) at (4.4670,   -3.2801); 
\coordinate (rz-28) at (4.4064,   -3.4043);
\coordinate (rz-29) at (4.4062,   -3.4048); 
\coordinate (rz-30) at (4.3773,   -3.4799); 
\coordinate (rz-31) at (4.3115,   -3.6148); 
\coordinate (rz-32) at (4.2457,   -3.7497); 
\coordinate (rz-33) at (4.1799,   -3.8847); 
\coordinate (rz-34) at (4.4918,   -3.9092); 
\coordinate (rz-35) at (4.5687,   -3.8276);
\coordinate (rz-36) at (4.6456,   -3.7460); 
\coordinate (rz-37) at (4.8215,   -3.7090); 
\coordinate (rz-38) at (4.9982,   -3.7414); 
\coordinate (rz-39) at (5.1496,   -3.8382); 
\coordinate (rz-40) at (5.2529,   -3.9852); 
\coordinate (rz-41) at (5.2628,   -4.1244); 
\coordinate (rz-42) at (5.2727,   -4.2636);
\coordinate (rz-43) at (5.5650,   -4.5559); 
\coordinate (rz-44) at (5.5650,   -4.4026); 
\coordinate (rz-45) at (5.5650,   -4.2494); 
\coordinate (rz-46) at (5.5650,   -4.0962); 
\coordinate (rz-47) at (5.5720,   -3.9961); 
\coordinate (rz-48) at (5.5720,   -3.9956); 
\coordinate (rz-49) at (5.5720,   -3.8960);
\coordinate (rz-50) at (5.5720,   -3.8950); 
\coordinate (rz-51) at (5.6008,   -3.7024); 
\coordinate (rz-52) at (5.6842,   -3.5265); 
\coordinate (rz-53) at (5.8150,   -3.3823); 
\coordinate (rz-54) at (5.9821,   -3.2822); 
\coordinate (rz-55) at (6.1710,   -3.2350); 
\coordinate (rz-56) at (6.3655,   -3.2446);

\coordinate (vv-1) at (3.5396  , 0.0307  )   ;
\coordinate (vv-2) at ( 3.5396 ,  0.2278 )   ;
\coordinate (vv-3) at ( 3.5396 ,  0.4250 )   ;
\coordinate (vv-4) at ( 3.5396 ,  0.6221 )   ;
\coordinate (vv-5) at ( 3.5396 ,  0.8192 )   ;
\coordinate (vv-6) at ( 3.5396 ,  1.0164 )   ;
\coordinate (vv-7) at ( 3.5396 , 1.2135  )   ;
\coordinate (vv-8) at (3.5396  ,  1.4106 )   ;
\coordinate (vv-9) at ( 3.5396 ,  1.6078 )   ;
\coordinate (vv-10) at ( 3.5396 ,  1.8049 )  ;
\coordinate (vv-11) at ( 3.5396 ,  2.0020 )  ;
\coordinate (vv-12) at ( 3.5396 ,  2.1992 )  ;
\coordinate (vv-13) at ( 3.5396 , 2.3963  )  ;
\coordinate (vv-14) at ( 3.5396 ,  2.5935 )  ;
\coordinate (vv-15) at (3.5396  ,  2.7906 )  ;
\coordinate (vv-16) at ( 3.5396 ,  2.9877 )  ;
\coordinate (vv-17) at ( 3.5396 ,  3.1849 )  ;
\coordinate (vv-18) at ( 3.5396 ,  3.3820 )  ;
\coordinate (vv-19) at ( 3.5396 , 3.5791  )  ;
\coordinate (vv-20) at ( 3.5514 ,  3.7753 )  ;
\coordinate (vv-21) at ( 3.5869 ,  3.9686 )  ;
\coordinate (vv-22) at (3.6453  ,  4.1562 )  ;
\coordinate (vv-23) at ( 3.7260 ,  4.3354 )  ;
\coordinate (vv-24) at ( 3.8276 ,  4.5035 )  ;
\coordinate (vv-25) at ( 3.9489 , 4.6582  )  ;
\coordinate (vv-26) at ( 4.0878 ,  4.7971 )  ;
\coordinate (vv-27) at ( 4.2425 ,  4.9183 )  ;
\coordinate (vv-28) at ( 4.4107 ,  5.0200 )  ;
\coordinate (vv-29) at (4.5899  ,  5.1006 )  ;
\coordinate (vv-30) at ( 4.7775 ,  5.1590 )  ;
\coordinate (vv-31) at ( 4.9708 , 5.1944  )  ;
\coordinate (vv-32) at ( 5.1670 ,  5.2062 )  ;
\coordinate (vv-33) at ( 5.3631 ,  5.1943 )  ;
\coordinate (vv-34) at ( 5.5564 ,  5.1589 )  ;
\coordinate (vv-35) at ( 5.7440 ,  5.1004 )  ;
\coordinate (vv-36) at (5.9232  ,  5.0197 )  ;
\coordinate (vv-37) at ( 6.0913 , 4.9180  )  ;
\coordinate (vv-38) at ( 6.2474 ,  4.8102 )  ;
\coordinate (vv-39) at ( 6.4034 ,  4.7024 )  ;
\coordinate (vv-40) at ( 6.5595 ,  4.5946 )  ;
\coordinate (vv-41) at ( 6.7156 ,  4.4868 )  ;
\coordinate (vv-42) at ( 6.8716 ,  4.3790 )  ;
\coordinate (vv-43) at (7.0315  , 4.2640  )  ;
\coordinate (vv-44) at ( 7.1867 ,  4.1429 )  ;
\coordinate (vv-45) at ( 7.3371 ,  4.0158 )  ;
\coordinate (vv-46) at ( 7.4825 ,  3.8830 )  ;
\coordinate (vv-47) at ( 7.6226 ,  3.7446 )  ;
\coordinate (vv-48) at ( 7.7572 ,  3.6009 )  ;
\coordinate (vv-49) at ( 7.8861 , 3.4521  )  ;
\coordinate (vv-50) at (8.0092  ,  3.2984 )  ;
\coordinate (vv-51) at ( 8.1262 ,  3.1400 )  ;
\coordinate (vv-52) at ( 8.2370 ,  2.9772 )  ;
\coordinate (vv-53) at ( 8.3413 ,  2.8102 )  ;
\coordinate (vv-54) at ( 8.4391 ,  2.6393 )  ;
\coordinate (vv-55) at ( 8.5302 , 2.4647  )  ;
\coordinate (vv-56) at ( 8.6145 ,  2.2867 )  ;
\coordinate (vv-57) at (8.6918  ,  2.1056 ) ;
\coordinate (vv-58) at ( 8.7619 ,  1.9216 ) ;
\coordinate (vv-59) at ( 8.8249 ,  1.7351 ) ;
\coordinate (vv-60) at ( 8.8764 ,  1.5529 ) ;
\coordinate (vv-61) at ( 8.9166 , 1.3679  ) ;
\coordinate (vv-62) at ( 8.9453 ,  1.1809 ) ;
\coordinate (vv-63) at ( 8.9638 ,  1.0073 ) ;
\coordinate (vv-64) at (8.9768  ,  0.8333 ) ;
\coordinate (vv-65) at ( 8.9844 ,  0.6590 ) ;
\coordinate (vv-66) at ( 8.9865 ,  0.4845 ) ;
\coordinate (vv-67) at ( 8.9832 , 0.3100  ) ;
\coordinate (vv-68) at ( 8.9744 ,  0.1357 ) ;
\coordinate (vv-69) at ( 8.9602 , -0.0382 ) ;
\coordinate (vv-70) at ( 8.9364 , -0.2088 ) ;
\coordinate (vv-71) at (8.8985  , -0.3768 ) ;
\coordinate (vv-72) at ( 8.8468 , -0.5411 ) ;
\coordinate (vv-73) at ( 8.7817 ,-0.7006  ) ;
\coordinate (vv-74) at ( 8.7023 , -0.8743 ) ;
\coordinate (vv-75) at ( 8.6228 , -1.0480 ) ;
\coordinate (vv-76) at ( 8.5433 , -1.2217 ) ;
\coordinate (vv-77) at ( 8.4638 , -1.3954 ) ;
\coordinate (vv-78) at (8.3844  , -1.5691 ) ;
\coordinate (vv-79) at ( 8.3049 ,-1.7428  ) ;
\coordinate (vv-80) at ( 8.2254 , -1.9166 ) ;
\coordinate (vv-81) at ( 8.1459 , -2.0903 ) ;
\coordinate (vv-82) at ( 8.0665 , -2.2640 ) ;
\coordinate (vv-83) at ( 7.9870 , -2.4377 ) ;
\coordinate (vv-84) at ( 7.9075 , -2.6114 ) ;
\coordinate (vv-85) at (7.8280  ,-2.7851  ) ;
\coordinate (vv-86) at ( 7.7486 , -2.9588 ) ;
\coordinate (vv-87) at ( 7.6691 , -3.1325 ) ;
\coordinate (vv-88) at ( 7.5896 , -3.3062 ) ;
\coordinate (vv-89) at ( 7.5102 , -3.4799 ) ;
\coordinate (vv-90) at ( 7.4307 , -3.6536 ) ;
\coordinate (vv-91) at ( 7.3458 ,-3.8222  ) ;
\coordinate (vv-92) at (7.2486  , -3.9840 ) ;
\coordinate (vv-93) at ( 7.1396 , -4.1380 ) ;
\coordinate (vv-94) at ( 7.0195 , -4.2835 ) ;
\coordinate (vv-95) at ( 6.8888 , -4.4196 ) ;
\coordinate (vv-96) at ( 6.7483 , -4.5456 ) ;
\coordinate (vv-97) at ( 6.5988 ,-4.6608  ) ;
\coordinate (vv-98) at ( 6.4411 , -4.7645 ) ;
\coordinate (vv-99) at (6.2762  , -4.8561 ) ;
\coordinate (vv-100) at ( 6.1048 , -4.9352 ) ;
\coordinate (vv-101) at ( 5.9281 , -5.0013 ) ;
\coordinate (vv-102) at ( 5.7469 , -5.0540 ) ;
\coordinate (vv-103) at ( 5.5623 ,-5.0931  ) ;
\coordinate (vv-104) at ( 5.3753 , -5.1183 ) ;
\coordinate (vv-105) at ( 5.1869 , -5.1295 ) ;
\coordinate (vv-106) at (4.9897  , -5.1217 ) ;
\coordinate (vv-107) at ( 4.7949 , -5.0899 ) ;
\coordinate (vv-108) at ( 4.6054 , -5.0345 ) ;
\coordinate (vv-109) at ( 4.4242 ,-4.9564  ) ;
\coordinate (vv-110) at ( 4.2538 , -4.8567 ) ;
\coordinate (vv-111) at ( 4.0969 , -4.7369 ) ;
\coordinate (vv-112) at ( 3.9558 , -4.5989 ) ;
\coordinate (vv-113) at (3.8327  , -4.4446 ) ;
\coordinate (vv-114) at ( 3.7293 , -4.2765 ) ;
\coordinate (vv-115) at ( 3.6472 ,-4.0970  ) ;
\coordinate (vv-116) at ( 3.5877 , -3.9088 ) ;
\coordinate (vv-117) at ( 3.5516 , -3.7148 ) ;
\coordinate (vv-118) at ( 3.5396 , -3.5178 ) ;
\coordinate (vv-119) at ( 3.5396 , -3.3206 ) ;
\coordinate (vv-120) at (3.5396  , -3.1235 ) ;
\coordinate (vv-121) at ( 3.5396 ,-2.9264  ) ;
\coordinate (vv-122) at ( 3.5396 , -2.7292 ) ;
\coordinate (vv-123) at ( 3.5396 , -2.5321 ) ;
\coordinate (vv-124) at ( 3.5396 , -2.3349 ) ;
\coordinate (vv-125) at ( 3.5396 , -2.1378 ) ;
\coordinate (vv-126) at ( 3.5396 , -1.9407 ) ;
\coordinate (vv-127) at (3.5396  ,-1.7435  ) ;
\coordinate (vv-128) at ( 3.5396 , -1.5464 ) ;
\coordinate (vv-129) at ( 3.5396 , -1.3493 ) ;
\coordinate (vv-130) at ( 3.5396 , -1.1521 ) ;
\coordinate (vv-131) at ( 3.5396 , -0.9550 ) ;
\coordinate (vv-132) at ( 3.5396 , -0.7579 ) ;
\coordinate (vv-133) at ( 3.5396 ,-0.5607  ) ;
\coordinate (vv-134) at (3.5396  , -0.3636 ) ;
\coordinate (vv-135) at ( 3.5396 , -0.1665 ) ;

\coordinate (co-1) at (3.9431,   7.5741);
\coordinate (co-2) at (8.2851,   6.5398);
\coordinate (co-3) at (11.9919,  3.2752);
\coordinate (co-4) at (11.963,  -2.2336); 
\coordinate (co-5) at (8.3908,  -6.7269); 
\coordinate (co-6) at (4.334,   -7.4665);

\coordinate (so-1) at (1.696, -5.415 );
\coordinate (so-2) at (1.696, -3.6067);
\coordinate (so-3) at (1.696, -1.7983);
\coordinate (so-4) at (1.696,  1.8183);
\coordinate (so-5) at (1.696,  3.6267);
\coordinate (so-6) at (1.696,  5.4350);

 \draw[fill=blue!30] (vv-1)
--(vv-2)
--(vv-3)
--(vv-4)
--(vv-5)
--(vv-6)
--(vv-7)
--(vv-8)
--(vv-9)
--(vv-10)
--(vv-11)
--(vv-12)
--(vv-13)
--(vv-14)
--(vv-15)
--(vv-16)
--(vv-17)
--(vv-18)
--(vv-19)
--(vv-20)
--(vv-21)
--(vv-22)
--(vv-23)
--(vv-24)
--(vv-25)
--(vv-26)
--(vv-27)
--(vv-28)
--(vv-29)
--(vv-30)
--(vv-31)
--(vv-32)
--(vv-33)
--(vv-34)
--(vv-35)
--(vv-36)
--(vv-37)
--(vv-38)
--(vv-39)
--(vv-40)
--(vv-41)
--(vv-42)
--(vv-43)
--(vv-44)
--(vv-45)
--(vv-46)
--(vv-47)
--(vv-48)
--(vv-49)
--(vv-50)
--(vv-51)
--(vv-52)
--(vv-53)
--(vv-54)
--(vv-55)
--(vv-56)  
--(vv-57) 
--(vv-58) 
--(vv-59) 
--(vv-60) 
--(vv-61) 
--(vv-62) 
--(vv-63) 
--(vv-64) 
--(vv-65) 
--(vv-66) 
--(vv-67) 
--(vv-68) 
--(vv-69) 
--(vv-70) 
--(vv-71) 
--(vv-72) 
--(vv-73) 
--(vv-74) 
--(vv-75) 
--(vv-76) 
--(vv-77) 
--(vv-78) 
--(vv-79) 
--(vv-80) 
--(vv-81) 
--(vv-82) 
--(vv-83) 
--(vv-84) 
--(vv-85) 
--(vv-86) 
--(vv-87) 
--(vv-88) 
--(vv-89) 
--(vv-90) 
--(vv-91) 
--(vv-92) 
--(vv-93) 
--(vv-94) 
--(vv-95) 
--(vv-96) 
--(vv-97) 
--(vv-98) 
--(vv-99) 
--(vv-100)
--(vv-101)
--(vv-102)
--(vv-103)
--(vv-104)
--(vv-105)
--(vv-106)
--(vv-107)
--(vv-108)
--(vv-109)
--(vv-110)
--(vv-111)
--(vv-112)
--(vv-113)
--(vv-114)
--(vv-115)
--(vv-116)
--(vv-117)
--(vv-118)
--(vv-119)
--(vv-120)
--(vv-121)
--(vv-122)
--(vv-123)
--(vv-124)
--(vv-125)
--(vv-126)
--(vv-127)
--(vv-128)
--(vv-129)
--(vv-130)
--(vv-131)
--(vv-132)
--(vv-133)
--(vv-134)
--(vv-135)--cycle;
 
 \draw[fill=red!50] (rz-1)
--(rz-2)
--(rz-3)
--(rz-4)
--(rz-5)
--(rz-6)
--(rz-7)
--(rz-8)
--(rz-9)
--(rz-10)
--(rz-11)
--(rz-12)
--(rz-13)
--(rz-14)
--(rz-15)
--(rz-16)
--(rz-17)
--(rz-18)
--(rz-19)
--(rz-20)
--(rz-21)
--(rz-22)
--(rz-23)
--(rz-24)
--(rz-25)
--(rz-26)
--(rz-27)
--(rz-28)
--(rz-29)
--(rz-30)
--(rz-31)
--(rz-32)
--(rz-33)
--(rz-34)
--(rz-35)
--(rz-36)
--(rz-37)
--(rz-38)
--(rz-39)
--(rz-40)
--(rz-41)
--(rz-42)
--(rz-43)
--(rz-44)
--(rz-45)
--(rz-46)
--(rz-47)
--(rz-48)
--(rz-49)
--(rz-50)
--(rz-51)
--(rz-52)
--(rz-53)
--(rz-54)
--(rz-55)
--(rz-56)-- cycle;

\foreach \x in {1,2,3,4,5,6}{
   \fill (co-\x) circle (0.2cm);
}
\draw[] (co-3) node[above,yshift=.2cm] {Coils};

\draw[-, blue, very thick] (so-1) -- (so-2) ;
\draw[-, blue, very thick] (so-2) -- (so-3) ;
\draw[-, blue, very thick] (so-3) -- (so-4) ;
\draw[-, blue, very  thick] (so-4) -- (so-5) ;
\draw[-, blue, very  thick] (so-5) -- (so-6) ;

\foreach \x in {1,2,3,4,5,6}{
   \fill[blue] (so-\x) circle (0.1cm);
}

\draw[step=.2cm,gray, very thin] (3.19,-6) grid (9,6.);

\draw[] (so-6) node[above,yshift=.1cm] {Solenoids};

\node (wall) at (11,-4) {\small First wall};
\node (vv) at (11.7,-5.5) {\small Vacuum vessel};
\draw[->, thick] (wall) edge (7.4, -2);
\draw[->, thick] (vv) edge (6.5, -4);

\fill[fill=red!50]  (15,-2.5)--(15,3.1)--(22.5,2.)--(22.5, -2.5);
\draw[step=2.5cm, thick] (15,-2.51) grid (22.5,5);
\draw (15,3.1) --(22.5,2.);
\draw (5,4.2) rectangle (5.6,4.8);

\draw[dashed] (5,4.2) -- (15,-2.5) ;
\draw[dashed] (5,4.8) -- (15,5) ;

  \draw(16,-3.5) node[draw,fill=white,anchor=west] {\small Zoomed-in cut-cells};
\draw[] (6.2,.4) node[black, minimum width=3cm] {\small Plasma region};

\draw[|<->|, thick] ([xshift=-0.4cm]so-1) -- ([xshift=-0.4cm]so-6) node[midway, left]{\SI{10.85}{\metre}};


\end{tikzpicture}
\caption{ITER tokamak geometry. Left: ITER plasma region, external currents{,} and an exemplar cut-cell mesh.
Right: zoomed-in view. Note there are two types of external currents in ITER, i.e., six coils (point source) and solenoids of five segments (line source).
The geometry uses the actual dimensions from ITER.
}
\label{fig:intro}
\end{figure}

Computationally the free-boundary problem for the Grad-Shafranov
equation is more challenging as the plasma shape is not known {\em a
  priori}.  The general solution strategy involves an iterative
approach that combines a fixed-boundary solver with a Green's function
representation for the contribution to the poloidal magnetic flux on
the computational boundary from the coil current external to the
computational boundary.  The solution process is iterative by nature
and Picard iteration is usually deployed to drive the convergence of
the poloidal magnetic flux on the computational boundary.  The choice
of this fiducial computational boundary is constrained by two
considerations: (1) the computational boundary shall enclose all
chamber volume where current-carrying plasma can access; and (2) all~external currents should be outside the computational boundary.
\Cref{fig:intro} shows the plasma region (colored by red) and the external currents (coils and solenoids).
The most natural place for such a fiducial computational boundary in a
free-boundary Grad-Shafranov solver is the vacuum vessel  (the blue region in~\Cref{fig:intro}) of which the
poloidal field currents are all outside, or the first wall (as indicated in~\Cref{fig:intro}).
Both first wall and vacuum
vessel in modern tokamaks are strongly shaped and give rise to an
irregular computational domain, for which some variation of finite
element/volume on an unstructured grid would be a ready choice for the
underlying fixed-boundary Grad-Shafranov solver.  It is interesting to
note that the essential ingredients of free-boundary Grad-Shafranov
solvers, for example, both the Green's function formulation and
iterative solution procedure, were first developed with finite
difference discretization over a structured grid in a regular domain~\cite{lackner1976computation, johnson1979numerical, jardin1986dynamic, jeon2015development}.  This is partly due to the considerable freedom in the
placement of the fiducial computational boundary despite the two
general constraints noted earlier, and partly due to the fact that
many tokamaks have enough space between current-carry plasmas and
poloidal field coils for a rectangular computational boundary
in-between them. For some of the next generation tokamaks, it becomes
more desirable to have the computational boundary aligned with the
better flux conserver, for example the vacuum vessel in ITER~\cite{ariola-etal-ITER-shape-2000}, or with
the first wall, such as the case with the SPARC tokamak~\cite{creely-sparc-jpp-2020} where the coils are very close to the plasma region.
Also, note that it is desirable to represent the field on a structured grid, which can significantly accelerate the particle pusher in a particle-in-cell code.
We
find that for these applications, cut-cells (as indicated in~\Cref{fig:intro}) can be used in combination
with finite difference in a structured grid to accommodate an
irregular computational domain for a free-boundary Grad-Shafranov
solver.

Numerical solutions to elliptic problems with complex geometry, of
which the Grad-Shafranov equilibrium is one example,  have  been
considered by many approaches, including finite difference, finite
volume, and finite element using various meshing techniques.  Among
those numerical approaches, the cut-cell approach has a number of
advantages, which include simplifying the grid generation process for
complicated geometries, enabling fast computation of the solution in
parallel, and shifting the complexity of dealing with complex
geometries to the discretization approach.
For more details, one can consult the review
in~\cite{berger2017cut}.  The cut-cell approach, a.k.a.~the Cartesian
grid embedded boundary method, generates the mesh using a background
regular mesh and taking special care of cut-cells where the geometry
intersects the grid. A variety of work have addressed the successfully
use of the cut-cell approach for solving {elliptic equations} with
finite volume~\cite{johansen1998cartesian, schwartz2006cartesian,
  devendran2014higher, devendran2017fourth} and finite difference
discretization \cite{jomaa2005embedded,
  gibou2002second,gibou2005fourth, leveque1994immersed}.  When
discretized using the flux-divergence form, finite volume methods have
the advantage of keeping discretely
conservative~\cite{leveque1992numerical}, which is important in
solving many problems such as heat and mass transfer. In this work, we
employ the cut-cell approach with a conservative discretization to
address the Grad-Shafranov equation in irregular domains.

{The paper is organized as follows. \Cref{sec:GSequation} gives a
  brief description of the Grad-Shafranov equation and the associated fixed-boundary
  and free-boundary problems.~\Cref{sec:cut-cellalgorithm}
  discusses the cut-cell algorithm, which is a Cartesian grid embedded
  boundary method, as well as a level-set approach to improve the
  efficiency of the cut-cell algorithm.~This is followed by a detailed
  description of a free-boundary Grad-Shafranov solver, focusing on
  solving the free-boundary problem reformulated on a bounded domain
  in \Cref{sec:FreeBDGSsolver}. In~\Cref{sec:implementation},
  the parallel implementation is described and some details of the
  implementation of the free-boundary Grad-Shafranov solver are
  presented.  Finally, numerical tests are presented
  in~\Cref{sec:numericaltests}, demonstrating the accuracy and
  efficiency of the free-boundary Grad-Shafranov solver.}

\section{Grad-Shafranov equation}
\label{sec:GSequation} 
In this section, the Grad-Shafranov equation, fixed-boundary problem,  and free-boundary problem are briefly introduced.

The Grad-Shafranov equation is derived from the MHD equilibrium equations given by 
\begin{subequations}
\begin{align}
\nabla p &= \vect{J} \times\vect{B},\label{eq:momentumequation}\\
 \mu_0\vect{J} &= \nabla\times\vect{B}, \label{eq:amperelaw}\\
 \nabla \cdot \vect{B} &= 0, \label{eq:divergencefree}
\end{align}
\label{eq:MHD}
\end{subequations}
where $\vect{J}$ is the electric current, $\vect{B}$ is the magnetic field,  $\mu_0$ is
the magnetic permeability,  and $p$ is the plasma pressure. 
\Cref{eq:momentumequation} is the force balance
equation. \Cref{eq:amperelaw} is the Amp\`ere's law. \Cref{eq:divergencefree} is the divergence
free constraint for the magnetic field.
In the case of equilibrium, the above system~\eqref{eq:MHD} holds in the entire region $\mathbb{R}^3$.

The Grad-Shafranov equation is usually written in a cylindrical coordinate system, denoted by $(R,\phi, Z)$,
for an axisymmetric plasma such as those in a tokamak. 
It is well known that the magnetic field in the tokamak can be represented by
\begin{equation*}
{\bf B} =\nabla \phi \times \nabla \psi + g(\psi) \nabla \phi,
\end{equation*}
where $\psi$ is commonly referred to as the poloidal flux function and $g$ is a scalar function.
The  Grad-Shafranov equation can be derived from~\Cref{eq:MHD} as
\begin{equation}\label{eq:GSeq}
\Delta^* \psi=\mu_0RJ_{\phi}(R,\psi), 
\end{equation}
where
the toroidal elliptic operator is defined as 
 \begin{displaymath}
\Delta^* \psi \equiv  \frac{\partial^2 \psi}{\partial R^2}-\frac{1}{R}\frac{\partial \psi}{\partial R}+ \frac{\partial^2 \psi}{\partial Z^2},
\end{displaymath}
and the source term satisfies
\begin{equation*}
\mu_0RJ_{\phi}(R,\psi)=-\left[\mu_0R^2 \, \frac{d \, p(\psi)}{d\psi}+g(\psi) \, \frac{d \,g(\psi)}{d\psi}\right],
\end{equation*}
in which the plasma pressure, $p$, is a scalar function of $\psi$.
Note that the  toroidal elliptic operator can be rewritten in a divergence form
\begin{equation*} 
\Delta^*\psi =R  \tilde\nabla\cdot \Big(\frac 1 R \tilde\nabla \psi \Big),
\end{equation*}
where $\tilde\nabla\equiv [\partial_R, \partial_Z]^T$.
This conservative form will be used when discretizing the problem with the cut-cell approach.

Note that in a tokamak plasma, $J_\phi$ is carried by the magnetically
confined plasma that resides inside the magnetic separatrix in the
formulation of both fixed- and free-boundary Grad-Shafranov
equilibria, which we shall explain next.

\subsection{Fixed-boundary problem}

The fixed-boundary problem refers to the situation when the boundary
condition for the equilibrium problem is known.  Mathematically, the
problem can be simply~cast as a boundary value problem given by
\begin{equation}\label{eq:fixedbdeq}
  \begin{aligned}
   & \Delta^* \psi= \mu_0RJ_{\phi}(R,\psi), && (R,Z)\in \Omega.\\
  & \psi=\psi_D,&& (R,Z)\in \partial \Omega.
  \end{aligned}
\end{equation} 
where $\Omega$ is the physical domain with a Lipschitz boundary $
\partial \Omega$.  In tokamaks, the physical domain $\Omega$
corresponds to the cross section of the device bounded by the first
wall or more commonly, the last closed flux surface beyond which
$J_\phi$ vanishes.  In this work, the fixed-boundary problem is also
implemented as an important tool to verify the numerical scheme.

 \subsection{Free-boundary problem}
 
The free-boundary problem refers to the situation when the plasma
domain, denoted as $\mathcal{P}(\psi)$, is unknown.  In tokamak, the
plasma domain refers to the region inside the device where the
confinement actually happens, and thus $\mathcal{P}(\psi)$ is filled
with hot plasmas.  In a diverted tokamak plasma, $\mathcal{P}(\psi)$
is the magnetic separatrix that demarcates the region of nested closed
flux surfaces and the scrape-off layer that has magnetic field lines
intercept the divertor and first wall.  For a limited plasma,
$\partial \mathcal{P}(\psi)$  is the last closed flux surface beyond which the
flux surfaces are intercepted by one or multiple limiters that
protrude from the chamber wall.  It is important to note that
$\mathcal{P}(\psi)$ is only meaningful when such an MHD equilibrium
exists.  Such a constraint poses an extra challenge to numerical
algorithms, since the success of finding a numerical solution to the
free-boundary problem implicitly indicates the existence of such a
domain, while a perturbed numerical solution may indicate an
unrealistic plasma domain, resulting into the divergence of
algorithms.  Therefore, it is necessary to design the numerical
algorithm to be robust with respect to such a perturbation.
 
Assuming an unbounded domain $\mathcal{H}:=[0,\infty) \times (-\infty,
  \infty)$, the free boundary problem is given by
\begin{equation}\label{eq:freebdproblempart1}
 \Delta^* \psi=\left\{
 \begin{aligned}
 &\mu_0RJ_{\phi}(R,\psi), && (R,Z)\in \mathcal{P}(\psi),\\
 &\mu_0RI_i, && (R,Z)=(R_i^c,Z_i^c), \, i=1,2, \cdots, n_c.\\
 &0, && \text{otherwise}.
 \end{aligned}
 \right.
 \end{equation}
 with the asymptotic constraint
\begin{equation}\label{eq:freebdproblempart2}
\psi(0,Z)=0 \quad \text{and} \quad \lim_{\lVert (R,Z) \rVert \to \infty}\psi(R,Z)=0.
\end{equation}
Here $I_i$ is the external current density in the $i$-th poloidal
field coil located at $(R_i^c,Z_i^c)\in \mathcal{H}$ (see \Cref{fig:intro} for the locations of coils; solenoids can be handled similarly), and
$J_{\phi}(R,\psi)$ is the prescribed toroidal component of the plasma
current density, generally a non-linear function of $\psi$, in the
plasma domain $\mathcal{P}(\psi)$.  Note that the magnetic potential
outside of the plasma region simply satisfies Poisson's equation.  We
further introduce the limiter domain $ \mathcal{L}$, satisfying
$\mathcal{P}(\psi)\subset \mathcal{L} \subset \mathcal{H}$.  For this
purpose, the limiter corresponds to the first wall in the
tokamak device surrounding the confinement region.  The limiter domain
is the entire region accessible by the plasma, while the plasma domain
$\mathcal{P}(\psi)$ is the region bounded by the last closed poloidal
flux surface inside the limiter bounded domain, see the discussion
in~\cite{faugeras2017fem} for instance.

The free-boundary problem is still underdetermined because (a) the
source  term $J_\phi$ (more specifically,  $p$ and  $g$) is~typically provided as a function of  the normalized $\bar\psi \in [0, 1]$ from experiment (its normalization shall be explained later); (b) the coil currents
are not known and in fact they need to be optimized based on the
targeted shape of $\mathcal{P}(\psi)$, along with solving the
free-boundary problem.  Both of the uncertainties introduce more
challenges to solve the free-boundary problem.  However, the
uncertainties are closely related to the concept of plasma shape
control in tokamaks, which is a critical topic in tokamak discharge
design as well as in tokamak machine design.

To address (a), one can define a normalized flux function,
\begin{displaymath}
\bar{\psi} \equiv  \frac{\psi - \psi_o}{\psi_{ X} - \psi_o},
\end{displaymath}
where $\psi_{ X}$ is the value of $\psi$ at the plasma boundary,
$\partial \mathcal{P}(\psi)$, and $\psi_o$ is the value at the
magnetic axis.  In a diverted tokamak plasma, a magnetic
\emph{separatrix} demarcates the topologically distinct closed and open
  flux regions. Plasmas also reside in the open flux region outside
  the separatrix, which is known as the scrape-off layer, but it has
  relatively small pressure gradient and plasma current.  For the
  purpose of defining a Grad-Shafranov equilibrium, the plasma
  boundary is usually given by the separatrix or a closed flux surface
  in the immediate vicinity of the separatrix.  In this paper, we will
  set $\psi_{ X}$ to be the $\psi$ value of the saddle point that
  labels the magnetic separatrix.
  In a limited plasma, $\psi_{ X}$ labels the last closed flux surface beyond which
  magnetic field lines interact a limiter protruded from camber wall.
  Since $\psi$ is the poloidal flux,
  it must be a monotonic function between $\psi_o$ and $\psi_{ X}$ in a tokamak.
  The normalized flux $\bar{\psi}$ is therefore between 0 and 1 for
  all $\psi$ values inside $\mathcal{P}(\psi)$.  Numerical techniques
  to efficiently locate those points will be discussed in the
  following section.  Note it is expected that the overall source term
  can have jumps on the right hand side
  in~\cref{eq:freebdproblempart1}.  The existence of solutions
  to~\cref{eq:freebdproblempart1} under certain normalization and jump
  source terms has been considered in some PDE theory studies,
  see~\cite{ambrosetti1990remarks} for instance.

In addressing (b), one usually imposes some constraints based
on a targeted plasma shape $\mathcal{P}_{\rm target}$.  The
targeted plasma shape is predetermined (up to some small variation) when the
fusion device is designed.  The standard approach is to impose a few control
points along the desired shape such that the computed poloidal flux
on those points is equal to the value on the separatrix, i.e.,
\begin{displaymath}
\psi(R_k,Z_k)= \psi_X, \quad k=1,\cdots, n_p, 
\end{displaymath}
where $(R_k,Z_k)$ is selected along the plasma boundary
$\partial\mathcal{P}_{\rm target}$ and $\psi_X$ is determined
numerically during the normalization.  It essentially means we
minimize the distance between the separatrix from the numerical
solution and the targeted plasma shape.  In the following section, we
will discuss how to impose a proper optimization problem to determine
the coil current based on the above constraints.

\section{Cut-cell algorithm} 
\label{sec:cut-cellalgorithm}
\def\xv{\mathbf{x}}

A cut-cell algorithm for the free-boundary Grad-Shafranov problem is described in this section.
The basic algorithm closely follows the cut-cell algorithm proposed for Poisson's equation in~\cite{johansen1998cartesian}.
We further introduce a level set function to improve the efficiency of the cut-cell algorithm and suit the specific needs of the Grad-Shafranov equation.

Consider a uniform grid defined by 
\[
\xv_{i,j} = (R_{i}, Z_{j}) =  (R_{\min}+i\Delta R, \, Z_{\min} + j\Delta Z).
\]
Let $\psi_{i,j}$ be a grid-point approximation to $\psi(\mathbf{x}_{i,j})$ and rewrite \Cref{eq:GSeq} as 
\begin{equation*}
\tilde\nabla \cdot \Big (\frac{1}{R} \tilde\nabla \psi \Big)=\mu_0J_{\phi}(R,\psi)
\end{equation*}
Integrating it on the cell corresponding to the grid point $(i,j)$  and applying the divergence theorem give
\begin{equation}\label{eq:conservationlawforGS}
\frac{1}{\Delta R \Delta Z}\int_{\partial C} \Big(\frac{1}{R} \tilde\nabla \psi \Big)\cdot \vect{n} \, ds=\frac{1}{\Delta R \Delta Z} \int_{C}\mu_0J_{\phi}(R,\psi) \,dRdZ, \quad C\in \Omega
\end{equation}
where $\vect{n}$ denotes the unit outward normal direction vector of $\partial C$. 
Assuming $C$ is a full cell, \Cref{eq:conservationlawforGS} can be discretized as 
\begin{equation}\label{eq:fullcellGSeq}
\begin{aligned}
\frac{1}{\Delta R \Delta Z}\Big[ F_{i,j+\frac 1 2}-F_{i,j-\frac 1 2}+F_{i+\frac 1 2,j}-F_{i-\frac 1 2,j}\Big]={\mu_0 J_{\phi}(R_{i,j}, \psi_{i,j})},
\end{aligned}
\end{equation}
where $F_{r,s}$ is the flux along the cell edge, given by
$$
F_{i,j+\frac 1 2}= \frac{\Delta R}{R_{i,j+\frac 1 2}}\frac{\psi_{i,j+1}-\psi_{i,j}}{\Delta Z}, \qquad F_{i+\frac 1 2,j}=\frac{\Delta Z}{R_{i+\frac 1 2,j}}\frac{\psi_{i+1,j}-\psi_{i,j}}{\Delta R}.
$$
Note that~\cref{eq:fullcellGSeq} on the full cell is a standard second-order discrete operator.

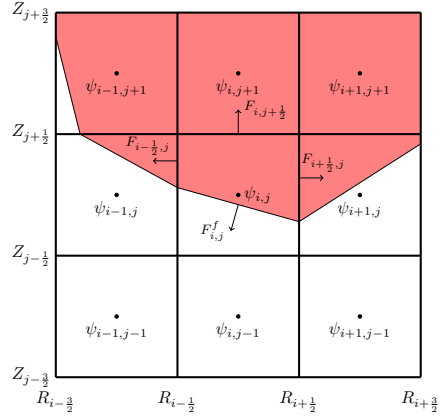
\begin{figure}
  \centering
\resizebox{6cm}{!}{
\begin{tikzpicture}

\fill[fill=red!50]  (5,7)--(5.5,5)--(7.5,3.9)--(10,3.2)--(12.5,4.8)--(12.5,7.5)--(12.5,2.5)--(12.5,7.5)--(5,7.5);

      \filldraw [black] (6.25,1.25) circle (1pt) node[anchor=north] {$\psi_{i-1,j-1}$};
      \filldraw [black] (8.75,1.25) circle (1pt)node[anchor=north] { $\psi_{i,j-1}$};
        \filldraw [black] (11.25,1.25) circle (1pt)node[anchor=north] {$\psi_{i+1,j-1}$}; 
         \filldraw [black] (6.25,3.75) circle (1pt)node[anchor=north] {$\psi_{i-1,j}$};
       \filldraw [black] (8.75,3.75) circle (1pt)node[anchor=west] {$\psi_{i,j}$};
        \filldraw [black] (11.25,3.75) circle (1pt)node[anchor=north] {$\psi_{i+1,j}$};
          \filldraw [black] (6.25,6.25) circle (1pt)node[anchor=north] {$\psi_{i-1,j+1}$};
       \filldraw [black] (8.75,6.25) circle (1pt)node[anchor=north] {$\psi_{i,j+1}$};
        \filldraw [black] (11.25,6.25) circle (1pt)node[anchor=north] {$\psi_{i+1,j+1}$};
     
    \node (A) at (5,-0.5) {$R_{i-\frac{3}{2}}$};
     \node (A) at (7.5,-0.5) {$R_{i-\frac{1}{2}}$};
 \node (A) at (10,-0.5) {$R_{i+\frac{1}{2}}$};
 \node (A) at (12.5,-0.5) {$R_{i+\frac{3}{2}}$};
 
  \node (A) at (4.5,0) {$Z_{j-\frac{3}{2}}$};
     \node (A) at (4.5,2.5) {$Z_{j-\frac{1}{2}}$};
 \node (A) at (4.5,5) {$Z_{j+\frac{1}{2}}$};
 \node (A) at (4.5,7.5) {$Z_{j+\frac{3}{2}}$};

    \draw[very thick] (5,0) --(5,7.5);
    \draw[very thick] (7.5,0) --(7.5,7.5);
    \draw[very thick] (10,0) --(10,7.5);
    \draw[very thick] (12.5,0) --(12.5,7.5);
        
    \draw[very thick] (5,0) --(12.5,0) ;
    \draw[very thick] (5,2.5) --(12.5,2.5);
    \draw[very thick](5,5) --(12.5,5);
    \draw[very thick] (5,7.5) --(12.5,7.5);
    
    \draw (5,7)--(5.5,5);
     \draw (5.5,5)--(7.5,3.9);
      \draw (7.5,3.9)--(10,3.2);
      \draw (10,3.2)--(12.5,4.8);
      
       \draw [->] (7.5,4.45)--(7,4.45) node[anchor=south, xshift=-3pt] {\footnotesize $F_{i-\frac 1 2,j}$};
       \draw [->] (8.75,5)--(8.75,5.5) node[anchor=west] {\footnotesize $F_{i,j+\frac 1 2}$};
       \draw [->] (10,4.1)--(10.5,4.1)node[anchor=south] {\footnotesize $F_{i+\frac 1 2,j}$}; 
       \draw [->] (8.75,3.55)--(8.6,3.014)node[anchor=east] {\footnotesize $F_{i,j}^f$};

\end{tikzpicture}
}
 \caption{Grid points and cut-cells. The red region is the computational domain. }
 \label{fig:cutcelldiagram}
\end{figure}

If $C$ is a cut-cell as shown in Figure~\ref{fig:cutcelldiagram}, introducing a volume fraction $\Lambda_{i,j} \in [0,1]$, the resulting discretization of \cref{eq:conservationlawforGS} can be written as 
\begin{equation}\label{eq:cutcellGSeq}
\begin{aligned}
\frac{1}{\Delta R \Delta Z  \Lambda_{i,j}}\Big[ F_{i,j+\frac 1 2}-F_{i,j-\frac 1 2}+F_{i+\frac 1 2,j}-F_{i-\frac 1 2,j}- F^f_{i,j}\Big]={\mu_0 J_{\phi}( R_{i,j}, \psi_{i,j} )},
\end{aligned}
\end{equation}
where $F^f_{i,j}$ is the flux along the boundary in the cut-cell. 
Note that $F_{r,s}$ can be  0 if the corresponding edge is not involved in the control volume.

{We further introduce a level set function for efficiently distinguishing interior cell, exterior cell, and cut-cell as well as {determining} whether a cell edge is a partial edge or a full edge.} 
We first construct a level set function $\rho(\footnotesize R,Z)$ given by 
\begin{equation}\nonumber
\rho(R,Z)=\left\{
\begin{array}{rl}
+d,&(R,Z)\in \Omega^c,\\
0,&(R,Z)\in \partial \Omega,\\
-d,&(R,Z)\in \Omega,
\end{array}
\right.
\end{equation}
where $d$ is the shortest distance from the point $(R,Z)$ to the boundary $\partial \Omega$. 
For convenience, let $\rho_{i+\frac 1 2,j+\frac 1 2}$ denote $\rho(\xv_{i+\frac 1 2,j+\frac 1 2})$, which is evaluated on the corners of each cell.

For the boundary described by an analytical expression, defining a level set value on a given point is  straightforward and therefore it is ignored. However, in a practical problem of a tokamak geometry,  
the geometry is described by a set of points on the boundary, and thus it is less obvious to define the value of level set function on a point. Here, we present how to define the level set function for a general irregular domain $\Omega$. Given a set of data $\{(r_1,z_1), (r_2,z_2), (r_3,z_3), \cdots, (r_n,z_n)\}$ which describes the boundary of the irregular domain, we connect all the points as a polygon. Let
$$
\bar{r}=\sum_{i=1}^{n}\frac{r_i}{n}, \quad \bar{z}=\sum_{i=1}^{n}\frac{z_i}{n}, \quad \bar{P}=(\bar{r},\bar{z}), \quad P_{i}=(r_{i}, z_{i}),  \quad i=1, 2,\cdots, n.
$$
To compute $d$ for a given point $P=(r,z)$, we locate the boundary point $P_i$ which has the shortest distance to  $P$. Then we determine $d$ as the shorter distance between $P$ and the two boundary segments connected with $P_i$.

The next step is to decide whether $P$ is in or out of domain if $d\neq 0$. 
We choose $(r_1,z_1)$ as a starting point and compute the rotated angle $\theta_i$ from $\overrightarrow{\bar{P}P_1}$ to $\overrightarrow{\bar{P}P_i}$ where $\theta_i\in (-\pi,\pi]$.
For any given point $P=(r,z)$, by checking the value of the angle $\theta$ rotating from $\overrightarrow{\bar{P}P_1}$ to $\overrightarrow{\bar{P}P}$, we can get the interval $[\theta_i, \theta_{i+1}]$ where $\theta$ belongs. 
We denote $P_{intersect}$ as the intersection point of line passing through $P$ and $\bar{P}$ and the line passing through $P_i$ and $P_{i+1}$. By comparing values of  the distance between $\bar{P}$ and $P_{intersect}$ and the distance between $\bar{P}$ and $P$, we can determine whether the point $P$ is in or out of domain. Some care is needed for the corner point where the curve becomes not convex.

The precomputed level set function can be used to  distinguish the cell category. 
 For example, in \cref{fig:cutcelldiagram}, since $\rho$  at each corner of the cell $(i,j-1)$  is positive, 
 we determine the cell $(i,j-1)$ is an exterior cell. Similarly, the cell $(i,j+1)$ is an interior cell in the domain $\Omega$, as $\rho $ values at all the corners of the cell are negative. The cell $(i,j)$ is a cut-cell which contains part of the boundary as signs of $\rho $ value  at each corner are not the same. 
The precomputed level set function can also be used to evaluate some geometry quantities associated with the cut-cells~\cite{liu2018numerical, liu2020numerical}.
 For a partial edge, we can use the information of the level set function~$\rho$ at the two endpoints to approximate its aperture. 
Suppose $\xv_f = (R_{i+\frac 1 2},  Z_f)\in \partial \Omega$, assuming two points $\xv_{i+\frac 1 2,j+\frac 1 2}\in \Omega$ and $\xv_{i+\frac 1 2, j-\frac 1 2} \in \Omega^c$ border $\xv_f $, we use the formula given in~\cite{chen1997simple}
\begin{subequations}
\begin{align*}
Z_{i+\frac 1 2, j+\frac 1 2}-Z _f&= \frac{\rho_{i+\frac 1 2,j+\frac 1 2}\triangle Z}{\rho_{i+\frac 1 2,j+\frac 1 2}-\rho_{i+\frac 1 2,j-\frac 1 2}}, \\
Z_f-Z_{i+\frac 1 2, j-\frac 1 2}&= \frac{-\rho_{i+\frac 1 2,j - \frac 1 2}\triangle Z}{\rho_{i+\frac 1 2,j+\frac 1 2}-\rho_{i+\frac 1 2,j-\frac 1 2}}, 
 \end{align*}
\end{subequations}
to evaluate the distances between $\xv_f$ and two points $\xv_{i+\frac 1 2,j+\frac 1 2}$ and $\xv_{i+\frac 1 2, j-\frac 1 2}$.
Therefore, as an example, the fraction of the right edge of the partial cell $(i,j)$ in \Cref{fig:cutcelldiagram}, denoted as $a_{i+\frac 1 2,j}$, is given by
\[
a_{i+\frac 1 2,j} = \frac{\rho_{i+\frac 1 2,j+\frac 1 2}}{\rho_{i+\frac 1 2,j+\frac 1 2}-\rho_{i-\frac 1 2,j+\frac 1 2}}.
\]

After calculating the apertures of each cell edge, the area of the front, denoted as $A^f$, in the cut-cells can be evaluated based on the relationship
\begin{equation}\nonumber
A^f_{i,j}{\mathbf n^{in}_{i,j}}=\Delta R(a_{i+\frac 1 2,j}-a_{i-\frac 1 2,j})\hat{i}+\Delta Z(a_{i,j+\frac 1 2}-a_{i,j-\frac 1 2})\hat{j},
\end{equation}
where ${\mathbf n}^{in}$ is the inward-facing normal at the boundary in cut-cell. Note that for interior cells, since all aperture values are 1, $A^f$ is also applicable to the interior cells as $A^f=0$.

Next, we discuss the discretization of \Cref{eq:cutcellGSeq} on cut-cells.
$F_{i,j+\frac 1 2}$, $F_{i,j-\frac 1 2}$, $F_{i+\frac 1 2,j}${,} and $F_{i-\frac 1 2,j}$ are evaluated at the midpoint of cell edge covered by $\Omega$. To evaluate flux on a partial edge, in order to  keep a second-order accurate approximation of the fluxes through cell edges, a linear interpolation between values at the midpoints of full edges is used. For example, in \cref{fig:cutcelldiagram}, the flux $F_{i+\frac 1 2,j}$ is evaluated at the midpoint of the partial edge centered at $\xv_{i+\frac 1 2,j}$ by linearly interpolating the flux at neighboring full edge centered at $\xv_{i+\frac{1}{2},j+1}$, which is stated as the following 
 \begin{equation}\nonumber
F_{i+\frac 1 2,j}=\frac{1}{R_{i+\frac 1 2,j}}\Delta Z \left[\frac{1+a_{i+\frac 1 2,j}}{2}\frac{(\psi_{i+1,j}-\psi_{i,j})}{\Delta R}+\frac{1-a_{i+\frac 1 2,j}}{2}\frac{(\psi_{i+1,j+1}-\psi_{i,j+1})}{\Delta R}\right].
\end{equation}
The flux $F^f$ is evaluated at the midpoint of the boundary covered by the cut-cell. A three-point gradient stencil is proposed in~\cite{johansen1998cartesian} to evaluate  the normal component of the gradient of $\psi$ and deployed in the current work. 
The readers are referred to~\cite{johansen1998cartesian} for the details of the discretization and other techniques such as small cell ignorance, which are also deployed here.

Finally, we discuss how to evaluate the source term of the Grad-Shafranov equation along with the cut-cells. Note that ${\mu_0 J_{\phi}(R_{i,j},\psi_{i,j})}$ should be evaluated at the centroid of the covered part in each cut-cell. There are three types of cut-cells covered by the physical domain, triangle, trapezoid, and pentagon, based on the piecewise-linear representation of the boundary in each cut-cell. With the geometry information gathered, we can now turn the problem into computing the centroid of a convex and closed polygon. Suppose a convex and closed polygon is defined by the ordered vertices $\{(r_1,z_1), (r_2,z_2),  \cdots, (r_{n+1},z_{n+1})\}$ with $(r_{n+1}, z_{n+1}) = (r_1, z_1)$, the centroid of a convex polygon is computed as the following
\begin{displaymath}
C_{polygon}=\frac{1}{3}\left (\frac{\sum_{i=1}^n(r_i+r_{i+1})(r_i z_{i+1}-r_{i+1}z_i)}{\sum_{i=1}^n(r_i z_{i+1}-r_{i+1} z_i)},\frac{\sum_{i=1}^n(z_i+z_{i+1})(r_i z_{i+1}-r_{i+1} z_i)}{\sum_{i=1}^n(r_i z_{i+1}-r_{i+1} z_i)}\right ),
\end{displaymath}
which will be used during the evaluation of the source term.

\section{Free-boundary Grad-Shafranov solver}
\label{sec:FreeBDGSsolver}

A common approach to solve the free-boundary
problem~\cref{eq:freebdproblempart1} is to reformulate the problem in
a bounded computational domain $\Omega$, see
\cite{jardin2010computational, faugeras2017fem, jeon2015development}
for instance.  The requirement for $\Omega$ is that it should enclose
the plasma domain $\mathcal{P}(\psi),$ but all the poloidal field
coils, including those for vertical stability control and divertor
flux shaping, should be outside $\Omega.$ This is to accommodate a
Green's function approach in calculating the coil current contribution
to the Grad-Shafranov equilibrium through the poloidal flux $\psi$ on
$\partial\Omega.$ As the plasma domain $\mathcal{P}(\psi)$ is unknown
\emph{a priori}, the bounded domain $\Omega,$ in practice, should
contain at least the limiter bounded domain $\mathcal{L}$. The precise
choice for $\partial\Omega$ constrains the numerical discretization
for the Grad-Shafranov solver and motivates the cut-cell approach in
section~\ref{sec:cut-cellalgorithm}.

In this section we will focus on the overall scheme for the
free-boundary Grad-Shafranov solver. 
The main algorithm
is based on a Picard iteration
\begin{alignat*}{2}
& \Delta^* \psi = \mathcal{S}( \psi^{\rm old}), \quad && (R, Z) \in \Omega, \\
& \psi =  \mathcal{D}(\psi^{\rm old}), \quad && (R, Z) \in \partial\Omega,
\end{alignat*}
where the operators $\mathcal{S}$ and $\mathcal{D}$ stand for the
steps to determine the source term and boundary values from the old
solution $\psi^{\rm old}$, respectively.  The details of $\mathcal{S}$
and $\mathcal{D}$ will be given in the following discussion.  Note
that unlike a fixed-boundary problem~\cite{dpgGS}, Newton's method is
challenging to employ here, which will be  further discussed in \Cref{sec:bc}.

The discussions will start from determining the source term from an
old solution, i.e., the operator $\mathcal{S}$, which includes a
search algorithm to locate $\psi_o$ and $\psi_X$ and a barycentric
interpolation to interpolate the given source profile.  After the
source term is determined, the approaches to determine the boundary
condition, i.e., the operator $\mathcal{D}$, are discussed.  Then a
minimization problem is introduced to optimize the coil current
density based on the old solution.  Finally, the full algorithm based
on Picard iteration with Aitken's acceleration is summarized.  Note
that all the approaches in this section are applicable to a
free-boundary solver on a rectangular domain and some of them were
actually initially proposed for the rectangular domain.

\subsection{Evaluate the source term} 
Evaluating the source term in the free-boundary problem is a
multi-step process.  Since $\psi_o$ and $\psi_{ X}$ are to be found
from the solution, Grad-Shafranov equilibria are usually specified by
prescribing the source term as functions of $\bar{\psi},$ which is
always well-defined for any trial solution that has varying $\psi_o$
and $\psi_{ X}.$ The most straightforward case has known
$g(\bar\psi)$ and $p(\bar\psi)$ profile, often supplemented by
integral constraints of specified plasma beta and plasma toroidal
current. A practically useful alternative is to specify the safety factor profile, $q(\bar{\psi}),$
for $\bar{\psi} \in [0, 0.95].$

\paragraph*{{Search  for magnetic axis and saddle points.}}
Since the trial solution of Grad-Shafranov equation gives $\psi$ on
grid points, it is necessary to locate $\psi_o$ and $\psi_X$ first.
As discussed in~\Cref{sec:GSequation}, $\psi_o$ and $\psi_X$
correspond to the critical points where $\| \nabla \psi \|=0$.
Therefore, a simple minimization problem can be defined to find them,
\begin{equation}
\label{eqn:searchPoints}
\min_{R, Z} \| \nabla \psi \|^2.
\end{equation}
To simplify the procedure, we use a two-step search algorithm to locate those points.
We first evaluate numerical gradient on all the grid points in order to identify several candidate points with small gradient values. We choose 10 candidate points in our solver (here 10 points are chosen to locate all the saddle points and magnetic axis and they often converge to one of those points. During the searching process, if the point locates out of the computational domain, we delete the corresponding candidate).
Newton's method for optimization is then used to
solve~\Cref{eqn:searchPoints} using the candidate points as the
initial guess.  During the optimization, the Newton's method needs a
continuous representation of $\psi$ for any $(R,Z)$, for which a
barycentric interpolation formula is used. Since the initial
guesses are very good, the Newton's method typically only needs a few
iterations to converge.  After all the critical points are found, a
further check is performed to remove the duplicate critical points.
Finally, $\psi_o$ and $\psi_X$ are determined based on the Hessian of
$\psi$.  If the critical point has $\partial_{rr}\psi \,
\partial_{zz}\psi -( \partial_{rz} \psi)^2>0$, it is a local minimum
or maximum, which corresponds to the magnetic axis ($\psi_o$), and if
it has $\partial_{rr}\psi \, \partial_{zz}\psi -( \partial_{rz}
\psi)^2<0$, it is a saddle point corresponding to the separatrix
($\psi_X$). We comment that the above algorithm does not guarantee to find a global minimum for a general solution.
However, for a physical MHD equilibrium, the solution has a well-defined local minimum problem as we shall see in the numerical section,
and thus the above algorithm works well.
 In the case when multiple local minimum/maximum points
are found, we determine the solution is invalid.  When that happens in
the free-boundary solver, it means the solver has diverged and the
iteration will stop.  In the case when multiple saddle points are
found, if the poloidal flux is convex inside the plasma region, we
choose the smallest (or largest if $\psi$ is concave) value among all
the saddle points.  Note the algorithm described above is a critical
piece in the free-boundary solver and the subroutine has to be called
routinely in the solver.
As long as all the saddle points and magnetic axis are successfully located, this ensures an accurate evaluation of the source term.
Therefore, the number of candidate points (10 in our case) has no impact on the accuracy.

Given the plasma pressure profile and the toroidal field function
profile as functions of the poloidal flux on a set of known points, finding a
proper interpolation method to describe the source term in the
Grad-Shafranov equation is very important.  We note however, at least
on equispaced case, polynomial interpolation is not proper, no matter
what form we use, due to Runge Phenomenon of high order. It is
reasonable to shift to rational interpolation
\cite{berrut1988rational, berrut2005recent, floater2007barycentric,
  schneider1986some}.  In our work, the barycentric rational
interpolation (a 4th-order rational interpolation) is employed to
evaluate the plasma pressure and the toroidal field function of the
poloidal flux. This should allow highly accurate representation of the source term
on grid points from numerical data profiles of $p(\bar{\psi})$ and $g(\bar{\psi}).$

\subsection{Determine  $\psi_b$ on $\partial\Omega$}\label{sec:bc}
The approach to determine the boundary condition value from the external coil source is well known~\cite{jardin2010computational, faugeras2017fem, jeon2015development}.
We give a brief discussion on two choices we have experimented.
Note that there is a known Green's function for the toroidal
elliptic operator $\Delta^*$ 
\begin{displaymath}
   G(\vect{R};\vect{R}')=\frac{1}{2\pi}\frac{\sqrt{RR'}}{k}[(2-k^2)K(k)-2E(k)],
\end{displaymath}
where $K(k)$ and $E(k)$ are complete elliptic integrals of the first
and second kind~\cite{polyanin2006handbook} and $k$ is given by
\begin{displaymath}
  k^2=\frac{4RR'}{(R+R')^2+(Z-Z')^2}.
\end{displaymath}

A straightforward approach based on the Green's third identity defines the boundary value as
\begin{equation}
\label{eqn:volumeInt}
  \psi_b(R',Z')=-\int_\Omega \mu_0G(R,Z;R',Z') \, \tilde J_{\phi}(R,Z) \, drdz-\sum_{i=1}^{n_c}\mu_0G(R_i^c,Z_i^c;R',Z')I_i,
\end{equation}
where the source is given by $\tilde J_{\phi}(R,Z)=J_{\phi}(R,\psi^{\rm old})$.
A more efficient approach that only computes a line integral was proposed in~\cite{hagenow1975}.  
It involves solving the elliptic problem with the same source term but a homogenous boundary condition,
\begin{equation}\label{eq:zerobdproblem}
\begin{aligned}
& \Delta^* U=\mu_0R \, \tilde J_{\phi}(R,Z), &&(R,Z)\in \Omega, \\
& U=0, && (R,Z)\in \partial\Omega.
\end{aligned}
\end{equation}
Integrating the following identity
\begin{equation*}
    \nabla \cdot [U\frac{1}{R^2}\nabla G(\xv;\xv')]-\nabla \cdot [G(\xv; \xv')\frac{1}{R^2}\nabla U]
    =\frac{1}{R^2}U\Delta^*G(\xv;\xv')-\frac{1}{R^2}G(\xv;\xv')\Delta^*U
\end{equation*}
 over the computational domain $\Omega$  leads to the formulations  for both boundary and interior points, 
 \begin{align}\label{eq:bdexpression}
  \psi_b(R',Z') &=-\int_{\partial \Omega}\frac{dl}{R}G(R,Z;R',Z')\frac{\partial U}{\partial n}-\sum_{i=1}^{n_c}\mu_0G(R_i^c,Z_i^c;R',Z')I_i, \\
\label{eq:interiorexpression}
  \psi_{interior}(R',Z')&=-\int_{\partial \Omega}\frac{dl}{R}G(R,Z;R',Z')\frac{\partial U}{\partial n}-\sum_{i=1}^{n_c}\mu_0G(R_i^c,Z_i^c;R',Z')I_i+U(R',Z').
 \end{align}

The second approach comes at the cost of inverting an elliptic operator,
   which only accounts for an extra cost of $O(N\log N)$ thanks to the
   algebraic multigrid (AMG) preconditioner we used.
   Thus, we typically choose the second approach in our solver for its efficiency.
   The second
   approach can also provide an initial guess on the entire domain
   through~\cref{eq:interiorexpression} in the initialization stage,
   which is useful to provide a better initial guess and speed up the
   convergence.  One issue is that the Green's function is singular
   when $\vect{R}'=\vect{R}$, which is more problematic in the line
   integration of the second approach.  In the implementation, when
   $\vect{R}'$ is too close to $\vect{R}$, we perturb $\vect{R}'$ with a
   small distance of $\epsilon$, which typically happens in one or two
   cells.  It appears to be sufficient to produce a smooth boundary
   condition.  
   
Note that when Newton's method is applied to the free-boundary problem,
\cref{eqn:volumeInt} is more suitable to use in its nonlinear residual. 
This, however, results into a globally coupled system due to the global operator in~\cref{eqn:volumeInt}.
The Jacobian matrix assembly and its inversion are thus expensive. 
Moreover, AMG typically cannot solve such a system involving non-local boundary conditions.
All those facts contribute to the challenge of employing Newton's method here.
To resolve the issue due to $\psi_b$, Ref.~\cite{heumann2015quasi}
proposed to extend the computational domain to include all the external currents. 
This results into a much less banded system, which is easier to assemble,  and the success of Newton's method but at the cost of solving on a domain much larger than that in the current work.

\subsection{Determine coil currents}
The last important piece in the algorithm is to determine the coil
current after a solution $\psi$ is given.  It is important to
adjust the coil current based on the solution so that the resulting
plasma has the intended position and shape. We therefore optimize the current density through a
minimization problem based on some constraints imposed through the
targeted plasma shape.  Given $n_c$ coils in the device and the
desired plasma region $\mathcal{P}_{\rm target}$, we first select $n_p$
points along its boundary $\partial\mathcal{P}_{\rm target}$ such that
the points are evenly spaced and $n_p>n_c$.
\Cref{eq:interiorexpression} provides a good approximation to $\psi$
at those points for any current density.  The search algorithm on the
other hand provides the current separatrix value $\psi_X$.  Minimizing
the $\psi$ values at those points with respect to $\psi_X$ gives a
minimization problem to find the required coil currents,
\begin{equation}\label{eq:leastsquare}
\begin{aligned}
\mathbf{I}_c = 
 \argmin_{\mathbf{I}_c}& \bigg\{ \gamma\sum_{k=1}^{n_c} {I_k}^2  + \\ 
 &  \sum^{n_p}_{j=1} \Big[-\int_{\partial \Omega}\frac{dl}{R}G(R,Z;R'_j,Z'_j)\frac{\partial U}
 {\partial     n}-\sum_{i=1}^{n_c}\mu_0G(R_i^c,Z_i^c;R'_j,Z'_j)I_i +U(R'_j,Z'_j) -\psi_{X}\Big]^2 \bigg\},
 \end{aligned}
\end{equation}
where $\mathbf{I}_c$ is a vector of size $n_c$, consisting of the current values $I_i$, and the first term is a penalty term to prevent the system from becoming ill-posed.  

We found that the resulting current density is sensitive to the choice
of $\gamma$.  Part of the reason is that the problem setup is based on
practical units and the density values can be dramatically different
from each other in different coils.  To address that, we found the Generalized Cross validation (GCV) approach  in~\cite{golub1979generalized} provides a good estimate of $\gamma$.
To utilize GCV, the minimization problem~\cref{eq:leastsquare} is written in the form of 
\begin{align*}
\mathbf{I}_c= 
 \argmin_{\mathbf{I}_c} \bigg ( \gamma \, \| \mathbf{I}_c \|^2 +   \frac 1 {n_p} \,  \| X \, \mathbf{I}_c - \mathbf{y} \|^2 \bigg),
 \end{align*}
 where \cref{eq:leastsquare} indicates $X$ is  a matrix of size $n_p\times n_c$, and $\mathbf{y}$ is a vector of size $n_p$.
 The GVC approach suggests the best $\gamma$ should be computed as the minimizer of $V(\gamma)$, given by 
 \[
 V(\gamma) =  \frac 1 {n_p} \| (I - A(\gamma)) \mathbf{y}\|^2 \bigg/ \bigg[  \frac 1 {n_p} \, {\rm Trace} (I - A(\gamma))\bigg]^2,
 \]
 where $A(\gamma) \equiv X(X^T X + n_p \gamma  I)^{-1} X^T$ and $I$ is an identity matrix of size $n_p\times n_p$.
 In the initial experiment, we implement a simple search algorithm to locate the optimized value of $\gamma$ in each iteration for determining the current. 
 We found that the optimized $\gamma$ falls within the range of $[10^{-16}, 10^{-13}]$ for the practical ITER equilibrium problem.
 Therefore, to avoid recalculating $\gamma$ in each iteration, we fix $\gamma =10^{-15} $ in all of our numerical examples, which is sufficient to provide less sensitive current values
in the final solver. 
The current values of the equilibrium are found in the range of  $[10^3, 10^7]$, and thus the regularization term is not negligible with the given $\gamma$ value.

\subsection{Full algorithm}\label{sec:algorithm}
With all the important pieces laid out in the previous sections, a
full algorithm for the free-boundary Grad-Shafranov equation will be
described.  In particular, two versions of the algorithm will be
discussed and an additional acceleration technique will be adopted to
improve the robustness and efficiency of the solver.

\tikzstyle{decision} = [ellipse, draw, fill=white!20, 
    text width=2.4cm,  minimum height=2em, text badly centered, node distance=1.2cm, inner sep=0pt]
\tikzstyle{block1} = [rectangle, draw, fill=blue!20, 
    text width=12em, text centered, rounded corners, minimum height=2em]
\tikzstyle{block2} = [rectangle, draw, fill=blue!20, 
    text width=9em, text centered, rounded corners, minimum height=2em]
   \tikzstyle{block3} = [rectangle, draw, fill=blue!20, 
    text width=2.7cm, text centered, rounded corners, minimum height=2em]
  \tikzstyle{block4} = [rectangle, draw, fill=blue!20, 
    text width=3em, text centered, rounded corners, minimum height=2em]
\tikzstyle{line} = [draw, -latex']
\tikzstyle{cloud1} = [draw, rectangle,fill=blue!20, rounded corners, node distance=3cm,
    minimum height=2em]
\tikzstyle{cloud2} = [draw, rectangle,fill=blue!20, rounded corners, node distance=4.3cm,
    minimum height=2em]
\tikzstyle{cloud3} = [draw, rectangle,fill=blue!20, rounded corners,node distance=5.6cm,
    minimum height=2em]

\begin{figure}
  \centering
  \begin{subfigure}{.6\textwidth}
    \centering
{\footnotesize
\begin{tikzpicture}[node distance = 1.1cm, auto]
 \useasboundingbox (-4,-8.3) rectangle (1,.4);
    \node [block1] (init) {Given initial equilibrium current $\psi_0$};
    \node [block1, below of=init] (currentdensities) {Calculate coil currents $I_i^k$};
    \node [block3, below of=currentdensities] (Boundary) {Determine $\psi_b^m$};
    \node [block3, below of=Boundary] (GS) {Solve $\Delta^* \psi^n=\mathcal{S}(\psi^{n-1})$};
    \node [decision, below of=GS] (innerloop) {$\psi^n$ Converged?};
    \node [cloud1, left of=innerloop] (Itersininner) {$n=n+1$};
    \node [decision, below of=innerloop] (middleloop) {$\psi_b^m$ Converged?};
    \node [cloud2, left of=middleloop] (Itersinmid) {$m=m+1$};
    \node [decision, below of=middleloop] (outerloop) {$I_i^k$ Converged?};
     \node [cloud3, left of=outerloop] (Itersinouter) {$k=k+1$};
     \node [block4, below of=outerloop] (stop) {Stop};
    \path [line] (init) --node [, color=black] {Initialize($k=0$)} (currentdensities);
    \path [line] (currentdensities) -- (Boundary);
    \path [line] (Boundary) -- (GS);
    \path [line] (GS) -- (innerloop);
    \path [line] (innerloop) -- node [, color=black] {Yes} (middleloop);
     \path [line] (middleloop) --node [, color=black] {Yes} (outerloop);
    \path [line] (outerloop) --node [, color=black] {Yes}  (stop);
    \path [line] (innerloop) -- node [, color=black] {No}(Itersininner);
    \path [line] (middleloop) --node [, color=black] {No} (Itersinmid);
   \path [line] (outerloop) --node [, color=black] {No} (Itersinouter);
   \path [line] (Itersininner) |-(GS);
    \path [line] (Itersinmid)|-(Boundary);
   \path [line] (Itersinouter) |-(currentdensities);
\end{tikzpicture}
}
   \caption{Three-loop algorithm}
      \label{fig:loop3}
  \end{subfigure}%
  \begin{subfigure}{.4\textwidth}
  \centering
  {\footnotesize
\begin{tikzpicture}[node distance = 1.4cm, auto]
 \useasboundingbox (-3,-7.7) rectangle (1,1);
    \node [block1] (init) {Given initial equilibrium current $\psi_0$};
    \node [block1, below of=init] (Boundary) {Calculate coil currents $I_i^m$, determine boundary value $\psi_b^m$};
    \node [block3, below of=Boundary] (GS) {Solve $\Delta^* \psi^n=\mathcal{S}(\psi^{n-1})$};
    \node [decision, below of=GS] (innerloop) {$\psi^n$ Converged?};
    \node [cloud1, left of=innerloop] (Itersininner) {$n=n+1$};
    \node [decision, below of=innerloop] (middleloop) {$\psi_b^m$ Converged?};
    \node [cloud2, left of=middleloop] (Itersinmid) {$m=m+1$};
     \node [block4, below of=middleloop] (stop) {Stop};
    \path [line] (init) --node [, color=black] {Initialize ($m=0$)} (Boundary);
    \path [line] (Boundary) -- (GS);
    \path [line] (GS) -- (innerloop);
    \path [line] (innerloop) -- node [, color=black] {Yes} (middleloop);
     \path [line] (middleloop) --node [, color=black] {Yes} (stop);
    \path [line] (innerloop) -- node [, color=black] {No}(Itersininner);
    \path [line] (middleloop) --node [, color=black] {No} (Itersinmid);
   \path [line] (Itersininner) |-(GS);
    \path [line] (Itersinmid)|-(Boundary);
\end{tikzpicture}
}
   \caption{Two-loop algorithm}
      \label{fig:loop2}
  \end{subfigure}
  \caption{Flowcharts for the free-boundary Grad-Shafranov solvers. }
  \label{fig:loopDiagrams}
\end{figure}
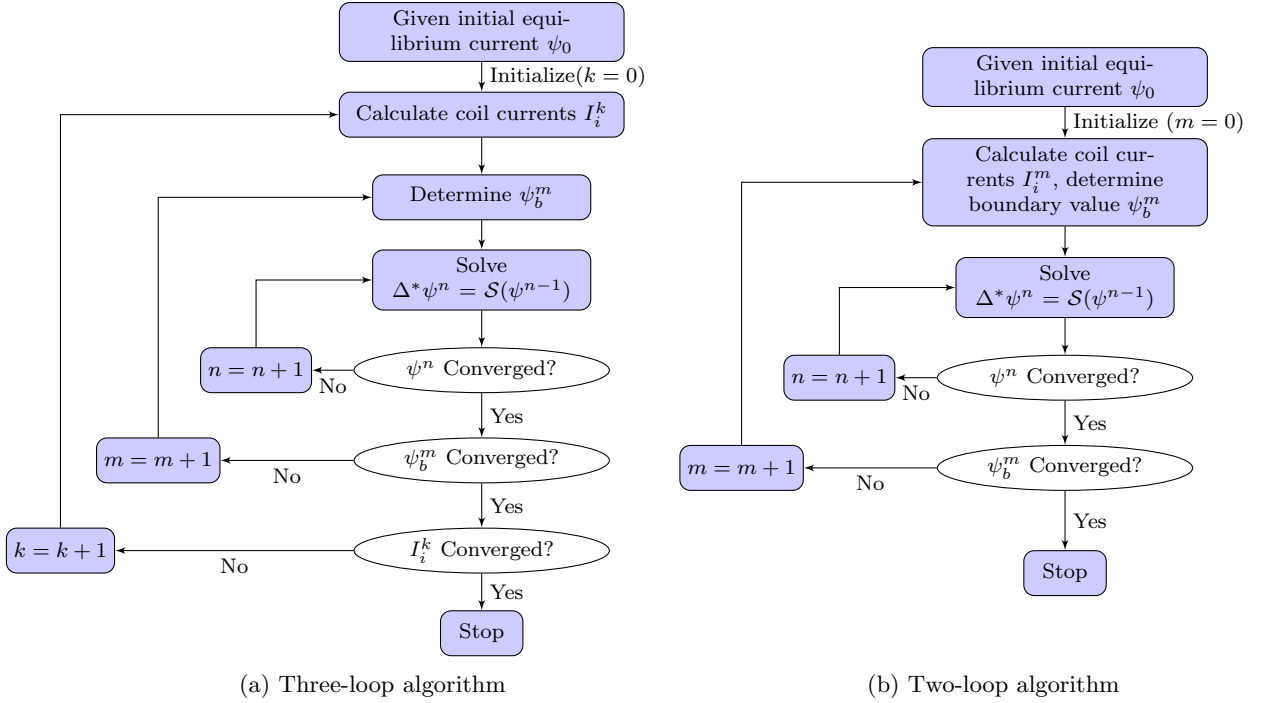

\paragraph*{Three-loop algorithm.}
The first version of the algorithm would consist of three loops as the
main steps~(\Cref{fig:loop3}): the most inner loop inverts the operator $\Delta^*$ with a
given boundary condition $\psi_b$; the outer loop will adjust $\psi_b$
until some convergence threshold is reached for some fixed current
density $\{I_i\}$~(see \Cref{fig:intro} for the current locations); then the most outer loop would adjust the coil
currents such that the magnetic separatrix will align with the control
points.  Although this idea appears to be plausible and in fact it is
the first algorithm we implemented, we found that this version of the
algorithm is sensitive to a small perturbation of the current, which
could easily result into the failure of convergence of the two inner
loops.  One of the fundamental issues in this approach is that there
is an underlying assumption that for any given current density, there
exists a corresponding equilibrium, despite not satisfying the plasma
domain constraint.  This, however, may not be true in general, and
even if there exists such an equilibrium, its $g$ and $p$ profiles
could be very different from the profiles being used, and therefore
the evaluated source term could be very inaccurate, which could
contribute to the divergence of the solver.  Therefore, we do not
pursue further using this version of the algorithm in the current
work.

\paragraph*{Two-loop algorithm~\cite{johnson1979numerical}.}
The second version of the algorithm consists of two loops as the main
steps~(\Cref{fig:loop2}):
the inner loop still inverts the operator $\Delta^*$ with a given
boundary condition $\psi_b$ while the outer loop simultaneously
updates $\psi_b$ and $\{I_i\}$ until some threshold with respect to
$\psi_b$ is reached. This version turns out to be much more robust
than the previous version, which is the main algorithm we used in the
free-boundary solver.  Its robustness can be further improved by the
acceleration techniques discussed below.

\paragraph*{Aitken's acceleration}
As previously discussed, the general procedure of the algorithm is a
Picard iteration, and its drawback is its slow convergence.  A common
approach to accelerate its convergence is through under-relaxation,
i.e.,
$$
  \psi^{new}=(1-\alpha)\psi^{old}+\alpha\tilde\psi^{{new}}.
$$
where $\alpha$ is the under-relaxation coefficient, typically between 0 and 1. 
The optimal value of $\alpha$, however, is problem dependent and for nonlinear problems it can vary during the iterations. The optimal relaxation value can be approximated through the well known Aitken's acceleration~\cite{mok2001accelerated}. 
The version we adopted is given by
\begin{equation}\label{eq:Aitkenpart3}
\begin{aligned}
D \psi^{n+1}&=\psi^{n}-\tilde \psi^{{n+1}},\\
\lambda^{n+1}&=\lambda^{n}+(\lambda^n-1)\frac{(D \psi^{n}-D \psi^{n+1})^T D \psi^{n+1}}{\|D \psi^{n}-D \psi^{n+1}\|^2},\\
\alpha &=1-\max ( \min(\lambda^{n+1}, \lambda_{\max}),\lambda_{\min}),\\
 \psi^{n+1}&=(1-\alpha)\psi^{n}+\alpha\tilde\psi^{n+1},
\end{aligned}
\end{equation}
where $\tilde\psi^{n+1}$ is the solution directly solved from the Picard iteration, and $\lambda_{\min}=0$ and $\lambda_{\max}=0.95$ are typically chosen. 
Ref.~\cite{borazjani2008curvilinear} suggests a preset value of $\lambda= 0.3$ in the first iteration, which is adopted in the current work. 
Note that the technique fits into the current framework very well due to its simplicity and computational expedience,
and we use it to accelerate both the inner and outer loops.
It is found that the Aitken's acceleration greatly improves the convergence of the full algorithm,
which will be demonstrated in the numerical section. 
We note however that there are other techniques such as Anderson's acceleration and nonlinear GMRES available for the Picard iteration. 
We plan to explore more acceleration techniques in future work.

Finally, we summarize the full algorithm in \cref{alg:2loops}, which
consists of two main loops, Aitken's accelerations and all the
important steps described in the previous sections.  Note that Step 2
corresponds to all the steps related to the cut-cell algorithm, which
will be addressed carefully in the next section. In the nonlinear solver, we
 set the relative difference of the boundary values $\epsilon_{out}= \num{1e-4}$  as the convergence criterion for the outer loop  and
 the relative difference of the solution values $\epsilon_{in}=\num{1e-3}$ as  the convergence criterion for the inner loop. For Krylov linear solvers in PETSc, we set the relative tolerance of convergence as
  $\num{1e-5}$, and the absolute tolerance of convergence as
  $\num{1e-8}$.

{
\begin{algorithm}
\caption{Free-boundary Grad-Shafranov solver}
\label{alg:2loops}
\begin{algorithmic}[1]
\small
\STATE{Choose an irregular domain $\Omega$ to reformulate the free boundary problem \cref{eq:freebdproblempart1} and \cref{eq:freebdproblempart2}}
\STATE{Construct a level set function $\tilde{\phi}$ and prepare a cut-cell mesh}
\STATE{Given an initial equilibrium $\psi^0$, solve \Cref{eq:zerobdproblem} for $U(R,Z)$}
\STATE{Calculate coil currents $I_i^0$ according to \Cref{eq:leastsquare}}
\STATE{Determine boundary value $\psi_b^0$ according to \Cref{eq:bdexpression}}
\smallskip
\STATE{\tt //Outer loop}
\FOR{$m =0~{\rm to }~m_{\max}$}
\smallskip
\STATE{\tt //Inner loop}
  \FOR{$n =0~{\rm to }~n_{\max}$} 
       \STATE{Search $\psi_o$ and $\psi_X$ in $\psi^n$}
   \STATE{Solve  $\Delta_h^* \tilde\psi^{{n+1}}= \mathcal{S}(\psi^n)$} with the boundary condition $\psi_b^m$ 
     \IF {n=0} 
     \STATE{$\lambda_{in}^0$=0.3}
     \STATE{$\alpha_{in}=1-\lambda_{in}^0$} 
     \ELSE
     \STATE{$\lambda_{in}^{n+1}=\lambda_{in}^{n}+(\lambda_{in}^n-1)\frac{(D \psi^{n}-D \psi^{n+1})^T D \psi^{n+1}}{\|D \psi^{n}-D \psi^{n+1}\|^2}$ with $D \psi^{n+1}=\psi^{n}-\tilde \psi^{{n+1}}$}
      \STATE{$\alpha_{in} =1- \max ( \min(\lambda_{in}^{n+1}, \lambda_{\max}),\lambda_{\min})$}
     \ENDIF
    \STATE{$\psi^{n+1}=(1-\alpha_{in})\psi^n+\alpha_{in}\tilde\psi^{{n+1}}$}
    \STATE{Check convergence by $\Delta \psi= \frac{\parallel \psi^{n+1}-\psi^{n}  \parallel_{\infty, \Omega}}{\parallel \psi^{0} \parallel_{\infty, \Omega}}$}
    \IF{($\Delta \psi <\epsilon_{in}$)}
    \STATE{\bf break}
     \ENDIF
     \ENDFOR
     \smallskip
     \STATE{Given $\psi^{n+1}$, solve \Cref{eq:zerobdproblem} for $U(R,Z)$}
     \STATE{Update coil currents $I_i^m$ according to \Cref{eq:leastsquare}}
     \STATE{Determine boundary value $\tilde\psi_b^m$  according to \Cref{eq:bdexpression}}
      \IF {m=0}
      \STATE{$\lambda_{out}^0$=0.3}
     \STATE{$\alpha_{out}=1-\lambda_{out}^0$}
     \ELSE
     \STATE{$\lambda_{out}^{m+1}=\lambda_{out}^{m}+(\lambda_{out}^m-1)\frac{(D \psi_b^{m}-D \psi_b^{m+1})^T D \psi_b^{m+1}}{\|D \psi_b^{m}-D \psi_b^{m+1}\|^2}$ with $D \psi_b^{m+1}=\psi_b^{m}-\tilde \psi_b^{{m+1}}$}
     \STATE{$\alpha_{out} =1- \max ( \min(\lambda_{out}^{m+1}, \lambda_{\max}),\lambda_{\min})$}
     \ENDIF
     \STATE{$\psi_b^{m+1}=(1-\alpha_{out})\psi^m_b+\alpha_{out}\tilde\psi_b^{{m+1}}$}
     \STATE{Check convergence by $\Delta \psi_b= \frac{ \parallel \psi^{m+1}_b-\psi^{m}_b  \parallel_{\infty, \partial\Omega}}{\parallel \psi^{0}_b  \parallel_{\infty, \partial\Omega}}$}
     \IF{$\Delta \psi_b < \epsilon_{out}$}
     \STATE {\bf break}
    \ENDIF
\ENDFOR
\end{algorithmic} \end{algorithm}
}

\section{Parallel implementation}
\label{sec:implementation}
 Some aspects of the implementation of the free boundary
 Grad-Shafranov solver with the cut-cell algorithm will be discussed
 in this section.
 The free boundary
 Grad-Shafranov solver described in this work is implemented in
 parallel under the PETSc framework~\cite{balay2019petsc} using a
 standard domain decomposition approach.  All the vectors, arrays and
 matrices use the parallel distributed data structure provided by
 PETSc and its communication between sub-domains is based on the
 message-passing interface (MPI).  All the linear and nonlinear
 solvers described in this work are implemented through PETSc and the
 elliptic operator, $\Delta_h^*$, is preconditioned with algebraic
 multigrid preconditioners to improve its efficiency and scalability.

 The data structure is based on a Cartesian structured mesh (DMDA)
 provided by PETSc.  The solution is stored as grid points while its
 control volume is straightforwardly implied.  This is a common
 approach to adapt a finite volume algorithm in a finite difference
 code, see~\cite{brown1997overture} for instance.  To develop a
 cut-cell algorithm based on DMDA, all the grid points are
 distinguished into three categories of active points, inactive points
 and cut-cell points, based on its underlying control volume.  If the
 control volume has the full cell size, then it is an active point,
 otherwise it is either a cut-cell point (if the volume fraction is
 between 0 and 1) or an inactive point (if the volume fraction is 0).
 Note that the approach distinguishing active and inactive grid points
 is commonly used in overlapping grids, where some of the grid points
 are not involved in the discretization, see~\cite{banks2017stable,
   banks2018stable} for instance.  Since the scheme is based on a
 five-point stencil, the discrete operator on the active points is
 well-defined, and meanwhile a redundant equation is imposed on all
 the inactive points when discretizing the elliptic operator on those
 points (simply solving $a \psi_{\rm inactive} = 0$ where $a$ is a
 large constant).  The main difference between the free-boundary
 solvers on a standard Cartesian mesh and on cut cells thus lies in
 the additional treatment related to the cut-cell points.  Involving a
 redundant equation into the discretization leads to an operator of
 the same size as the standard operator on a Cartesian grid.  All the
 approaches described above guarantee many PETSc functions for DMDA
 can be directly used in the cut-cell algorithm, which eases some data
 structure allocations and implementations of linear or nonlinear
 operators.  It also provides a rather straightforward path for the
 future improvement through geometry multigrid.

 To further ease the treatment related to cut-cells, 
several 2D arrays are precomputed to store the geometry information and the information related to level set function.  
For example,  volume fractions and area fractions of each cell, apertures of each edge and midpoint values of each cut edge and interpolation points are precomputed before performing the Picard iterations, 
which greatly improves the solver efficiency.

There are some details associated with the main algorithm worthwhile to mention. 
The search algorithm in parallel needs first identify candidate points through performing search on the sub-domain
and then gather all the candidate points into the root processor to finalize the locations of $\psi_x$ and $\psi_o$. 
Some MPI collective communications are necessary to successfully locate those points. 
Another important aspect is to perform an efficient line integral along the boundary.
The efficiency of computing the line integral can be improved by precomputing certain weights associated with the Green's function and the quadratures. 
For instance, the terms changing in \Cref{eq:bdexpression}  in each iteration are $\frac{\partial U}{\partial n}$ and $I_i$, and we therefore precompute the weights 
associated with $\frac{dl}{R}G(R,Z;R',Z')$ and $\mu_0G(R_i^c,Z_i^c;R',Z')$ and store them in two distributed dense matrices. 
Parallel matrix-vector operations are called to evaluate the boundary values $\psi_b$ efficiently during the iterations.

{
\newcommand{\logLogSlopeTriangle}[5]
{

    \pgfplotsextra
    {
        \pgfkeysgetvalue{/pgfplots/xmin}{\xmin}
        \pgfkeysgetvalue{/pgfplots/xmax}{\xmax}
        \pgfkeysgetvalue{/pgfplots/ymin}{\ymin}
        \pgfkeysgetvalue{/pgfplots/ymax}{\ymax}

        \pgfmathsetmacro{\xArel}{#1}
        \pgfmathsetmacro{\yArel}{#3}
        \pgfmathsetmacro{\xBrel}{#1-#2}
        \pgfmathsetmacro{\yBrel}{\yArel}
        \pgfmathsetmacro{\xCrel}{\xArel}

        \pgfmathsetmacro{\lnxB}{\xmin*(1-(#1-#2))+\xmax*(#1-#2)} 
        \pgfmathsetmacro{\lnxA}{\xmin*(1-#1)+\xmax*#1} 
        \pgfmathsetmacro{\lnyA}{\ymin*(1-#3)+\ymax*#3} 
        \pgfmathsetmacro{\lnyC}{\lnyA+#4*(\lnxA-\lnxB)}
        \pgfmathsetmacro{\yCrel}{\lnyC-\ymin)/(\ymax-\ymin)} 

        \coordinate (A) at (rel axis cs:\xArel,\yArel);
        \coordinate (B) at (rel axis cs:\xBrel,\yBrel);
        \coordinate (C) at (rel axis cs:\xCrel,\yCrel);

        \draw[#5]   (A)-- node[pos=0.5,anchor=north] {1}
                    (B)-- 
                    (C)-- node[pos=0.5,anchor=west] {1}
                    cycle;
    }
}

\begin{figure}[htb]
\begin{center}
    \begin{tikzpicture}[scale=.8]
    \useasboundingbox (0, 0.0) rectangle (8.,7);  
\begin{scope}[xshift=1cm, yshift=1cm]
\begin{loglogaxis}[xlabel=\# of CPUs,ylabel=Time, grid=both,
legend entries={Free-boundary solver, Linear scaling}
]
\addplot coordinates {
	(1, 897.68)
    (2, 535.44)
	(4, 315.6)
    (8, 212)
    (16, 125)
    (32, 84.5)
};
\addplot[black,mark=*,mark options={fill=black}] coordinates {
	(1, 897.68)
    (2, 449)
	(4, 224.5)
    (8, 112.2)
    (16, 56.1)
    (32, 28)
};
\end{loglogaxis}
\end{scope}
\end{tikzpicture}
\end{center}
\caption{Strong scaling result of the free-boundary cut-cell Grad-Shafranov solver.
The computational times of one and up to 32 CPUs are presented. A good parallel scaling up to 32 processors is observed.
A cut-cell mesh based on a Cartesian mesh of $512 \times 1024$ is used in the strong scaling study.
Details of the problem setup can be found in the last free-boundary example in the numerical section.
}
  \label{fig:icDD}
\end{figure}
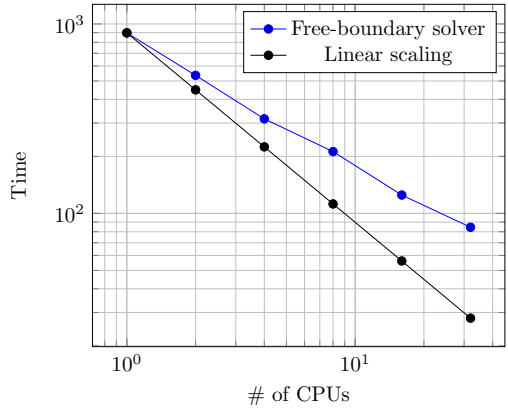
}

Finally, to verify the performance  of the full free-boundary Grad-Shafranov solver, a strong scaling study is performed.
In the study, a cut-cell mesh based on a Cartesian mesh of size $512 \times 1024$ is used. 
We verify the performance of the solver using different numbers of processors up to 32. 
\Cref{fig:icDD} shows the strong scaling result and a good scaling is observed.
The corresponding parallel efficiency for 2, 4, 8, 16{,} and 32 processors are 83.9\%, 71.1\%, 52.9\%, 44.9\%{,} and 33.2\%.
Note that since the algorithm is based on Picard iteration fundamentally, we do not expect such an algorithm would be scalable when the number of processors increase to some large number. 
Nevertheless, the performance of the algorithm is decent and 32 processors are likely sufficient for a 2D steady-state problem. 
In the scaling study, the full algorithm described in~\cref{alg:2loops} is used and  all the elliptic operators use an algebraic multigrid preconditioner provided by PETSc. 
Details of the problem setup can be found in the free-boundary cut-cell example in the numerical section later.

\section {Numerical results}
\label{sec:numericaltests}
The performances of the cut-cell algorithm, the free-boundary
Grad-Shafranov solver and Aitken's acceleration are demonstrated
through several numerical tests. We start with the accuracy test of
the cut-cell solver through two examples: Grad-Shafranov equation with
a linear source term and Grad-Shafranov equation with a nonlinear
source term.  Numerical examples of the free-boundary problem in a
rectangular domain and in the limiter-bounded domain $ \mathcal{L}$
are presented separately to show the performance of the free-boundary
Grad-Shafranov solver and cut-cell algorithm.  Finally, we compare the
performance of Picard iteration with different under-relaxation
coefficients and Aitken's acceleration.

\subsection{Convergence study for the fixed-boundary cut-cell solver}
\label{sec: convergencestudy}
In this section, two numerical examples are presented to demonstrate
the accuracy of the cut-cell algorithm.  Here we use two
fixed-boundary problems with cut-cells to verify its accuracy. These
cases are chosen for the known analytical solutions.  Note that
compared with a free-boundary solver, the fixed-boundary solver do not
have the surface integral and the barycentric interpolation, both of
which can be verified easily through standalone tests.

In the first example, we consider the linear Soloviev profiles from
\cite{cerfon2010one}. This test is used to compare the accuracy of the
solution $\psi$ solved by our fixed boundary Grad-Shafranov solver and
the analytical solution. The second example considers a manufactured
solution with a nonlinear source. We use it to demonstrate the
accuracy of our scheme in nonlinear problems.

\subsubsection{Linear case}
\label{sec:lineartest}
Consider a linear Grad-Shafranov equation of
$$
\Delta^* \psi=-\frac{1}{R}\frac{\partial
  \psi}{\partial R}+\frac{\partial ^2 \psi}{\partial
  R^2}+\frac{\partial ^2 \psi}{\partial
  Z^2}=R^2,
$$
which has an exact solution in the form of 
$$
\psi(R,Z)=\frac{R^4}{8}+D_1+D_2R^2+D_3(R^4-4R^2Z^2),
$$
where the parameters $D_1$, $D_2${,} and $D_3$ are determined so that the contour $\phi=0$ represents a reasonable plasma cross section. 

\begin{table}[h!]
\caption{\label{tab:21}$L_1$-errors, $L_2$-errors, $L_\infty$-errors  and corresponding convergence rates of $\psi$ in the linear accuracy test using the cut-cell mesh.
The grid points in the meshes represent all the grid points in a Cartesian grid including both active and inactive points.
}
\begin{center}
\begin{tabular}{ |c||c|c|c|c|c|c|  }
 \hline
 $N_x\times N_y$& $L_1$-Error&Order&$L_2$-Error&Order& $L_{\infty}$-Error &Order\\
 \hline
 \hline
 31$\times$41&2.310e-01   & &1.276e-02 &  & 1.519e-03    &  \\
 \hline
 61$\times$81&5.925e-02  & 1.96&3.361e-03&1.92 &5.513e-04   &1.46\\
 \hline
 121$\times$161&1.467e-02 &2.01&8.327e-04&2.01 &9.177e-05  &2.59 \\
 \hline
 241$\times$321&3.774e-03   & 1.96&2.155e-04&1.95 & 2.229e-05 &2.04\\
 \hline
 481$\times$641&1.086e-03  & 1.80&6.198e-05&1.80&5.422e-06    &2.04\\
 \hline
 \end{tabular}
\end{center}
\end{table}

\begin{figure}
  \centering
        \includegraphics[trim=0cm 0cm 3cm 0cm, clip=true, width=.36\textwidth]{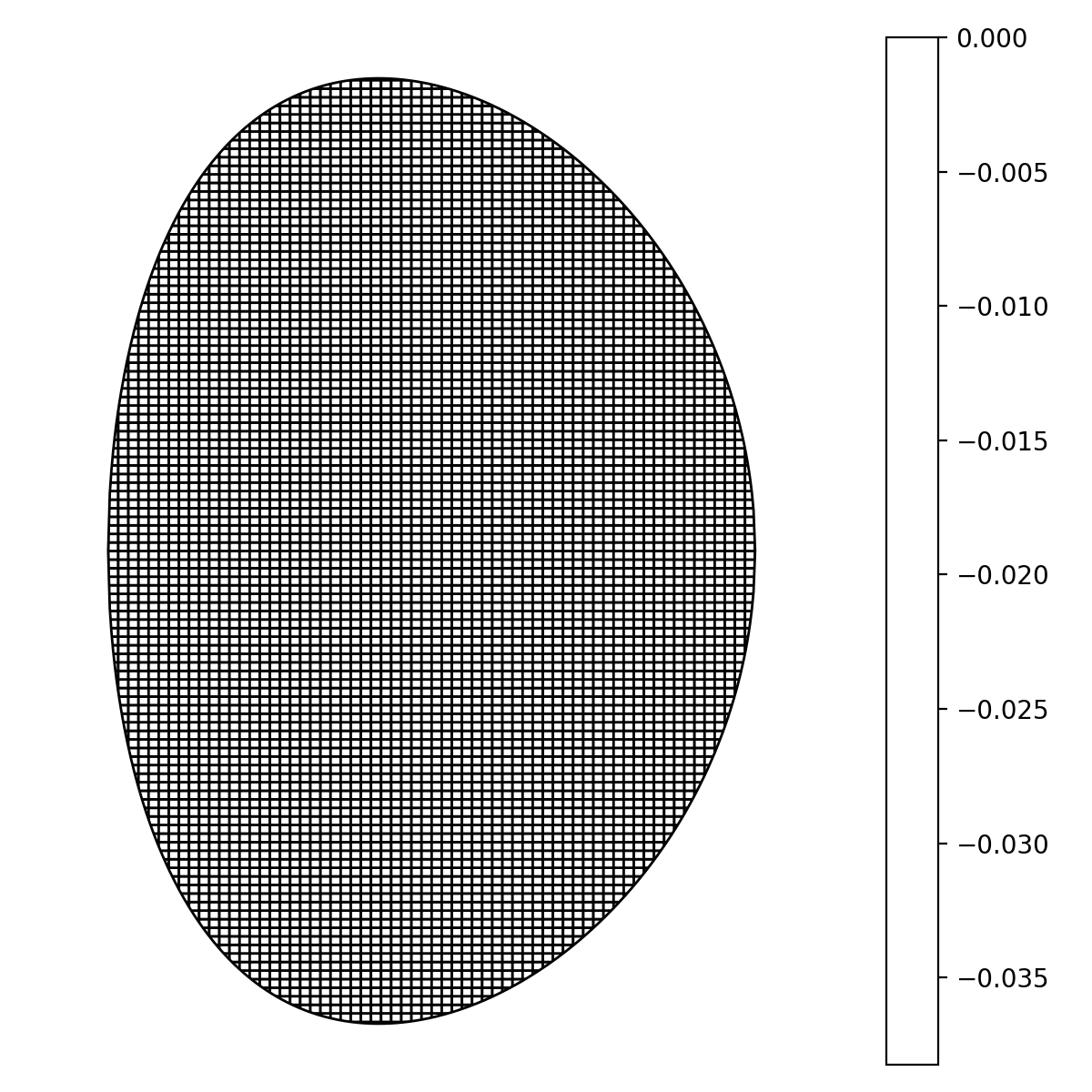}
      \includegraphics[scale=0.49]{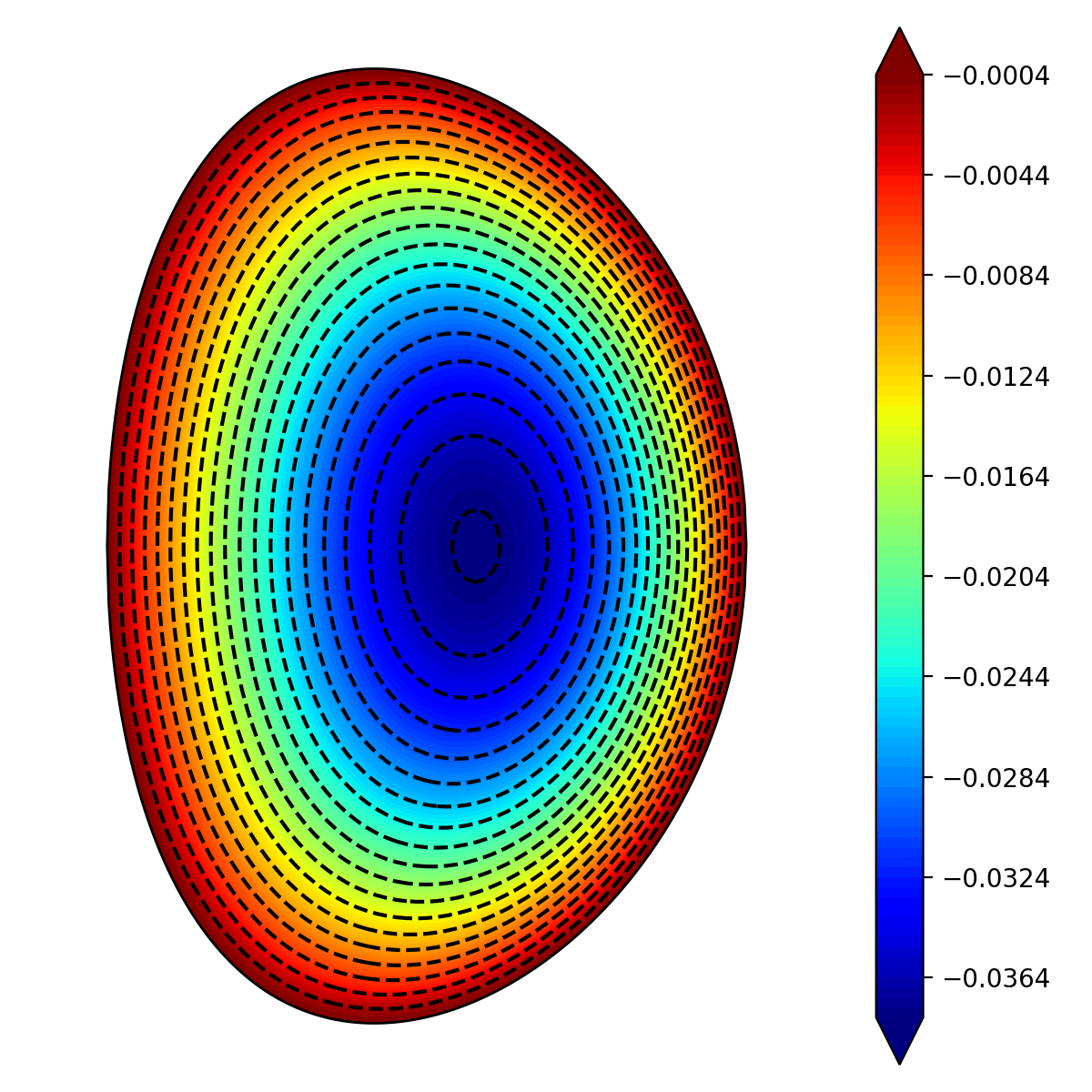}
   \caption{Linear accuracy test.  Left: cut-cell mesh on a base mesh of $121\times161$ (only active points are shown) and geometry.
   Right: corresponding numerical solution $\psi$. }
    \label{fig:cutcelllineartestnumericalsoln}
\end{figure}

Ref.~\cite{pataki2013fast} introduced three characteristic quantities describing the shape of the cross section in a magnetic confinement device: the inverse aspect ratio
$\varepsilon$, the elongation $\kappa${,} and the triangularity $\delta$. We adopt the same
manufactured solutions in the accuracy test. The parameters are
determined by the following linear system
\begin{align*}
\left[
\begin{array}{ccc}
1  & (1+\varepsilon)^2 &(1+\varepsilon)^4 \\
1  & (1- \varepsilon)^2 &(1- \varepsilon)^4 \\
1  & (1-\delta\varepsilon)^2 & (1-\delta\varepsilon)^4-4(1-\delta\varepsilon)^2\kappa^2\varepsilon^2\\
\end{array}\right]
\begin{bmatrix}
D_1\\
D_2\\
D_3
\end{bmatrix} 
&=-\frac{1}{8}\begin{bmatrix}
(1+\varepsilon)^4\\
(1-\varepsilon)^4\\
(1-\delta\varepsilon)^4
\end{bmatrix},
\end{align*}
of which the equations correspond to the boundary conditions of $\psi(1+\varepsilon, 0) = 0$, $\psi(1-\varepsilon, 0) = 0$, and $\psi(1-\delta\varepsilon, \kappa \varepsilon) = 0$, respectively.
The test is taken with the ITER-like configuration of $ \varepsilon=
0.32$, $\kappa= 1.7$, $\delta= 0.33$.  In particular, the
computational boundaries are described by Chebyshev nodes along the
$R$ direction and $Z$ coordinates are then determined by solving
$\psi=0$.  The cut-cell meshes are then accordingly generated.

Numerical errors  and corresponding convergence rates are reported in~\cref{tab:21}. The numerical
solution and an example cut-cell mesh are presented
in~\Cref{fig:cutcelllineartestnumericalsoln}.  We observed a
second-order accuracy for all the error norms.

\subsubsection{Nonlinear case}
\label{sec:nonlineartest}
We consider the same ITER geometry {as for the} linear test in the previous
section but with a nonlinear source term.  
A manufactured solution 
$$
\psi(R,Z)=\sin(K_R(R+R_0))\cos(K_ZZ),
$$
can be constructed for 
the nonlinear Grad-Shafranov equation of
$$
\Delta^* \psi=-F(R,Z,\psi),
$$
where  the source term $F(R,Z,\psi)$ is given by
\begin{align*}
F(R,Z,\psi) &=(K_R^2+K_Z^2)\psi+\frac{K_R}{R}\cos(K_R(R+R_0))\cos(K_ZZ)+R \Big[\sin^2(K_R(R+R_0))\cos^2(K_ZZ)\\
&-\psi^2+\exp(-\sin(K_R(R+R_0))\cos(K_ZZ))-\exp(-\psi)\Big],
\end{align*}
and the coefficients are $K_R=1.15\pi$, $K_Z= 1.15${,} and $R_0 = −0.5$.

\begin{table}[h!]
\caption{\label{tab:22}$L_1$-errors, $L_2$-errors, $L_\infty$-errors  and corresponding convergence rates of $\psi$ in the nonlinear accuracy test using cut-cell algorithm. The grid points in the meshes represent all the grid points in a Cartesian grid including both active and inactive points.}
\begin{center}
\begin{tabular}{ |c||c|c|c|c|c|c|  }
 \hline
 $N_x\times N_y$& $L_1$-Error&Order&$L_2$-Error&Order& $L_{\infty}$-Error &Order\\
 \hline
 \hline
 31$\times$41&1.486e+00   & &8.493e-02 &  & 1.359e-02    &  \\
 \hline
 61$\times$81&4.476e-01  & 1.73&2.584e-02&1.72 &4.480e-03   &1.60\\
 \hline
 121$\times$161&1.340e-01 &1.74&7.688e-03&1.75 &8.288e-04  &2.43 \\
 \hline
 241$\times$321&3.371e-02  & 1.99&1.912e-03&2.01& 2.154e-04 &1.94\\
 \hline
 481$\times$641&7.825e-03  & 2.11& 4.429e-04&2.11&4.148e-05    &2.38\\
 \hline
 \end{tabular}
\end{center}
\end{table}

\begin{figure}
  \centering
       \includegraphics[trim=0cm 0cm 3cm 0cm, clip=true, width=.36\textwidth]{Cutcell_lineartest_numerical_grid.png}
      \includegraphics[scale=0.49]{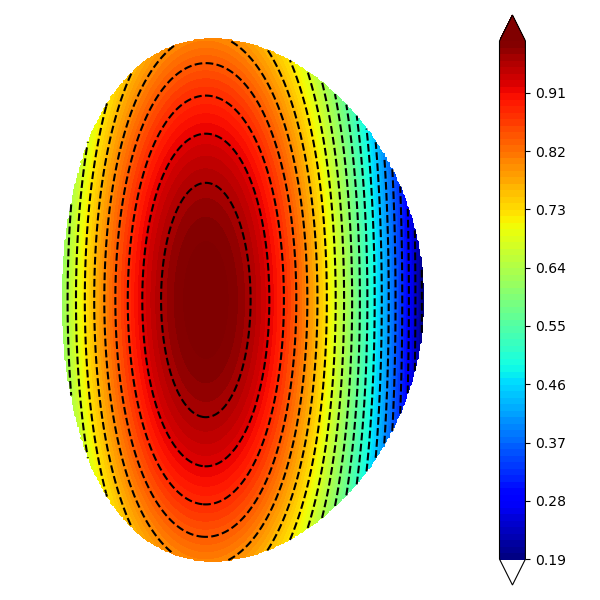}
         \caption{Nonlinear accuracy test.  Left: cut-cell mesh  on a base mesh of $121\times161$ (only active points are shown) and geometry.
   Right: corresponding numerical solution $\psi$. }
    \label{fig:cutcellnonlineartestnumericalsoln}
\end{figure}

A nonlinear solver (Picard iteration accelerated with Aitken's) is
used in this example. Numerical errors  and corresponding convergence rates are reported in~\cref{tab:22}.   The solution and an example cut-cell
mesh are presented in~\Cref{fig:cutcellnonlineartestnumericalsoln}.
Again, as expected, we observe a second-order accuracy of space
discretization.

\subsection{Examples of the free-boundary solver} 
 In this section, numerical solutions from the free-boundary Grad-Shafranov solver are presented.  One particular interest will be
 the performance of the solver in keeping the targeted shape of
 $\mathcal{P}(\psi)$, i.e., the shape control.  Numerical tests are
 scattered into the following two examples, a rectangular domain and a
 limiter-bounded domain.
 
 \begin{table}[h!]
\caption{\label{tab:11}Coil Information for the free-boundary problem. The data is based on the ITER configurations. 
The total current in the coil is set to current timing the value of turns.
 The type ``coil'' refers to a point source.
The type ``solenoid'' refers to a coil with a positive length for which the range of length is also given.
The units is meter.
}
\begin{center}
\begin{tabular}{|c||c|c|c|c|}
 \hline
 Coils & Type &$R$&$Z$&Turns\\
 \hline
 \hline
 PF1&Coil & 3.9431&7.5741 &  248.6  \\
 \hline
 PF2&Coil & 8.2851&6.5398 &  115.2  \\
  \hline
  PF3&Coil & 11.9919&3.2752 &  185.9  \\
   \hline
  PF4&Coil & 11.9630&-2.2336 &  169.9  \\
 \hline
  PF5&Coil & 8.3908&-6.7269 &  216.8  \\
 \hline
   PF6&Coil & 4.3340&-7.4665 &  459.4 \\
 \hline
 \end{tabular}
 \hspace{.2cm}
  \begin{tabular}{|c||c|c|c|c|}
  \hline
   Coils  & Type &$R$&($Z_{min}$, $Z_{max}$)&Turns\\
 \hline
 \hline
 CS1&Solenoid & 1.696&(-5.415, -3.6067)&  553  \\
 \hline
CS2&Solenoid & 1.696&(-3.6067, -1.7983) &  553  \\
  \hline
 CS3&Solenoid & 1.696&(-1.7983, 1.8183)&  1106  \\
    \hline
  CS4&Solenoid & 1.696&(1.8183, 3.6267) &  553  \\
 \hline
CS5&Solenoid & 1.696&(3.6267, 5.435)&  553  \\
 \hline
 \end{tabular}
\end{center}
\end{table}

 \begin{figure}
  \centering
       \centering
  \begin{subfigure}{.45\textwidth}
    \centering
    \includegraphics[trim=0cm 0cm 1cm 0cm, clip=true, width=\linewidth]{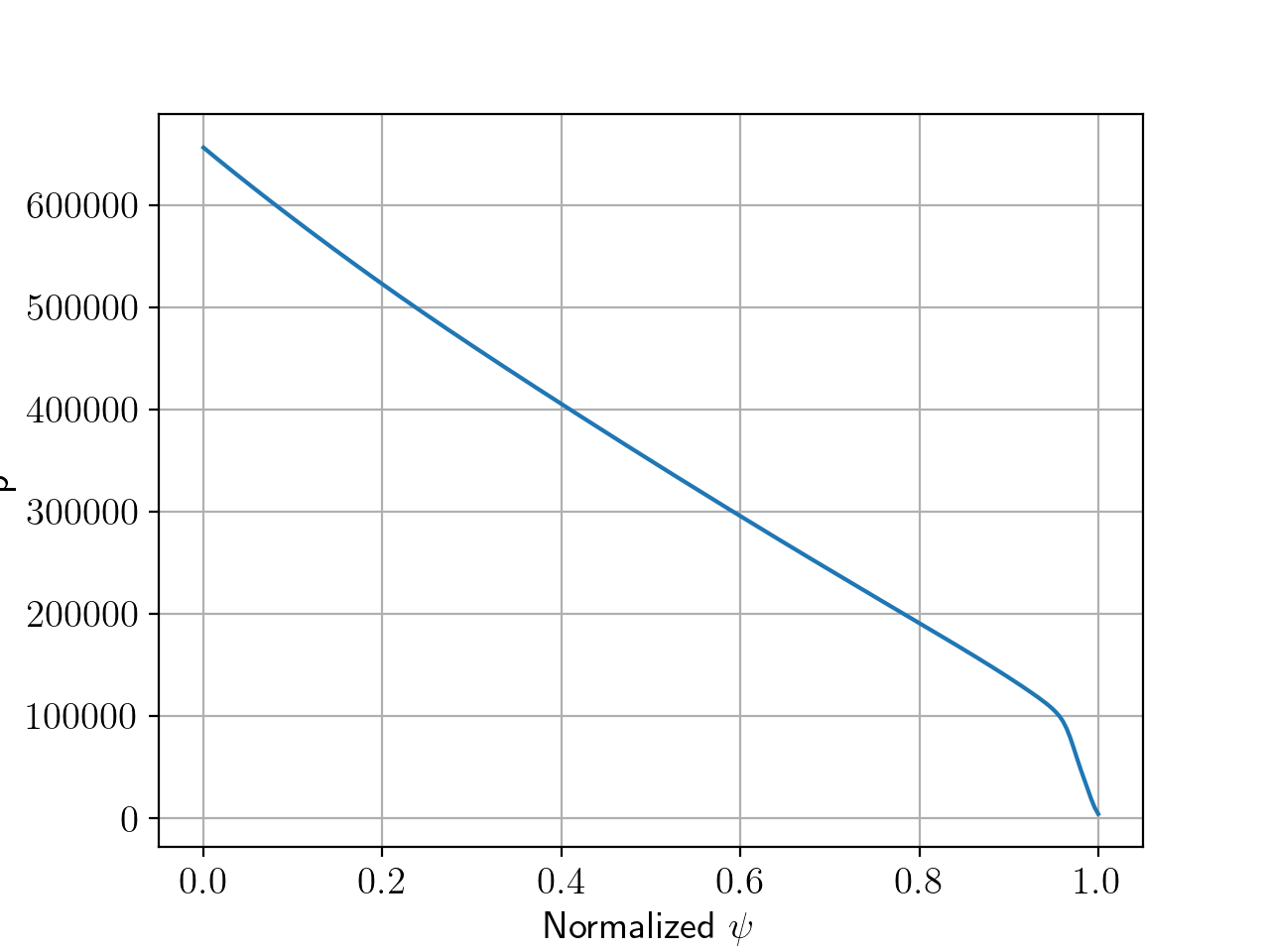}
     \caption{Pressure ($p$) profile of normalized $\psi$}
  \end{subfigure}%
  \hspace{.03\textwidth}
  \begin{subfigure}{.45\textwidth}
    \centering
    \includegraphics[trim=0cm 0cm 1cm 0cm, clip=true, width=\linewidth]{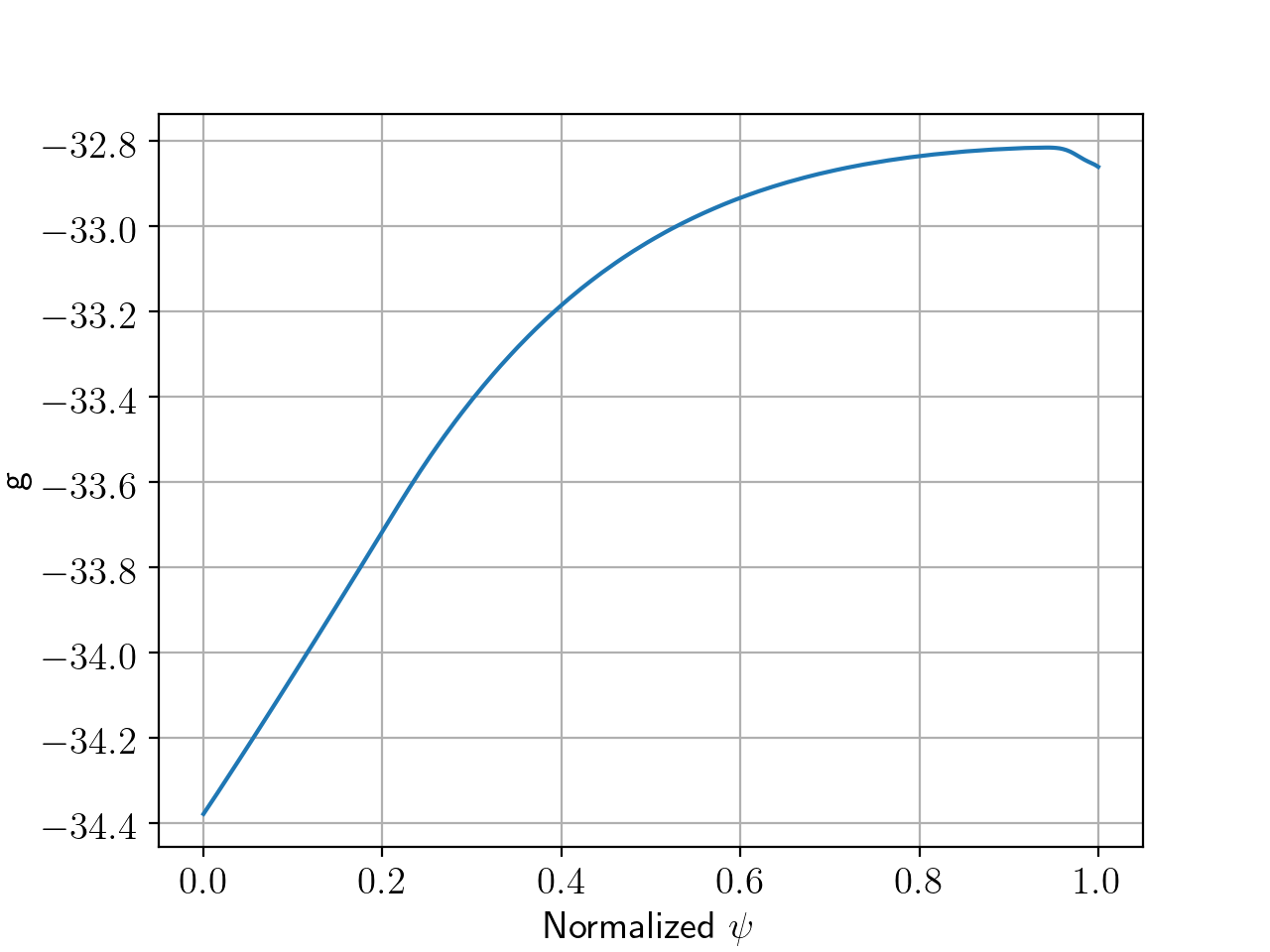}
    \caption{Poloidal current ($g$) profile of normalized $\psi$}
  \end{subfigure}%
  \caption{Numerical $p(\bar{\psi})$ and $g(\bar{\psi})$ profiles from the equilibrium data generated
    in Ref.~\cite{Liu-etal-NF-2015}.\label{fig:profiles}}
\end{figure}

 The problem is based on prescribed numerical $p(\bar{\psi})$ and
$g(\bar{\psi})$ profiles for a proposed ITER discharge at 15~MA toroidal
 plasma current~\cite{Liu-etal-NF-2015}, which carries the ITER
 reference number  ``ABT4ZL".  In addition to the numerical
 $p(\bar{\psi})$ and $g(\bar{\psi})$ profiles, there is also the
 numerical solution for $\psi(R,Z)$ from a different Grad-Shafranov
 solver~\cite{Liu-etal-NF-2015}, which we will use for the initial
 data.  The coil currents corresponding to the equilibrium were not
 given and they will be found by our free-boundary Grad-Shafranov
 solver in this test.  The ITER poloidal field coil locations,
 however, are fixed and given in~\Cref{tab:11}.  We will solve or more
 accurately, resolve this free-boundary Grad-Shafranov equilibrium
 problem in two ways.  One is based on a Cartesian mesh over a
 rectangular domain and the other one is based on a cut-cell mesh that
 has the first wall as the computational boundary.

\subsubsection{Rectangular domain} 
\label{sec:freeGSsolverinrectanguardomain}
The first example uses a rectangular computational domain that
contains the limiter regions but excludes the current regions.  This
is a conventional approach implemented in many practical
Grad-Shafranov solvers.  We present this case to verify our
implementation of the full solver as well as to compare the results
with the solutions from the cut-cell solver, which will be discussed
in the next section.

\begin{figure}
  \centering
  \begin{subfigure}{.33\textwidth}
    \centering
    \includegraphics[width=\linewidth]{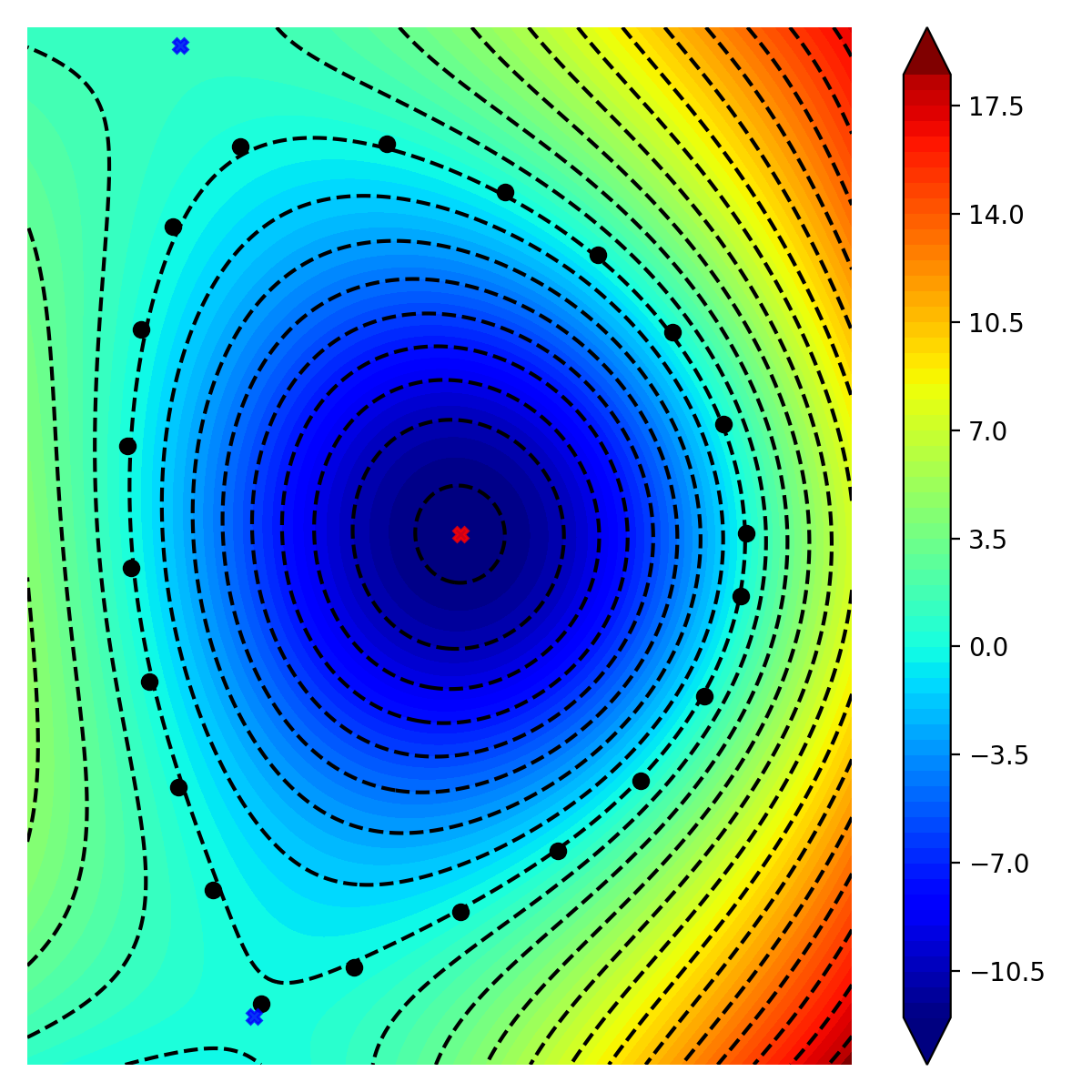}
     \caption{Converged $\psi$}
     \label{fig:p0convergedpsiVSinitialpsiimage1}
  \end{subfigure}%
  \begin{subfigure}{.33\textwidth}
    \centering
    \includegraphics[width=\linewidth]{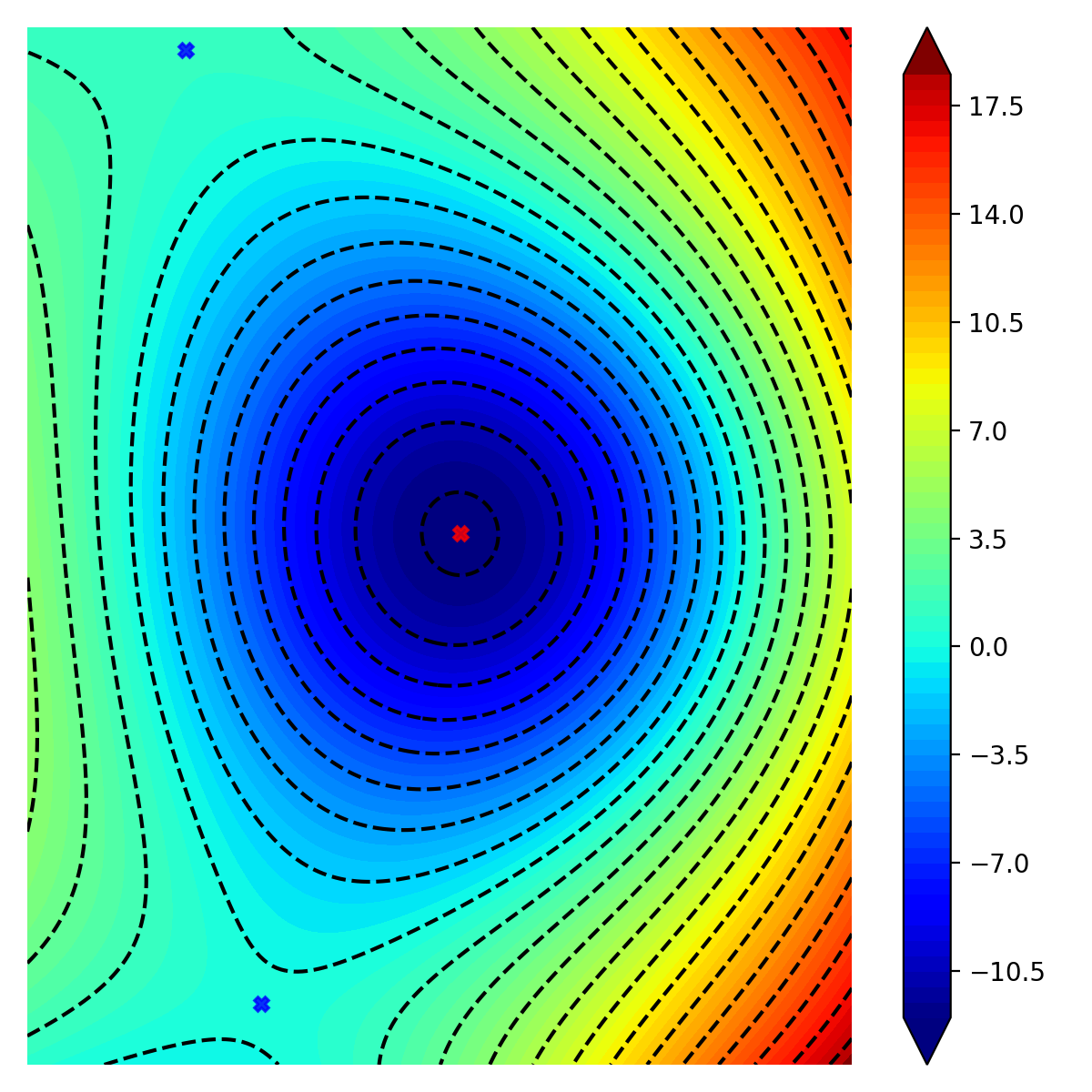}
     \caption{$\psi^0$ from initial data}
      \label{fig:p0convergedpsiVSinitialpsiimage2}
  \end{subfigure}%
  \begin{subfigure}{.33\textwidth}
    \centering
    \includegraphics[trim=1.4cm 1.3cm 0.2cm .3cm,clip,width=\linewidth]{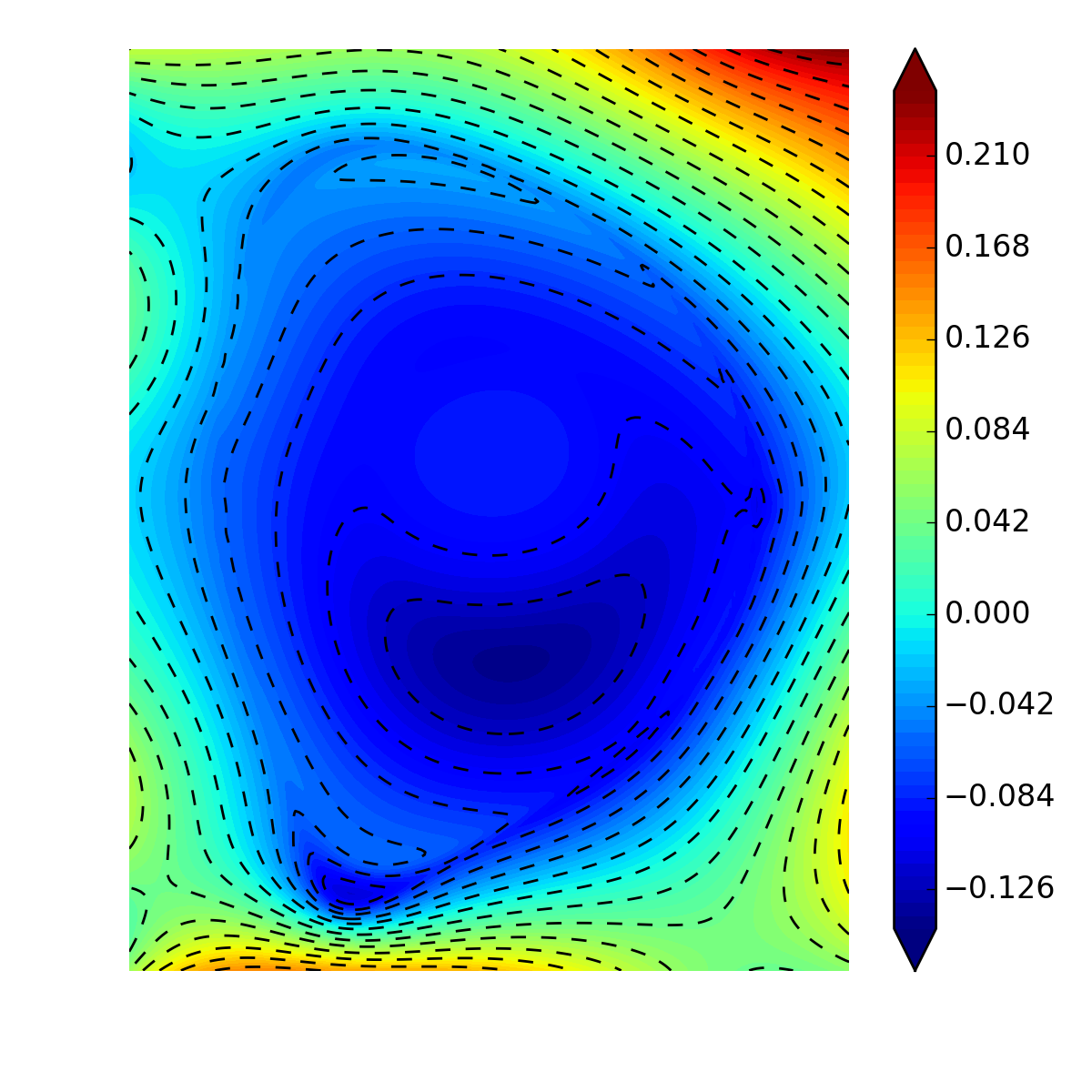}
     \caption{$\psi-\psi^0$}
      \label{fig:p0convergedpsiVSinitialpsiimage3}
  \end{subfigure}
  \caption{Rectangular geometry:  (a) converged magnetic flux from free-boundary Grad-Shafranov solver, (b) magnetic flux from initial equilibrium data file, (c) the difference between magnetic flux from free-boundary Grad-Shafranov solver and initial equilibrium data. }
  \label{fig:p0convergedpsiVSinitialpsi}
\end{figure}

The computational domain is $\Omega=\{(R,Z)\in \mathcal{H} \, \lvert
\, 3.55\leq R \leq 8.88 , -3.84 \leq Z \leq 4.92\}$ and the grid used
in this case is a Cartesian grid of size $196 \times 375$.  The
proposed free-boundary Grad-Shafranov solver is used to solve the
problem on the rectangular domain $\Omega$.  The converged magnetic
flux from the free-boundary Grad-Shafranov solver and magnetic flux
from the equilibrium data file are presented
in~\cref{fig:p0convergedpsiVSinitialpsi}.  Black dots in 
\cref{fig:p0convergedpsiVSinitialpsiimage1}  represent the shape control
points on targeted magnetic separatrix of the fixed plasma shape from
the initial data.  Twenty one points are selected along
$\psi=\psi_X$ in the initial data as the control points.  We can
observe from \cref{fig:p0convergedpsiVSinitialpsiimage1} that the converged
magnetic flux from free-boundary Grad-Shafranov solver keeps the
predetermined plasma shape very well.  The solution is comparable to
the initial data and its difference is presented in~\cref{fig:p0convergedpsiVSinitialpsiimage3}.  Note that the largest
difference is around the corner points.  This is expected since the
computational domain is very close to two current coils at $(8.3908,
-6.7269)$ and $(8.2851, 6.5398)$.  The magnetic axes $\psi_o$ are
presented in~\cref{fig:p0convergedpsiVSinitialpsiimage1} and ~\cref{fig:p0convergedpsiVSinitialpsiimage2} as red cross
points, which are clearly local minima in the solutions.  There are
two candidate of $\psi_X$ points shown
in~\cref{fig:p0convergedpsiVSinitialpsiimage1} and ~\cref{fig:p0convergedpsiVSinitialpsiimage2}, using blue cross points.  It
is clear that those are saddle points of the solutions.  Based on the
criteria we discussed in the algorithm section, the point in the
bottom portion of the domain is chosen as $\psi_X$.  It is seen that
the selected $\psi_X$ point in the solution is very close to the
control points.  Both $\psi_o$ and $\psi_X$ points are found to change
slightly throughout the iterations.

\begin{figure}
  \centering
  \begin{subfigure}{.33\textwidth}
    \centering
    \includegraphics[width=\linewidth]{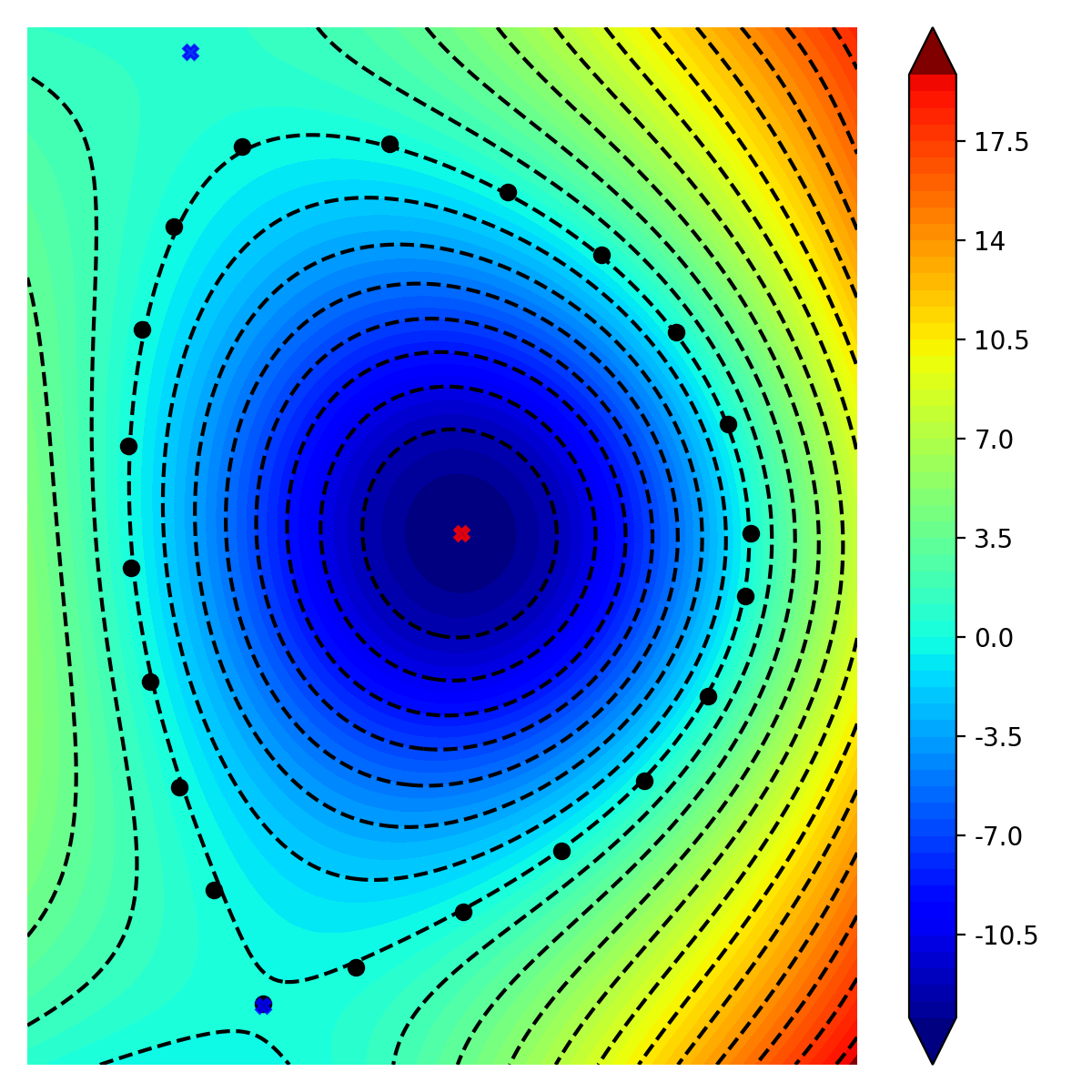}
     \caption{Converged $\psi$}
  \end{subfigure}%
  \begin{subfigure}{.33\textwidth}
    \centering
    \includegraphics[width=\linewidth]{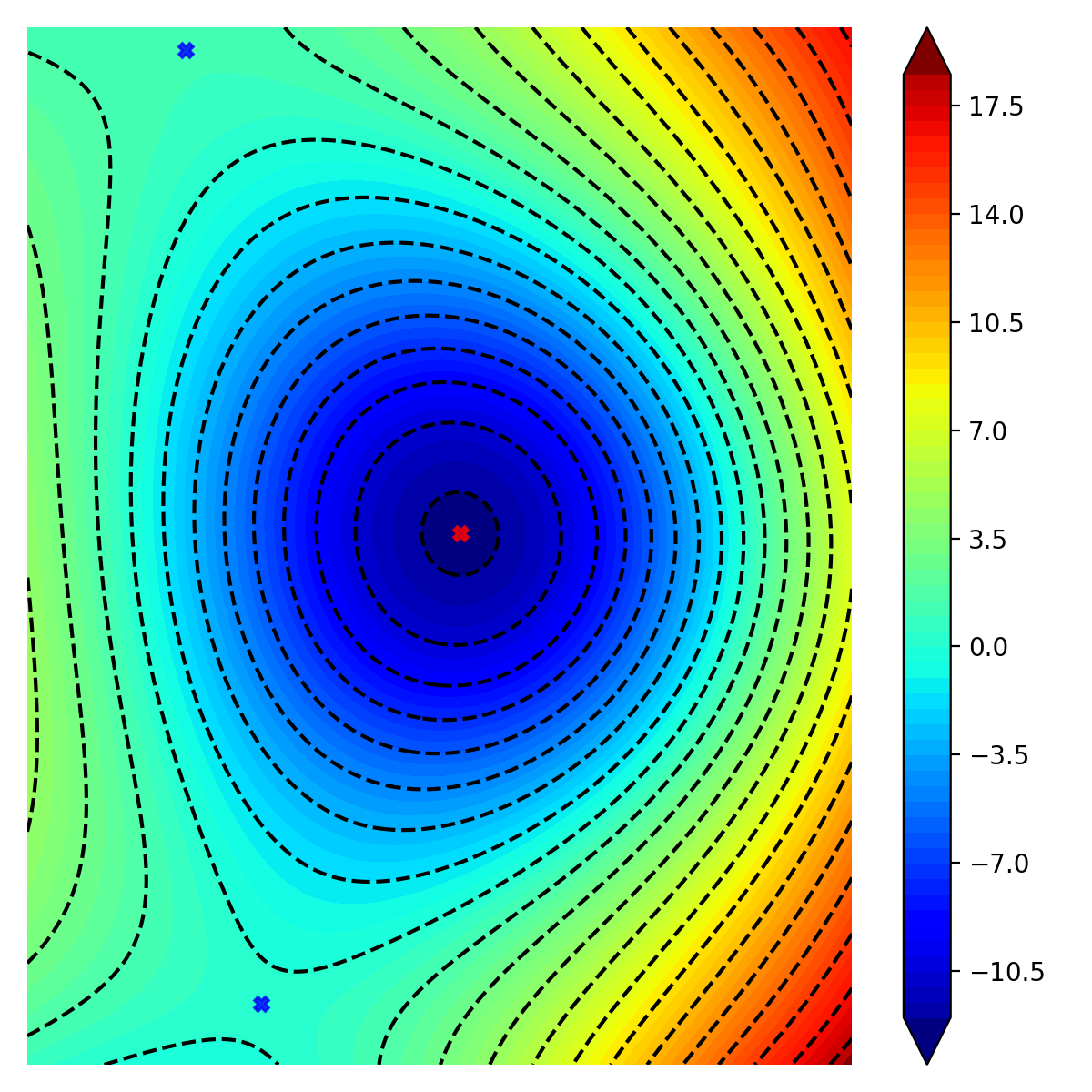}
     \caption{$\psi^0$ from initial data}
  \end{subfigure}%
  \begin{subfigure}{.33\textwidth}
    \centering
    \includegraphics[trim=1.4cm 1.3cm 0.2cm .3cm,clip,width=\linewidth]{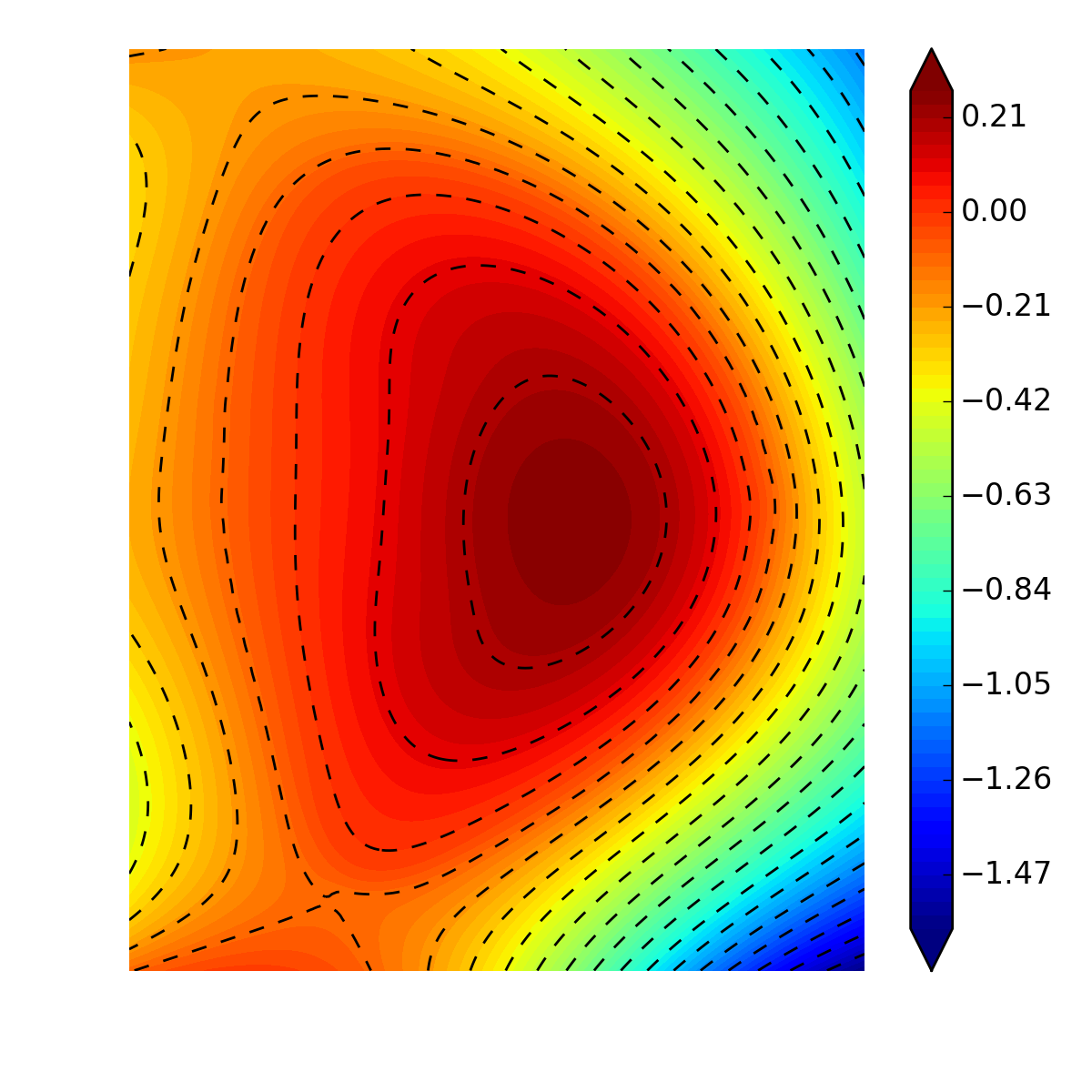}
     \caption{$\psi-\psi^0$}
  \end{subfigure}
  \caption{Rectangular geometry:  (a) converged magnetic flux from free-boundary Grad-Shafranov solver, (b) magnetic flux from initial equilibrium data file, (c) the difference between magnetic flux from free-boundary Grad-Shafranov solver and the initial data. The pressure is dropped to 80\% of the given profile.}
  \label{fig:08p0psiVSinitialpsi}
\end{figure}

Another interesting question to investigate is whether the converged
magnetic flux from the free-boundary Grad-Shafranov solver can keep
the predetermined plasma shape if we modify the source term of the
Grad-Shafranov equation.  This is a test to verify the effectiveness
of the shape control techniques we implemented in the solver.  As an
example, we drop the pressure to 80\% of the original pressure profile
and solve the same free-boundary problem using the solver.  Its
results are presented in~\cref{fig:08p0psiVSinitialpsi}.  The
converged magnetic flux from the free-boundary Grad-Shafranov solver
and magnetic flux from the equilibrium data file are also compared
in~\cref{fig:08p0psiVSinitialpsi}. It indicates that the converged
magnetic flux still keeps the predetermined plasma shape. It also
suggests that the free-boundary Grad-Shafranov solver performs well in
keeping the targeted shape of $\mathcal{P}(\psi)$ in the rectangular
domain.

\subsubsection{Limiter-bounded domain $ \mathcal{L}$}

In the second example, by combining cut-cell algorithm, we directly choose the irregular limiter-bounded domain $ \mathcal{L}$  as $\Omega$ to solve the free boundary problem. 
The domain  $\mathcal{L}$ consists of a bunch of straight lines due to the design of ITER, which will be used to construct the cut-cell mesh.
This example is used to verify the full algorithm we proposed in this work.

\begin{figure}
  \centering
  \begin{subfigure}{.25\textwidth}
    \centering
    \includegraphics[trim=1cm 0cm 0cm 0cm,clip, width=\linewidth]{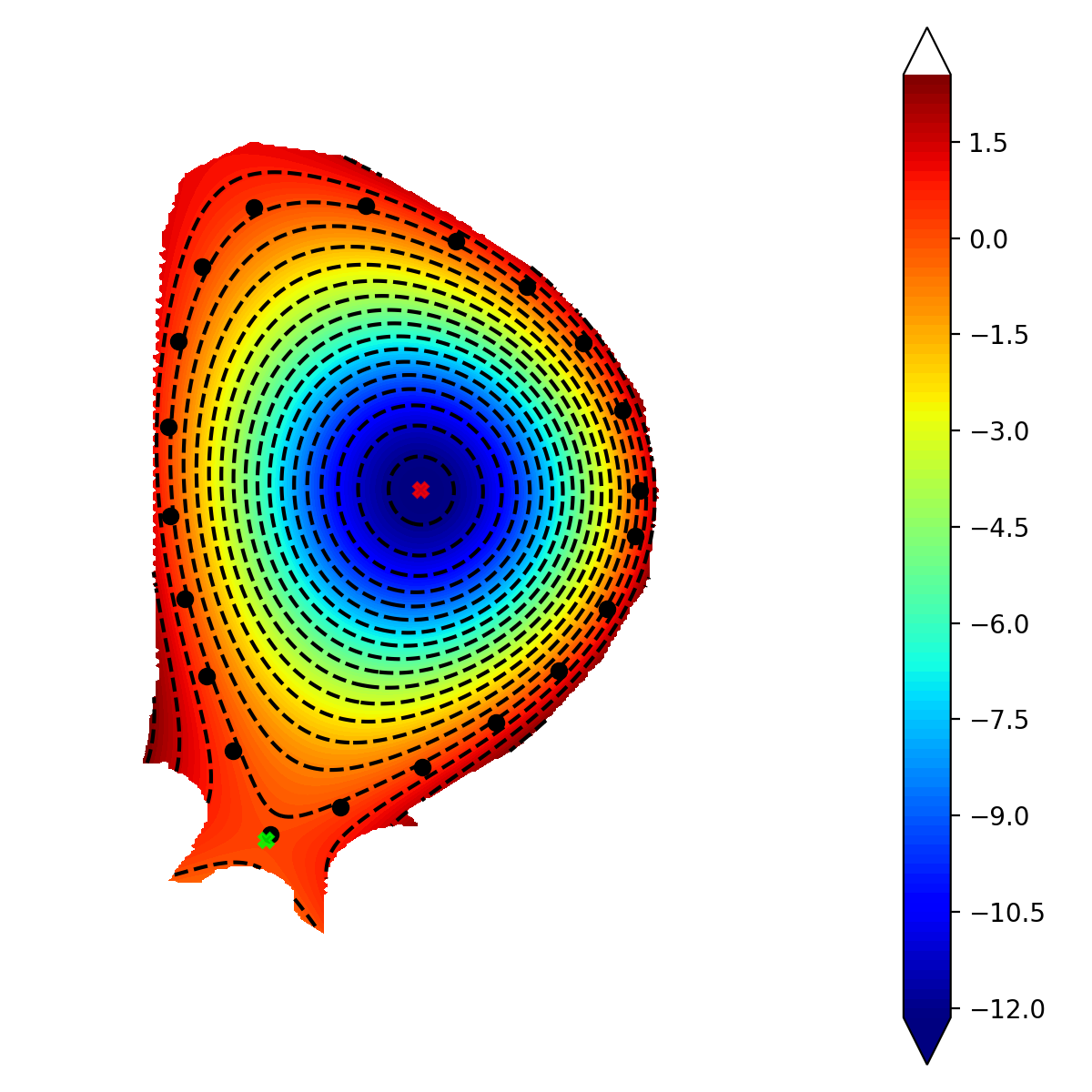}
     \caption{Converged $\psi_{\rm cut-cell}$}
       \label{fig:p0convergedpsiVSinitialpsiincutcellimage1}
  \end{subfigure}%
  \begin{subfigure}{.25\textwidth}
    \centering
    \includegraphics[trim=1cm 0cm 0cm 0cm,clip, width=\linewidth]{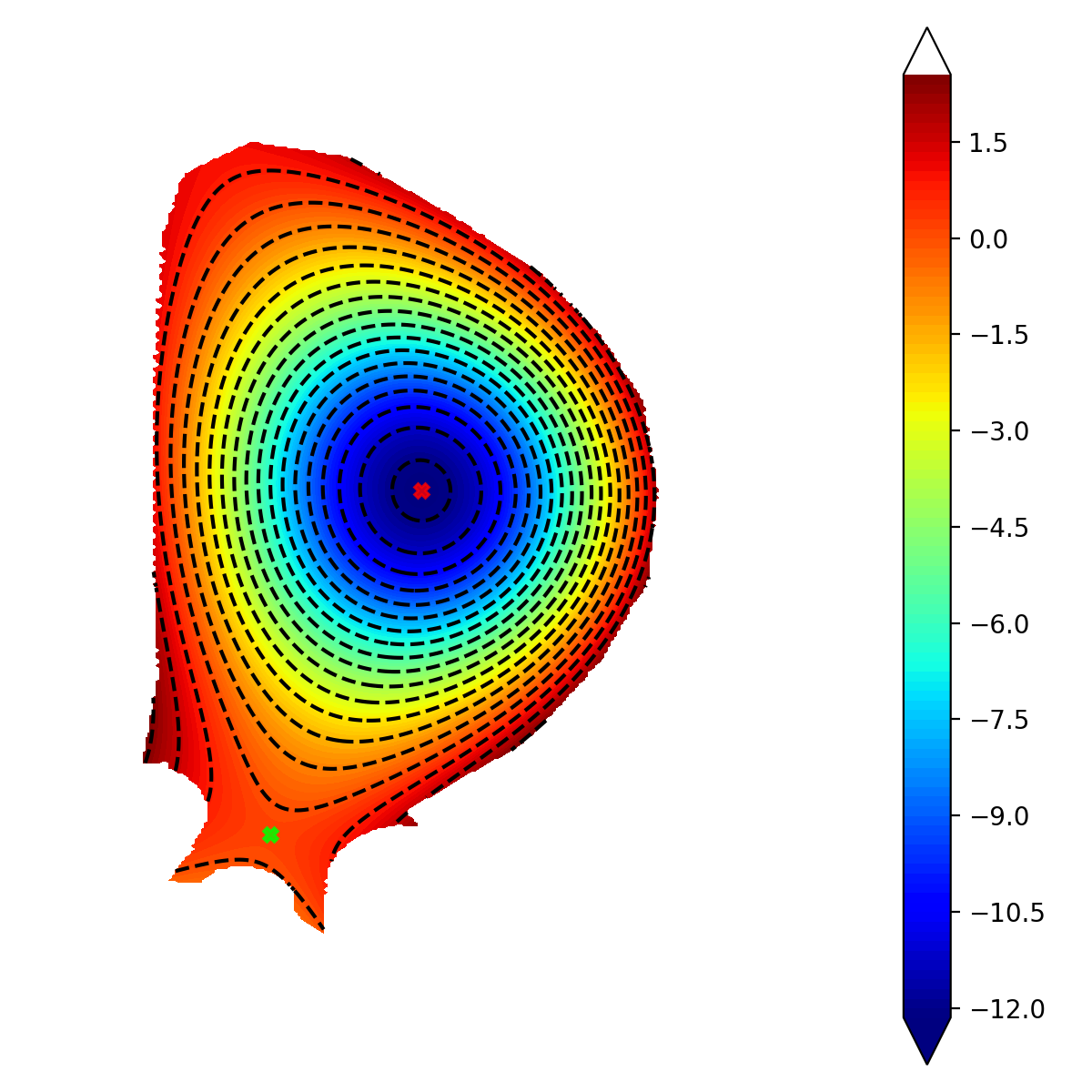}
     \caption{$\psi^0$ from initial data}
      \label{fig:p0convergedpsiVSinitialpsiincutcellimage2}
  \end{subfigure}%
  \begin{subfigure}{.25\textwidth}
    \centering
    \includegraphics[trim=1cm 0cm 0cm 0cm, clip,width=\linewidth]{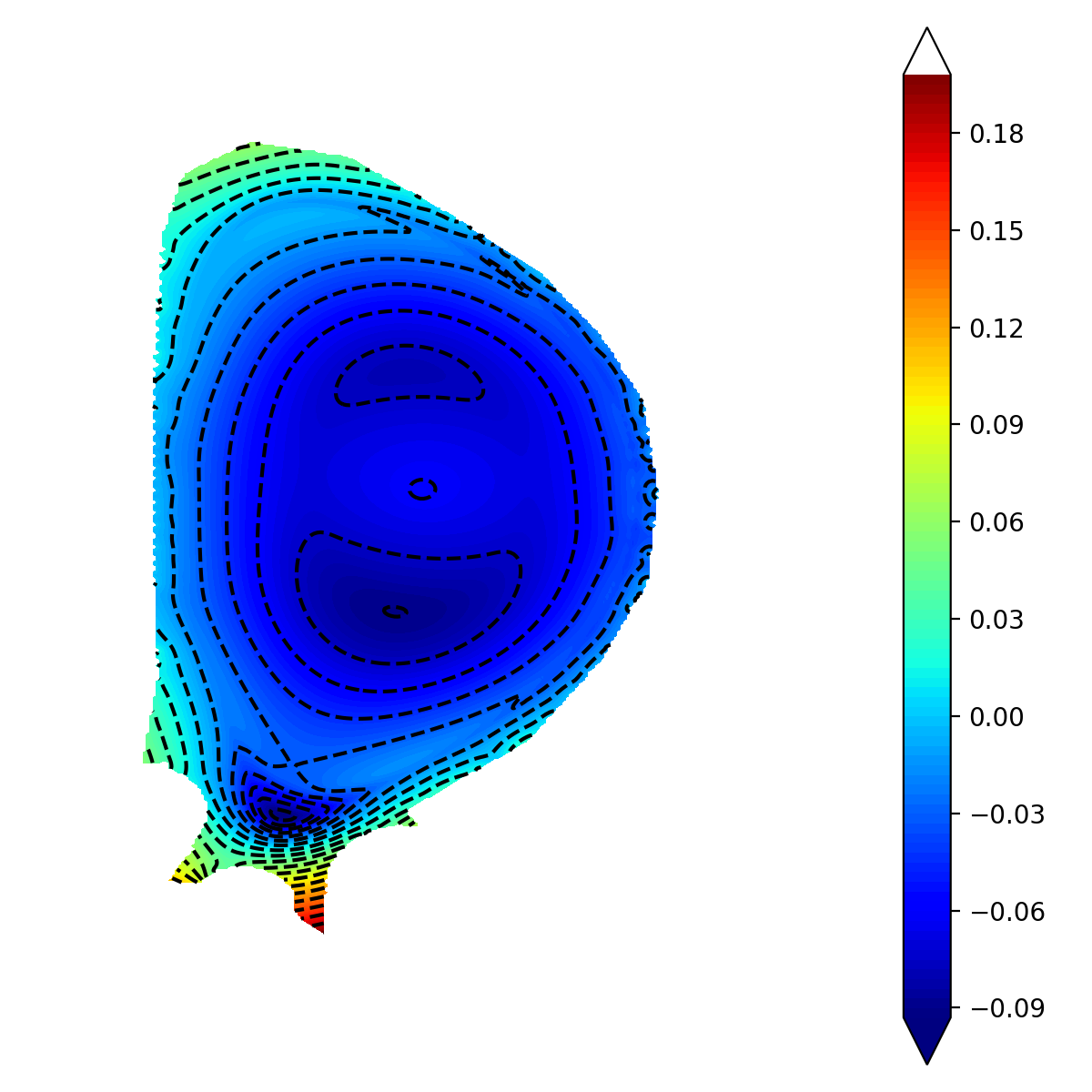}
     \caption{$\psi_{\rm cut-cell}-\psi^0$}
      \label{fig:p0convergedpsiVSinitialpsiincutcellimage3}
  \end{subfigure}%
  \begin{subfigure}{.25\textwidth}
    \centering
    \includegraphics[trim=1cm 0cm 0cm 0cm, clip,width=\linewidth]{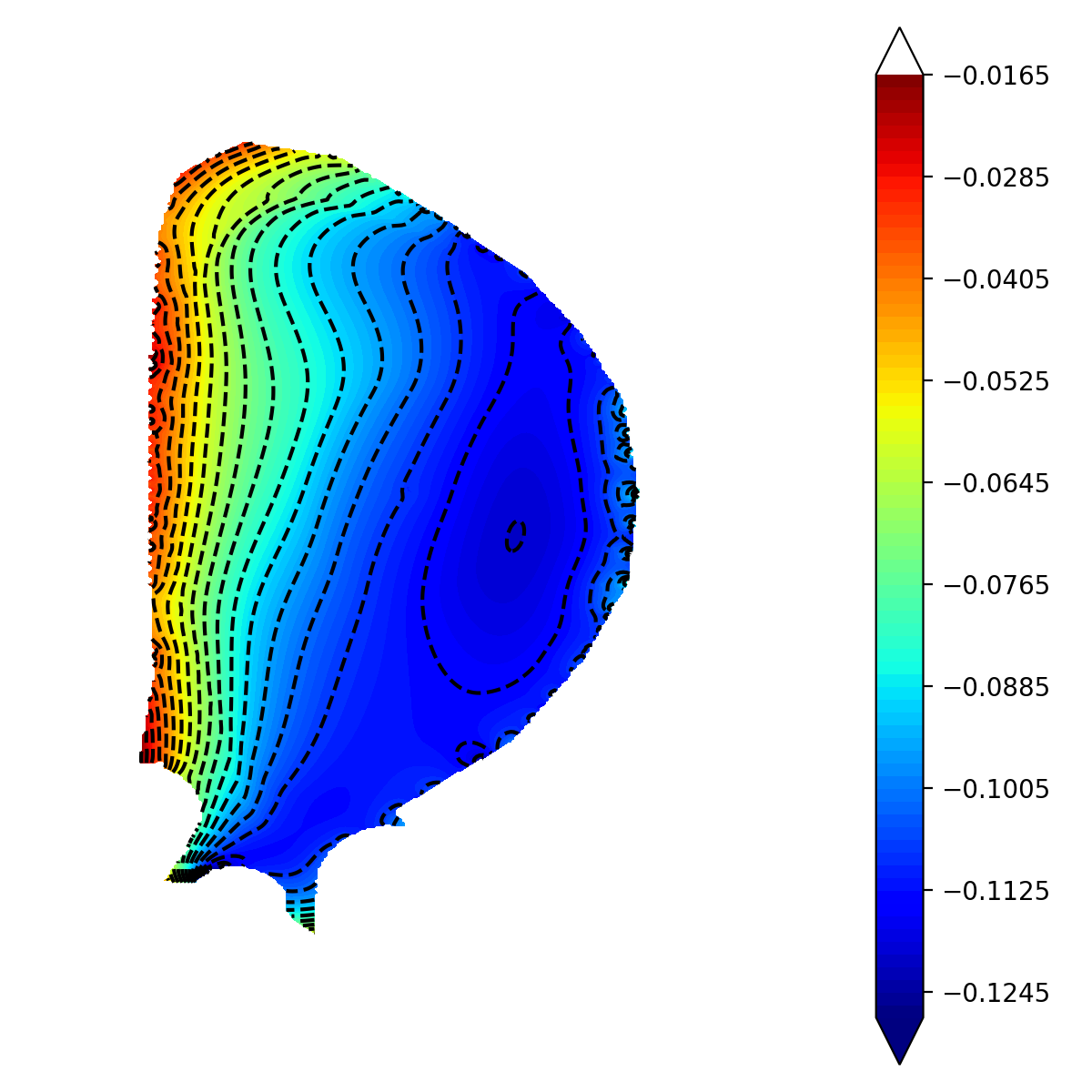}
     \caption{$\psi_{\rm box}-\psi_{\rm cut-cell}$}
      \label{fig:p0convergedpsiVSinitialpsiincutcellimage4}
  \end{subfigure}
  \caption{Limiter geometry: (a) converged magnetic flux from the free-boundary Grad-Shafranov solver with the cut-cell algorithm, (b) magnetic flux from the initial equilibrium data file, (c) difference between magnetic flux from the free-boundary Grad-Shafranov solver with the cut-cell algorithm and the initial equilibrium data file, (d) difference between the solutions of the rectangular domain and the limiter-bounded domain.}
  \label{fig:p0convergedpsiVSinitialpsiincutcell}
\end{figure}

We create a Cartesian mesh of size $216 \times493$ in the domain $\{(R,Z)\in \mathcal{H} \, \lvert \, 3.0\leq R \leq 8.88 ,  -5.53 \leq Z \leq 6.0\}$, which contains  the irregular limiter-bounded domain $ \mathcal{L}$. The cut-cell mesh is created accordingly, and  the total number of active cell  is 44668, which includes 1266 cut-cells. The proposed free-boundary Grad-Shafranov solver combining the cut-cell algorithm is used to solve the problem on the limiter-bounded domain $ \mathcal{L}$. The converged magnetic flux from the solver  and  magnetic flux  from the equilibrium data file are presented in~\cref{fig:p0convergedpsiVSinitialpsiincutcell}.  Black dots in \cref{fig:p0convergedpsiVSinitialpsiincutcellimage1} represent the shape control points on targeted magnetic separatrix of the fixed plasma shape from the initial data.  The same 21 points in \cref{sec:freeGSsolverinrectanguardomain} are selected along $\psi=\psi_X$ in the initial data as the control points. 
We can observe from  \cref{fig:p0convergedpsiVSinitialpsiincutcellimage1}  that the converged magnetic flux from free-boundary Grad-Shafranov solver keeps the predetermined plasma shape very well. 
The solution is comparable to the initial data and its difference is presented in~\cref{fig:p0convergedpsiVSinitialpsiincutcellimage3}.
Note that the largest difference is around the bottom corner, which is due to the non-smooth boundary in the computational domain. 
The magnetic axes $\psi_o$ are presented in~\cref{fig:p0convergedpsiVSinitialpsiincutcellimage1}  and~\cref{fig:p0convergedpsiVSinitialpsiincutcellimage2} as red cross points, which are clearly local minima in the solutions. Different from the example in the rectangular domain, only one $\psi_X$ point is shown in~\cref{fig:p0convergedpsiVSinitialpsiincutcellimage1}  and ~\cref{fig:p0convergedpsiVSinitialpsiincutcellimage2}, using a green cross point. This is caused by the choice of the limiter-bounded domain $ \mathcal{L}$, as the initial plasma domain is bounded by the last closed poloidal flux line inside the limiter-bounded domain. 
It is clear that the $\psi_X$ point is the saddle point of the solutions. Both $\psi_o$ and $\psi_X$ are found to change slightly throughout the iterations.

In~\cref{fig:p0convergedpsiVSinitialpsiincutcellimage4}, the converged magnetic fluxes of the free-boundary problem on the regular domain and the limiter-bounded domain are compared to each other. It is found that the difference is small, which indicates both solvers produce the same equilibrium.
When solving on the rectangular domain, the outer iteration is 12, and the average inner iteration is 32, while on the  limiter-bounded domain the outer iteration is 11, and the average inner iteration is 9.

\begin{figure}
  \centering
  \begin{subfigure}{.33\textwidth}
    \centering
    \includegraphics[width=\linewidth]{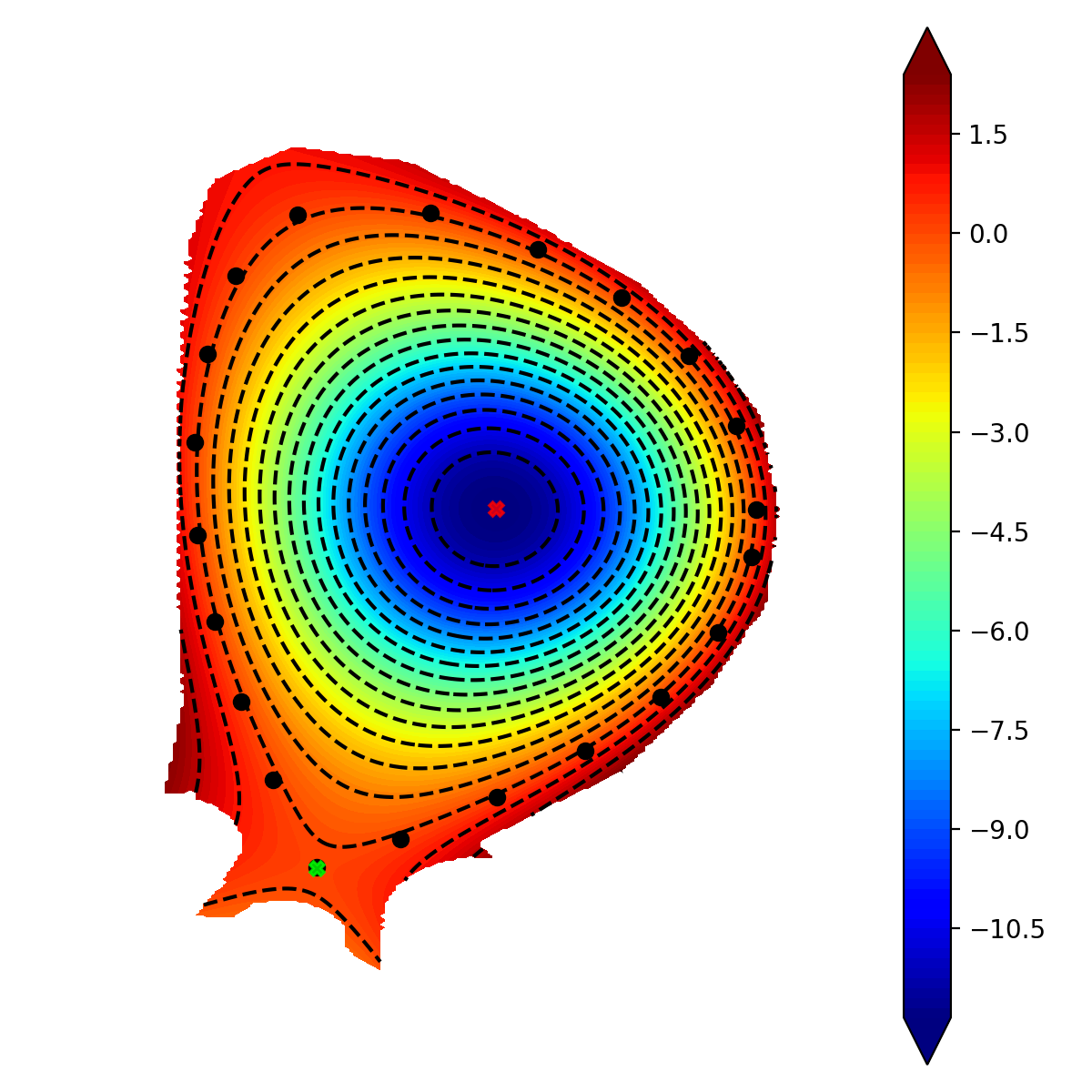}
     \caption{Converged $\psi_{\rm cut-cell}$}
       \label{fig:08p0convergedpsiVSinitialpsiincutcellimage1}
  \end{subfigure}%
  \begin{subfigure}{.33\textwidth}
    \centering
    \includegraphics[width=\linewidth]{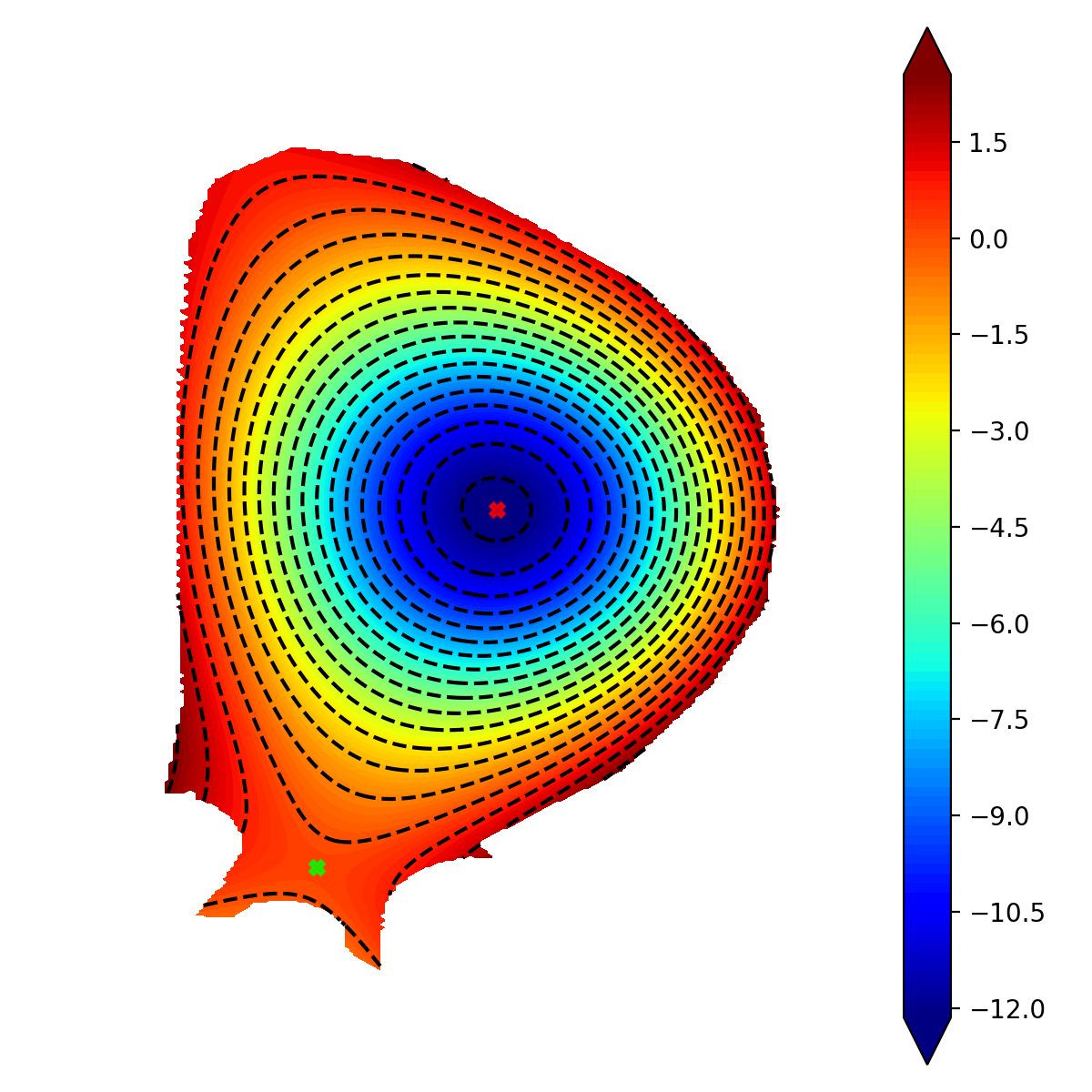}
     \caption{$\psi^0$ from initial data}
      \label{fig:08p0convergedpsiVSinitialpsiincutcellimage2}
  \end{subfigure}%
  \begin{subfigure}{.33\textwidth}
    \centering
    \includegraphics[width=\linewidth]{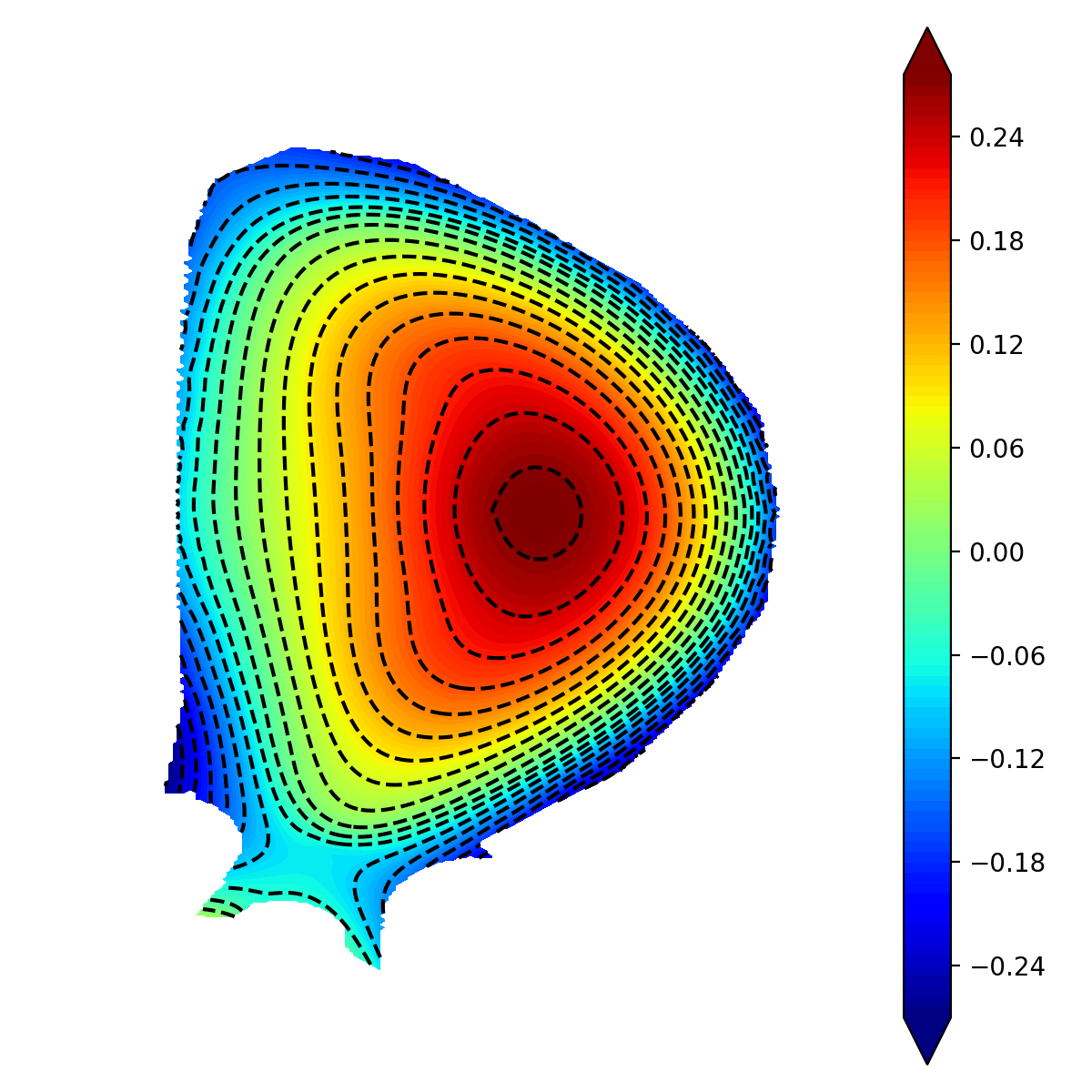}
     \caption{$\psi_{\rm cut-cell}-\psi^0$}
      \label{fig:08p0convergedpsiVSinitialpsiincutcellimage3}
  \end{subfigure}%
  \caption{Limiter geometry: (a) converged magnetic flux from the free-boundary Grad-Shafranov solver with the cut-cell algorithm, (b) magnetic flux from the initial equilibrium data file, (c) difference between magnetic flux from the free-boundary Grad-Shafranov solver with the cut-cell algorithm and the initial equilibrium data file.The pressure is dropped to 80\% of the given profile.}
  \label{fig:08p0convergedpsiVSinitialpsiincutcell}
\end{figure}

{In this example, we again investigate whether the converged magnetic flux from the solver can keep the predetermined plasma shape in  the limiter-bounded domain $ \mathcal{L}$ when we modify the source term of the Grad-Shafranov equation. Similarly, we drop the pressure to 80\% of the original pressure profile and solve the same free-boundary problem using the solver combining with the cut-cell algorithm. 
Its results are presented in~\cref{fig:08p0convergedpsiVSinitialpsiincutcell}. The converged magnetic flux from the free-boundary Grad-Shafranov solver and the magnetic flux  from the equilibrium data file are also presented in~\cref{fig:08p0convergedpsiVSinitialpsiincutcell}. It indicates that the converged magnetic flux still keeps the predetermined plasma shape. It also suggests that the  free-boundary Grad-Shafranov solver performs well in keeping the targeted shape of $\mathcal{P}(\psi)$ in  the limiter-bounded domain $ \mathcal{L}$.}

\subsection{Performance of Aitken's acceleration}
\label{sec: Aitken}
In the final examples, we focus on the demonstration of the efficiency of Aitken's acceleration. Results from Aitken's acceleration are compared to results from 
 Picard iteration with fixed under-relaxation coefficients $\alpha$ in the free-boundary Grad-Shafranov solver.
As presented in~\Cref{alg:2loops}, both inner and outer loops use under-relaxations, either through Aitken's acceleration or a fixed under-relaxed parameter. 
We consider both the rectangular geometry and the geometry with the limiter-bounded domain in the previous test.
Aitken's acceleration is found to converge much faster than Picard  iterations with fixed under-relaxation coefficients in  the both inner  and outer loops.

\subsubsection{Rectangular domain}
\label{sec: Aitkeninrectangulardomain}
In this example, we test the efficiency of Aitken's acceleration when compared with Picard iteration with fixed under-relaxation coefficients $\alpha$ on 
the rectangular domain with mesh size of   $196 \times 375$. Here we set the range $[\lambda_{min}, \lambda_{max}] = [0,0.95]$. 

\begin{figure}
  \centering
  \begin{subfigure}{.45\textwidth}
    \centering
     \includegraphics[trim=0cm 0cm 1cm 0cm, clip=true, width=\linewidth]{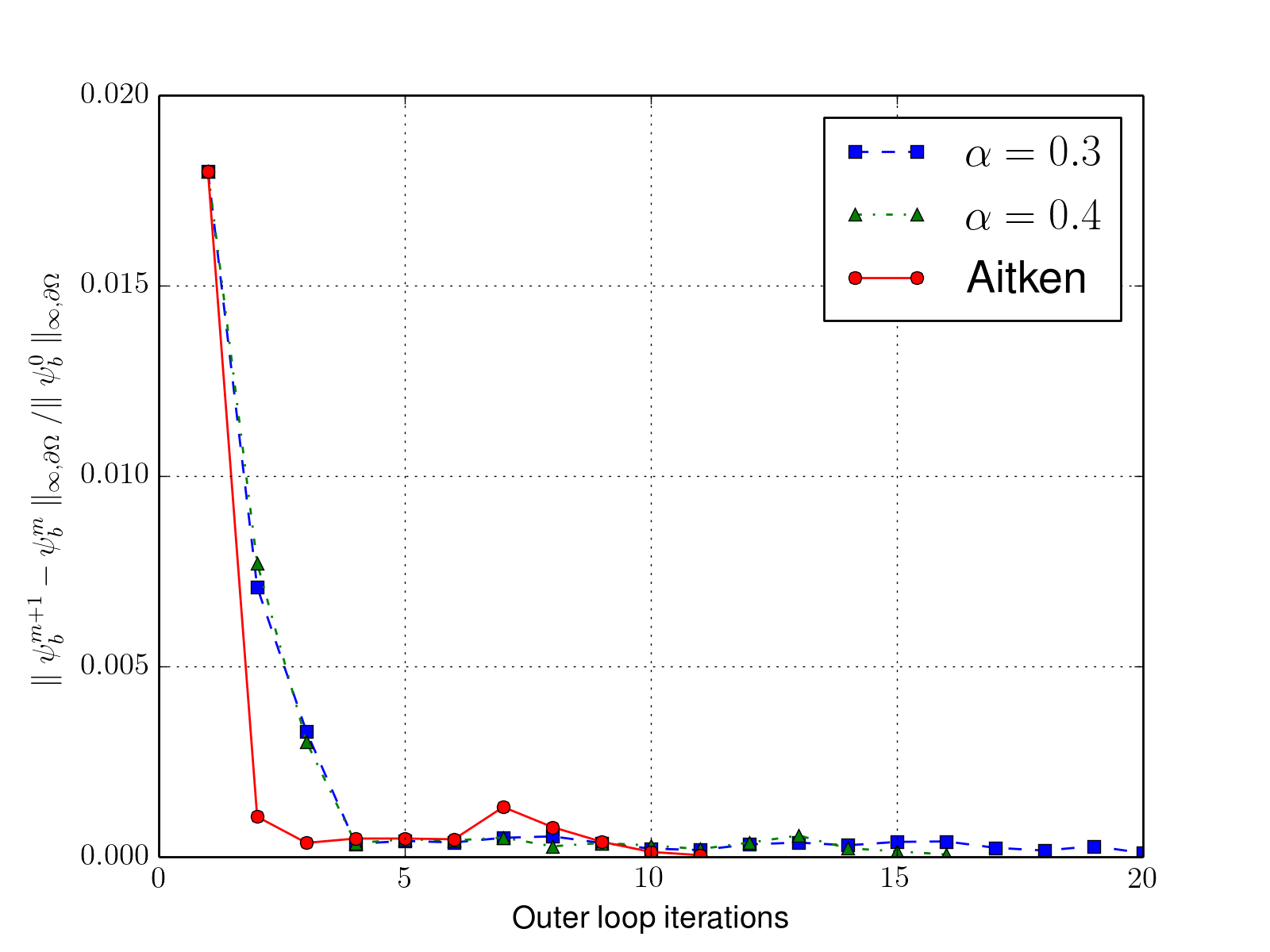}
  \end{subfigure}%
  \hspace{.5cm}
  \begin{subfigure}{.45\textwidth}
    \centering
      \includegraphics[trim=0cm 0cm 1cm 0cm, clip=true, width=\linewidth]{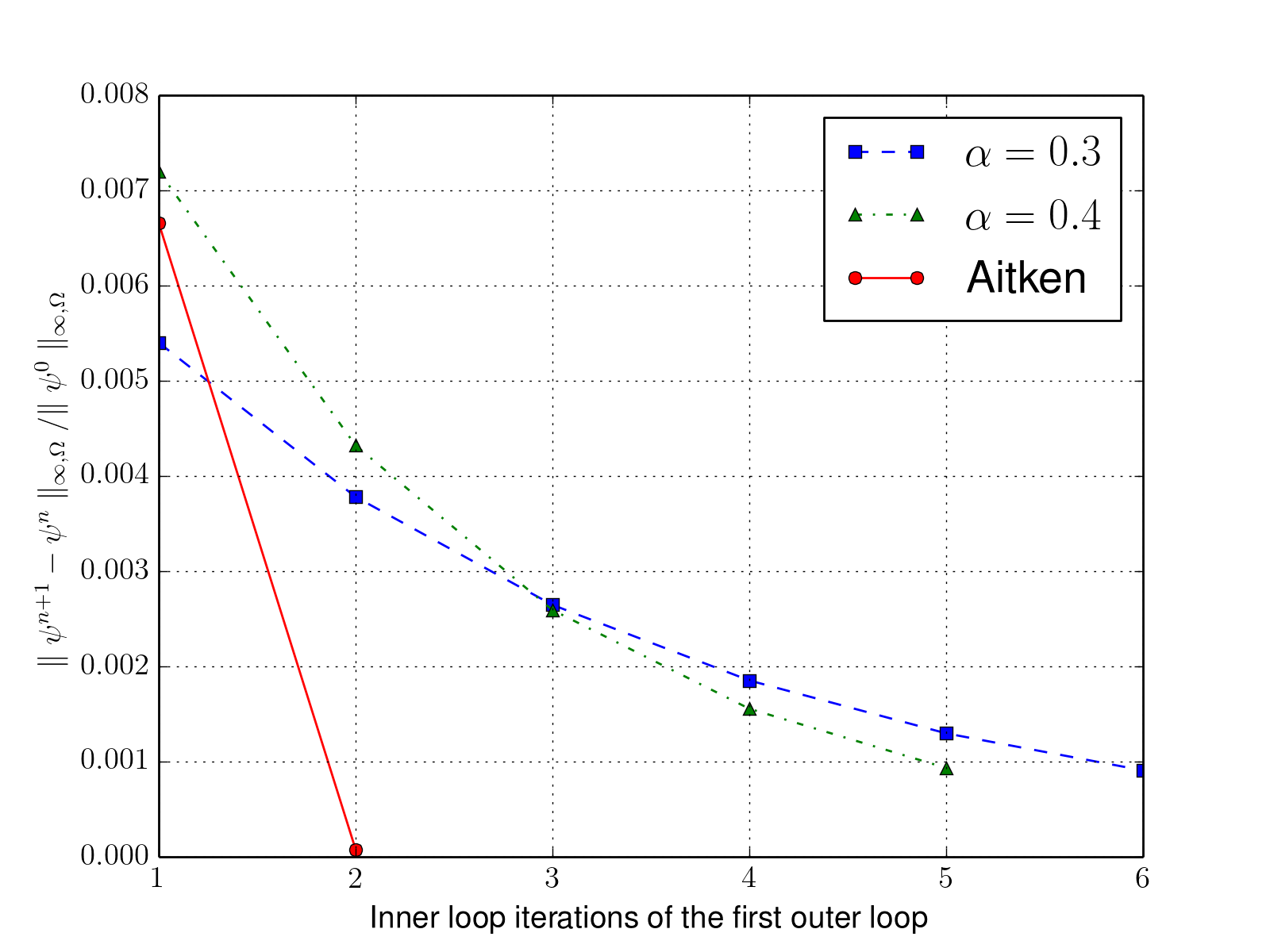}
  \end{subfigure}%
  \caption{Comparison of convergence history of the free-boundary solver algorithm with Aitken's acceleration and Picard iterations.
  Left: convergence history of the outer loop. 
  Right: convergence history of the first inner loop.
The convergence histories of two fixed under-relaxation parameter, $\alpha=0.3$ and $\alpha=0.4$, are also presented. 
The problem is solved on a rectangular computational domain with mesh size of  $196 \times 375$. Note that Aitken's acceleration significantly improve the convergence of Picard iterations with a predetermined under-relaxation. 
}
  \label{fig:aitkenVSpicardinrectangulardomain}
\end{figure}

\cref{fig:aitkenVSpicardinrectangulardomain} shows the convergence histories of Picard iteration with different under-relaxation coefficients $\alpha$ and  Aitken's acceleration.  
The results of two parameters, $\alpha =0.3$ and $\alpha =0.4$, are presented for comparison. 
The left figure shows that the free-boundary Grad-Shafranov solver algorithm with under-relaxed Picard iterations need to take 16 or 20 outer loops to satisfy the desired convergent threshold
for the boundary value $\psi_b$, 
while Aitken's acceleration only needs 11 iterations.
On the other hand, the right figure shows the convergence histories in the first inner loop. 
To have a fair comparison, the first inner loop is chosen since all the first inner loops use the same initial guess. It is found that the Aitken's acceleration only needs 2 iterations to satisfy the the desired convergence tolerance for the magnetic flux $\psi$, while the under-relaxed Picard iterations need to take 2 to 3 times more iterations. Therefore, the performance of Aitken's acceleration in the inner loops is also better than that  of Picard iterations.

\subsubsection{Limiter-bounded domain $ \mathcal{L}$}
\label{sec: Aitkenwithcutcell}
The final example focuses on 
the efficiency of Aitken's acceleration on the limiter-bounded domain $ \mathcal{L}$ with a cut-cell mesh of size  $216 \times493$ (base Cartesian mesh).
 Here we set the range $[\lambda_{min}, \lambda_{max}] = [0,0.7]$. 

\begin{figure}
  \centering
  \begin{subfigure}{.45\textwidth}
    \centering
     \includegraphics[trim=0cm 0cm 1cm 0cm, clip=true, width=\linewidth]{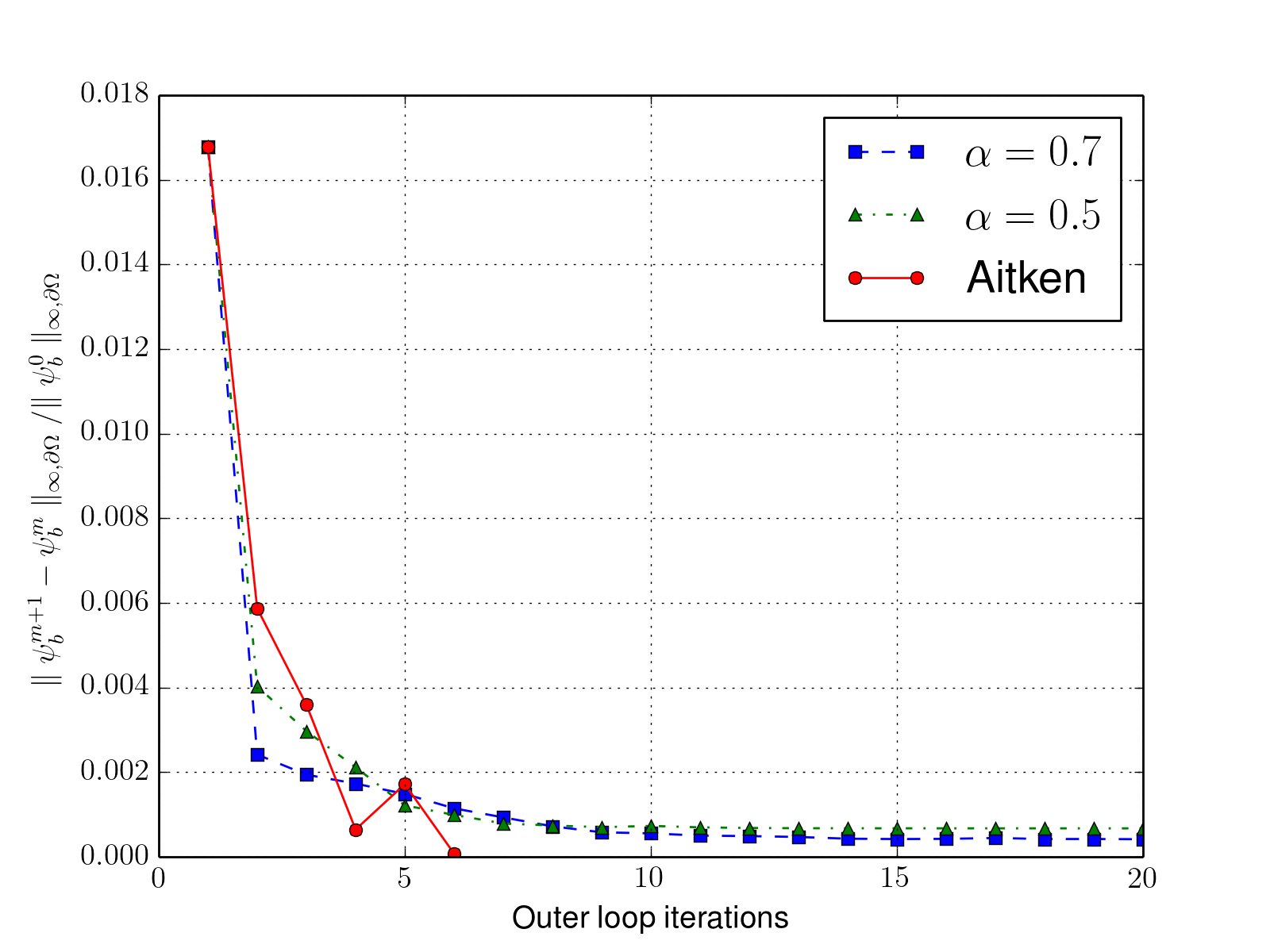}
  \end{subfigure}%
    \hspace{.5cm}
  \begin{subfigure}{.45\textwidth}
    \centering
      \includegraphics[trim=0cm 0cm 1cm 0cm, clip=true, width=\linewidth]{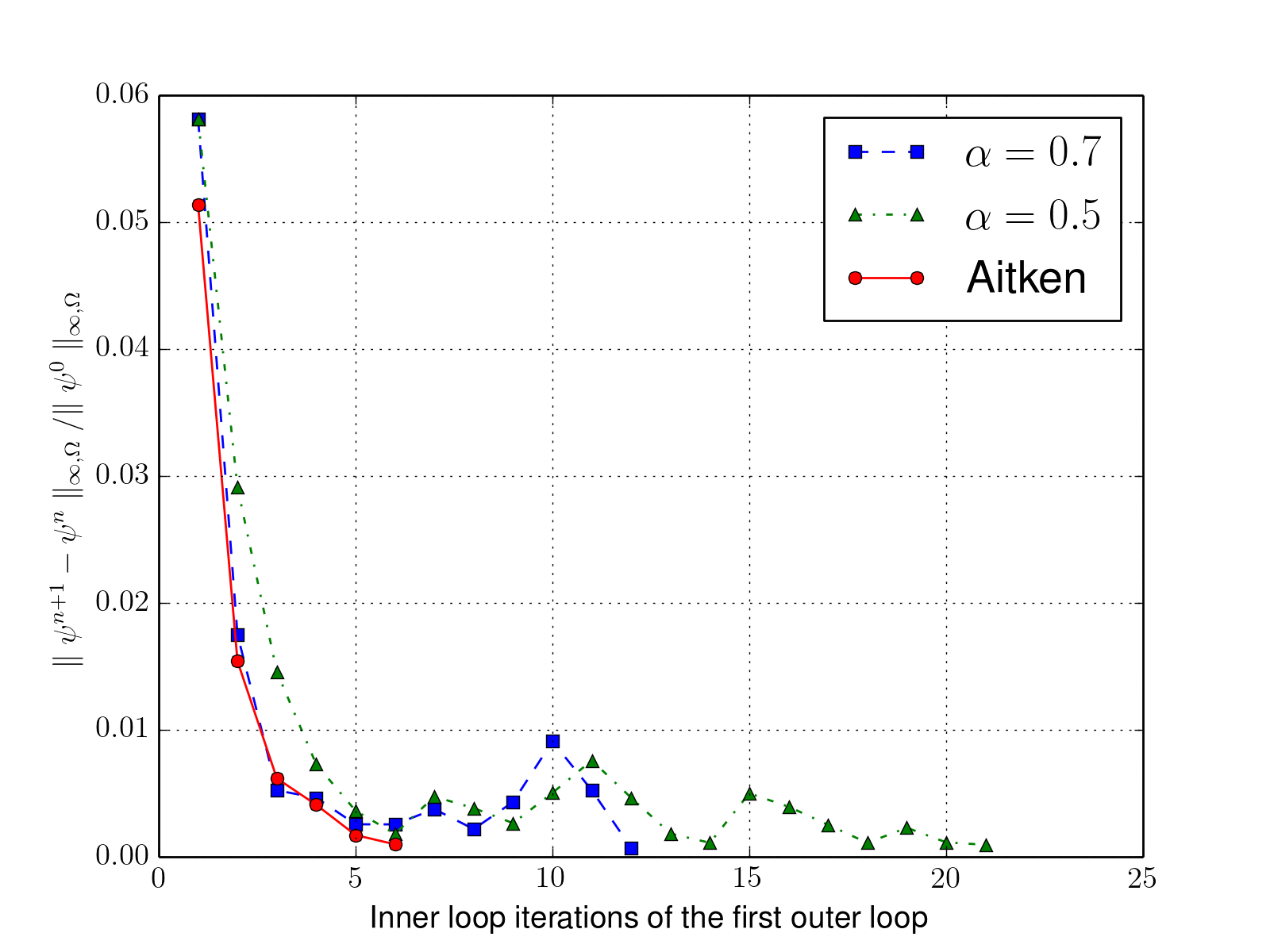}
  \end{subfigure}%
  \caption{
  Comparison of convergence history of the free-boundary solver algorithm with Aitken's acceleration and Picard iterations.
  Left: convergence history of the outer loop. 
  Right: convergence history of the first inner loop.
The convergence histories of two fixed under-relaxation parameter, $\alpha=0.7$ and $\alpha=0.5$, are also presented. The problem is solved on the limiter-bounded domain $ \mathcal{L}$ with a cut-cell mesh of size  $216 \times493$  (base Cartesian mesh). Note that Aitken's acceleration significantly improve the convergence of Picard iterations with a predetermined under-relaxation. }
  \label{fig:aitkenVSpicardinlimiterboundeddomain}
\end{figure}

\cref{fig:aitkenVSpicardinlimiterboundeddomain} shows the convergence histories of Picard iterations with different under-relaxation coefficients $\alpha$ and Aitken's acceleration.  
The results of two parameters, $\alpha =0.7$ and $\alpha =0.5$  are presented for comparison. 
We set the the maximum iteration to be 40 for the outer loop. The under-relaxed Picard iterations fail to reach the desired tolerance before reaching the maximum iteration,
while the Aitken's acceleration only needs 6 iterations to reach the same tolerance. The left figure presents the maximum relative difference of $\psi_b$ in the first 20 iterations. 
The right figure shows the first inner loop, for which the Aitken's acceleration only needs 6 iterations, while the under-relaxed Picard iterations take 12 or 21 iterations.

In conclusion, we find that the Aitken's acceleration could significantly
accelerate Picard iterations in both inner and outer loops. Other
acceleration techniques may also show similar performance, which will
be studied in the future work.

\section{Conclusions}

This work discusses the development of a parallel free-boundary
Grad-Shafranov solver.  As traditionally done, the free-boundary
problem is reformulated in a bounded computational domain that
encloses all current-carrying plasmas. The coil current contribution
to the plasma equilibrium enters through the value of the magnetic
flux that is evaluated by Green's function methods.  The main focus
here is to solve the free boundary problem of the nonlinear
Grad-Shafranov equation in an irregular computational domain by
combining the cut-cell algorithm with a regular mesh.  A key
application of free-bounday Grad-Shafranov solvers is to find the
required coil currents for precision plasma positioning and shaping.
This is cast into a proper optimization problem that determines the
coil current with targeted plasma shape and for which a cut-cell
algorithm is described.  To accelerate Picard iterations commonly used
in previous works, we also propose an improvement based on Aitken's
acceleration. This is applied to both the inner and outer loops of the
algorithm.  The proposed algorithm is implemented in parallel under
the PETSc framework.  Strong scaling study of the free-boundary
Grad-Shafranov solver with cut-cell algorithm demonstrates its
parallel performance.

A series of numerical tests is presented to demonstrate the accuracy
of the cut-cell algorithm and the efficiency of the free-boundary
Grad-Shafranov solver.  In particular, a refinement study based on
manufactured solutions confirms a second-order accuracy of the
cut-cell implementation and tokamak examples, including that of ITER,
are presented to verify the effectiveness and efficiency fo the full
algorithm.  Moreover, numerical results demonstrate that the converged
magnetic flux keeps the predetermined plasma domain shape even with
different pressure profile in the source term of the Grad-Shafranov
equation both in the rectangular domain and the irregular
limiter-bounded domain.  It is also found that the Aitken's
acceleration we employed in the algorithm is more efficient than
Picard iteration with fixed under-relaxation. In conclusion, numerical
results shows a good performance of the free-boundary cut-cell
Grad-Shafranov solver.

\section*{Acknowledgments}
The authors wish to thank Dr.~Yueqiang Liu of General Atomics for
sharing the Grad-Shafranov equilibrium of ITER reference number
``ABT4ZL''.  This research used resources of the National Energy
Research Scientific Computing Center (NERSC), a U.S. Department of
Energy Office of Science User Facility operated.

\bibliographystyle{siamplain}
\bibliography{references}

\begin{thebibliography}{10}

\bibitem{ambrosetti1990remarks}
{\sc A.~Ambrosetti, M.~Calahorrano, and F.~Dobarro}, {\em Remarks on the
  grad-shafranov equation}, Applied Mathematics Letters, 3 (1990), pp.~9--11.

\bibitem{ariola-etal-ITER-shape-2000}
{\sc M.~Ariola, A.~Pironti, and A.~Portone}, {\em Vertical stabilization and
  plasma shape control in the iter-feat tokamak}, in Proceedings of the 2000
  IEEE International Conference on Control Applications, Los Alamitos, CA, USA,
  sep 2000, IEEE Computer Society, pp.~401,402,403,404,405.

\bibitem{balay2019petsc}
{\sc S.~Balay, S.~Abhyankar, M.~Adams, J.~Brown, P.~Brune, K.~Buschelman,
  L.~Dalcin, A.~Dener, V.~Eijkhout, W.~Gropp, et~al.}, {\em Petsc users
  manual},  (2019).

\bibitem{banks2017stable}
{\sc J.~W. Banks, W.~D. Henshaw, D.~W. Schwendeman, and Q.~Tang}, {\em A stable
  partitioned fsi algorithm for rigid bodies and incompressible flow. part ii:
  General formulation}, Journal of Computational Physics, 343 (2017),
  pp.~469--500.

\bibitem{banks2018stable}
{\sc J.~W. Banks, W.~D. Henshaw, D.~W. Schwendeman, and Q.~Tang}, {\em A stable
  partitioned fsi algorithm for rigid bodies and incompressible flow in three
  dimensions}, Journal of Computational Physics, 373 (2018), pp.~455--492.

\bibitem{berger2017cut}
{\sc M.~Berger}, {\em Cut cells: Meshes and solvers}, in Handbook of Numerical
  Analysis, vol.~18, Elsevier, 2017, pp.~1--22.

\bibitem{berrut1988rational}
{\sc J.-P. Berrut}, {\em Rational functions for guaranteed and experimentally
  well-conditioned global interpolation}, Computers \& Mathematics with
  Applications, 15 (1988), pp.~1--16.

\bibitem{berrut2005recent}
{\sc J.-P. Berrut, R.~Baltensperger, and H.~D. Mittelmann}, {\em Recent
  developments in barycentric rational interpolation}, Journal of Computational
  and Applied Mathematics, 259 (2005), pp.~95--107.

\bibitem{borazjani2008curvilinear}
{\sc I.~Borazjani, L.~Ge, and F.~Sotiropoulos}, {\em Curvilinear immersed
  boundary method for simulating fluid structure interaction with complex 3d
  rigid bodies}, Journal of Computational physics, 227 (2008), pp.~7587--7620.

\bibitem{brown1997overture}
{\sc D.~L. Brown, W.~D. Henshaw, and D.~J. Quinlan}, {\em Overture: An
  object-oriented framework for solving partial differential equations}, in
  International Conference on Computing in Object-Oriented Parallel
  Environments, Springer, 1997, pp.~177--184.

\bibitem{cerfon2010one}
{\sc A.~J. Cerfon and J.~P. Freidberg}, {\em {``One size fits all'' analytic
  solutions to the Grad--Shafranov equation}}, Physics of Plasmas, 17 (2010),
  p.~032502.

\bibitem{chen1997simple}
{\sc S.~Chen, B.~Merriman, S.~Osher, and P.~Smereka}, {\em A simple level set
  method for solving stefan problems}, Journal of Computational Physics, 135
  (1997), pp.~8--29.

\bibitem{creely-sparc-jpp-2020}
{\sc A.~J. Creely, M.~J. Greenwald, S.~B. Ballinger, D.~Brunner, J.~Canik,
  J.~Doody, T.~Fülöp, D.~T. Garnier, R.~Granetz, T.~K. Gray, and et~al.},
  {\em Overview of the sparc tokamak}, Journal of Plasma Physics, 86 (2020),
  p.~865860502.

\bibitem{devendran2014higher}
{\sc D.~Devendran, D.~Graves, and H.~Johansen}, {\em A higher-order
  finite-volume discretization method for poisson's equation in cut cell
  geometries}, arXiv preprint arXiv:1411.4283,  (2014).

\bibitem{devendran2017fourth}
{\sc D.~Devendran, D.~Graves, H.~Johansen, and T.~Ligocki}, {\em A fourth-order
  cartesian grid embedded boundary method for poisson’s equation},
  Communications in Applied Mathematics and Computational Science, 12 (2017),
  pp.~51--79.

\bibitem{faugeras2017fem}
{\sc B.~Faugeras and H.~Heumann}, {\em Fem-bem coupling methods for tokamak
  plasma axisymmetric free-boundary equilibrium computations in unbounded
  domains}, Journal of Computational Physics, 343 (2017), pp.~201--216.

\bibitem{floater2007barycentric}
{\sc M.~S. Floater and K.~Hormann}, {\em Barycentric rational interpolation
  with no poles and high rates of approximation}, Numerische Mathematik, 107
  (2007), pp.~315--331.

\bibitem{gibou2005fourth}
{\sc F.~Gibou and R.~Fedkiw}, {\em A fourth order accurate discretization for
  the laplace and heat equations on arbitrary domains, with applications to the
  stefan problem}, Journal of Computational Physics, 202 (2005), pp.~577--601.

\bibitem{gibou2002second}
{\sc F.~Gibou, R.~P. Fedkiw, L.-T. Cheng, and M.~Kang}, {\em A
  second-order-accurate symmetric discretization of the poisson equation on
  irregular domains}, Journal of Computational Physics, 176 (2002),
  pp.~205--227.

\bibitem{golub1979generalized}
{\sc G.~H. Golub, M.~Heath, and G.~Wahba}, {\em Generalized cross-validation as
  a method for choosing a good ridge parameter}, Technometrics, 21 (1979),
  pp.~215--223.

\bibitem{heumann2015quasi}
{\sc H.~Heumann, J.~Blum, C.~Boulbe, B.~Faugeras, G.~Selig, P.~Hertout,
  E.~Nardon, J.-M. An{\'e}, S.~Br{\'e}mond, and V.~Grandgirard}, {\em
  Quasi-static free-boundary equilibrium of toroidal plasma with cedres++:
  Computational methods and applications}, Journal of Plasma Physics,  (2015),
  p.~35.

\bibitem{heumann2017finite}
{\sc H.~Heumann and F.~Rapetti}, {\em A finite element method with overlapping
  meshes for free-boundary axisymmetric plasma equilibria in realistic
  geometries}, Journal of Computational Physics, 334 (2017), pp.~522--540.

\bibitem{honda2010simulation}
{\sc M.~Honda}, {\em Simulation technique of free-boundary equilibrium
  evolution in plasma ramp-up phase}, Computer Physics Communications, 181
  (2010), pp.~1490--1500.

\bibitem{howell2014solving}
{\sc E.~Howell and C.~R. Sovinec}, {\em Solving the grad--shafranov equation
  with spectral elements}, Computer Physics Communications, 185 (2014),
  pp.~1415--1421.

\bibitem{huysmans1991isoparametric}
{\sc G.~Huysmans, J.~Goedbloed, W.~Kerner, et~al.}, {\em Isoparametric bicubic
  hermite elements for solution of the grad-shafranov equation}, International
  Journal of Modern Physics C, 2 (1991), pp.~371--376.

\bibitem{jardin2010computational}
{\sc S.~Jardin}, {\em Computational methods in plasma physics}, CRC Press,
  2010.

\bibitem{jardin1986dynamic}
{\sc S.~C. Jardin, N.~Pomphrey, and J.~Delucia}, {\em Dynamic modeling of
  transport and positional control of tokamaks}, Journal of computational
  Physics, 66 (1986), pp.~481--507.

\bibitem{jeon2015development}
{\sc Y.~M. Jeon}, {\em Development of a free-boundary tokamak equilibrium
  solver for advanced study of tokamak equilibria}, Journal of the Korean
  Physical Society, 67 (2015), pp.~843--853.

\bibitem{johansen1998cartesian}
{\sc H.~Johansen and P.~Colella}, {\em A cartesian grid embedded boundary
  method for poisson's equation on irregular domains}, Journal of Computational
  Physics, 147 (1998), pp.~60--85.

\bibitem{johnson1979numerical}
{\sc J.~L. Johnson, H.~Dalhed, J.~Greene, R.~Grimm, Y.~Hsieh, S.~Jardin,
  J.~Manickam, M.~Okabayashi, R.~Storer, A.~Todd, et~al.}, {\em Numerical
  determination of axisymmetric toroidal magnetohydrodynamic equilibria},
  Journal of Computational Physics, 32 (1979), pp.~212--234.

\bibitem{jomaa2005embedded}
{\sc Z.~Jomaa and C.~Macaskill}, {\em The embedded finite difference method for
  the poisson equation in a domain with an irregular boundary and dirichlet
  boundary conditions}, Journal of Computational Physics, 202 (2005),
  pp.~488--506.

\bibitem{lackner1976computation}
{\sc K.~Lackner}, {\em Computation of ideal mhd equilibria}, Computer Physics
  Communications, 12 (1976), pp.~33--44.

\bibitem{lee2015ecom}
{\sc J.~Lee and A.~Cerfon}, {\em Ecom: A fast and accurate solver for toroidal
  axisymmetric mhd equilibria}, Computer Physics Communications, 190 (2015),
  pp.~72--88.

\bibitem{leveque1992numerical}
{\sc R.~J. LeVeque}, {\em Numerical methods for conservation laws}, vol.~3,
  Springer, 1992.

\bibitem{leveque1994immersed}
{\sc R.~J. Leveque and Z.~Li}, {\em The immersed interface method for elliptic
  equations with discontinuous coefficients and singular sources}, SIAM Journal
  on Numerical Analysis, 31 (1994), pp.~1019--1044.

\bibitem{li2020solving}
{\sc H.~Li and P.~Zhu}, {\em Solving the grad-shafranov equation using spectral
  elements for tokamak equilibrium with toroidal rotation}, Computer Physics
  Communications,  (2020), p.~107264.

\bibitem{liu2020numerical}
{\sc S.~Liu, Y.~Du, and X.~Liu}, {\em Numerical studies of a class of
  reaction--diffusion equations with stefan conditions}, International Journal
  of Computer Mathematics, 97 (2020), pp.~959--979.

\bibitem{liu2018numerical}
{\sc S.~Liu and X.~Liu}, {\em Numerical methods for a two-species
  competition-diffusion model with free boundaries}, Mathematics, 6 (2018),
  p.~72.

\bibitem{Liu-etal-NF-2015}
{\sc Y.~Liu, R.~Akers, I.~Chapman, Y.~Gribov, G.~Hao, G.~Huijsmans, A.~Kirk,
  A.~Loarte, S.~Pinches, M.~Reinke, et~al.}, {\em Modelling toroidal rotation
  damping in iter due to external 3d fields}, Nuclear Fusion, 55 (2015),
  p.~063027.

\bibitem{lutjens1996chease}
{\sc H.~L{\"u}tjens, A.~Bondeson, and O.~Sauter}, {\em The chease code for
  toroidal mhd equilibria}, Computer physics communications, 97 (1996),
  pp.~219--260.

\bibitem{mok2001accelerated}
{\sc D.~Mok, W.~Wall, and E.~Ramm}, {\em Accelerated iterative substructuring
  schemes for instationary fluid-structure interaction}, Computational fluid
  and solid mechanics, 2 (2001), pp.~1325--1328.

\bibitem{pataki2013fast}
{\sc A.~Pataki, A.~J. Cerfon, J.~P. Freidberg, L.~Greengard, and M.~O’Neil},
  {\em A fast, high-order solver for the grad--shafranov equation}, Journal of
  Computational Physics, 243 (2013), pp.~28--45.

\bibitem{dpgGS}
{\sc Z.~Peng, Q.~Tang, and X.-Z. Tang}, {\em An adaptive discontinuous
  petrov--galerkin method for the grad--shafranov equation}, SIAM Journal on
  Scientific Computing, 42 (2020), pp.~B1227--B1249.

\bibitem{polyanin2006handbook}
{\sc A.~D. Polyanin and A.~V. Manzhirov}, {\em Handbook of mathematics for
  engineers and scientists}, CRC Press, 2006.

\bibitem{sanchez2019hybridizable}
{\sc T.~S{\'a}nchez-Vizuet and M.~E. Solano}, {\em A hybridizable discontinuous
  galerkin solver for the grad--shafranov equation}, Computer Physics
  Communications, 235 (2019), pp.~120--132.

\bibitem{sanchez2020adaptive}
{\sc T.~S{\'a}nchez-Vizuet, M.~E. Solano, and A.~J. Cerfon}, {\em Adaptive
  hybridizable discontinuous galerkin discretization of the grad-shafranov
  equation by extension from polygonal subdomains}, Computer Physics
  Communications,  (2020), p.~107239.

\bibitem{schneider1986some}
{\sc C.~Schneider and W.~Werner}, {\em Some new aspects of rational
  interpolation}, Mathematics of computation, 47 (1986), pp.~285--299.

\bibitem{schwartz2006cartesian}
{\sc P.~Schwartz, M.~Barad, P.~Colella, and T.~Ligocki}, {\em A cartesian grid
  embedded boundary method for the heat equation and poisson’s equation in
  three dimensions}, Journal of Computational Physics, 211 (2006),
  pp.~531--550.

\bibitem{hagenow1975}
{\sc K.~von Hagenow and K.~Lackner}, eds., {\em Proceedings of the 7th Conf. on
  the Numerical Simulation of Plasmas}, 1975.

\end{thebibliography}
\end{document}